\documentclass{ectj}

\usepackage{amsfonts,amssymb,graphics,epsfig,verbatim,bm,latexsym,amsmath,url,amsbsy}

\usepackage{booktabs}
\usepackage{graphicx}

\DeclareMathOperator{\Cov}{Cov}

\DeclareMathOperator{\cum}{cum}
\DeclareMathOperator{\Var}{Var}

\newcommand{\ee}{\ensuremath{\mathrm{e}}}
\newcommand{\ii}{\ensuremath{\mathrm{i}}}

\newcommand{\IC}{\ensuremath{\mathbb{C}}}

\newcommand{\E}{\ensuremath{\mathbb{E}}}

\newcommand{\IH}{\ensuremath{\mathbb{H}}}

\newcommand{\IL}{\ensuremath{\mathbb{L}}}
\newcommand{\IN}{\ensuremath{\mathbb{N}}}
\newcommand{\IP}{\ensuremath{\mathbb{P}}}
\newcommand{\IR}{\ensuremath{\mathbb{R}}}
\newcommand{\IZ}{\ensuremath{\mathbb{Z}}}

\newcommand{\bfrakB}{\ensuremath{\bm{\mathfrak{B}}}}
\newcommand{\bfrakf}{\ensuremath{\bm{\mathfrak{f}}}}
\newcommand{\bfrakR}{\ensuremath{\bm{\mathfrak{R}}}}

\newcommand{\bG}{\ensuremath{\bm{G}}}
\newcommand{\bI}{\ensuremath{\bm{I}}}

\newcommand{\bX}{\ensuremath{\bm{X}}}

\newcommand{\ff}[2]{\mathfrak{f}_{#1, #2}}

\newcommand{\BB}[2]{\bm{B}_{#1, #2}}
\newcommand{\HH}[2]{\IH_{#1, #2}}

\newcommand{\pkg}[1]{{\normalfont\fontseries{b}\selectfont #1}}

\renewcommand{\emph}[1]{#1}

\newtheorem{theorem}{Theorem}
\newtheorem{assumption}{Assumption}

\renewcommand{\thesection}{\arabic{section}}
\renewcommand{\theequation}{\arabic{section}.\arabic{equation}}

  \received{October 2018}
  \volume{20}
\title[Quantile Coherency]{Quantile Coherency: A General Measure for Dependence between Cyclical Economic Variables}

\author[J.~Barun\'{i}k and T.~Kley]{Jozef~Barun\'{i}k$^{\dagger}$ and
                       Tobias~Kley$^{\ddagger}$}

\address{$^{\dagger}$Econometric Department, IITA, The Czech Academy of Sciences\\ and Institute of Economic Studies,  Charles University in Prague.}
\email{barunik@fsv.cuni.cz}

\address{$^{\ddagger}$School of Mathematics, University of Bristol.}
\email{tobias.kley@bristol.ac.uk}

\def\AmSTeX{$\cal A$\kern-.1667em\lower.5ex\hbox{$\cal M$}\kern-.125em
            $\cal S$-\TeX}
\def\BibTeX{{\rm B\kern-.05em{\sc i\kern-.025em b}\kern-.08em
            T\kern-.1667em\lower.7ex\hbox{E}\kern-.125emX}}

\begin{document}

  \begin{abstract}
    In this paper, we introduce quantile coherency to measure general dependence structures emerging in the joint distribution in the frequency domain and argue that this type of dependence is natural for economic time series but remains invisible when only the traditional analysis is employed.
		We define estimators which capture the general dependence structure, provide a detailed analysis of their asymptotic properties and discuss how to conduct inference for a general class of possibly nonlinear processes. In an empirical illustration we examine the dependence of bivariate stock market returns and shed new light on measurement of tail risk in financial markets. We also provide a modelling exercise to illustrate how applied researchers can benefit from using quantile coherency when assessing time series models.


  \end{abstract}

\section{Dependence structures beyond second-order moments}

One of the fundamental problems faced by a researcher in economics is how to quantify the dependence between economic variables. Although correlated variables are rather commonly observed phenomena in economics, it is often the case that strongly correlated variables under study are truly independent, and what we measure is mere spurious correlation; see \cite{granger1974spurious}. Conversely, but equally deluding, \emph{uncorrelated} variables may possess dependence in different parts of the joint distribution, and/or at different frequencies. This dependence stays hidden when classical measures based on linear correlation and traditional cross-spectral analysis are used; see \cite{croux2001measure}, \cite{ning2009dependence} and \cite{fan2014copulas}. Hence, conventional models derived from averaged quantities as for example covariance-based measures may deliver rather misleading results.

In this paper, we introduce a new class of cross-spectral densities that characterise the dependence in quantiles of the joint distribution across frequencies (i.\,e., with respect to cycles). Subsequently, standardisation of the before-mentioned quantile spectra yields a related quantity to which we will refer to as quantile coherency. We define and motivate the quantile-based cross-spectral quantities in analogy to their traditional counterparts. Yet, instead of quantifying dependence in terms of joint moments (i.\,e., by averaging with respect to the joint distribution), the new measures are defined in terms of the probabilities to exceed quantiles. Hence, they are designed to detect \emph{any} general type of dependence structure that may arise between variables under study.

Such complex dynamics may arise naturally in many macroeconomic, or financial time series such as growth rates, inflation, housing markets, or stock market returns. In financial markets, extremely scarce and negative events in one asset can cause irrational outcomes and panics leading investors to ignore economic fundamentals and cause similarly extreme negative outcomes in other assets. In such situations, markets may be connected more strongly than in calm periods of small, or positive returns; cf. \cite{bae2003new}. Hence, the co-occurrences of large negative values may be more common across stock markets than co-occurrences of large positive values reflecting asymmetric behaviour of economic agents. Moreover, long-term fluctuations in quantiles of the joint distribution may differ from the ones in the short-term due to differing risk perception of economic agents over distinct investment horizons. This behaviour produces various degrees of persistence at different parts of the joint distribution, while on average the stock market returns remain impersistent. In univariate macroeconomic variables, researchers document asymmetric adjustment paths 
(cf. \cite{neftci1984economic} and \cite{enders1998unit}) as firms are more prone to an increase than to a decrease in prices. Asymmetric business cycle dynamics at different quantiles can be caused by positive shocks to output being more persistent than negative shocks. While output fluctuations are known to be persistent, \cite{beaudry1993recessions} document less persistence at longer horizons. Such asymmetric dependence at different horizons can be shared by multiple variables. Because classical, covariance-based approaches only take averaged information into account, these types of dependence fail to be identified by traditional means. Revealing such dependence structures, quantile cross-spectral analysis introduced in this paper can fundamentally change the way how we view the dependence between economic time series, and opens new possibilities for the modelling of interactions between economic and financial variables.

Quantile cross-spectral analysis provides a general, unifying framework for estimating dependence between economic time series. As noted in the early work of \cite{granger1966typical}, the spectral distribution of an economic variable has a typical shape which distinguishes long-term fluctuations from short-term ones. These fluctuations point to economic activity at different frequencies (after removal of trend in mean, as well as seasonal components). After \cite{granger1966typical} studied the behaviour of single time series, important literature using cross-spectral analysis to identify the dependence between variables quickly emerged (from \cite{granger1969investigating} to more recent \cite{croux2001measure}). Instead of considering only cross-sectional correlations, researchers started to use coherency (frequency dependent correlation) to investigate short-run and long-run dynamic properties of multiple time series, and identify business cycle synchronisation; see \cite{croux2001measure}. In one of his very last papers, \cite{granger2010some} hypothesised about possible cointegrating relationships in quantiles, leading to the first notion of general types of dependence that quantile cross-spectral analysis is addressing. The quantile cointegration developed by \cite{xiao2009quantile} partially addresses the problem, but does not allow to fully explore the frequency dependent structure of correlations in different quantiles of the joint distribution.

\begin{figure}[t!]
  \begin{center}
        \includegraphics[width=0.32\linewidth]{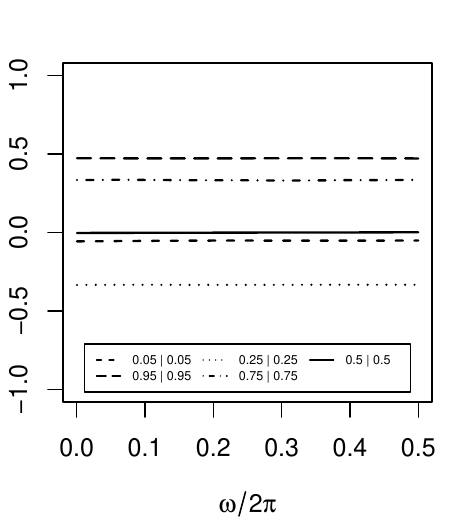} 
        \includegraphics[width=0.32\linewidth]{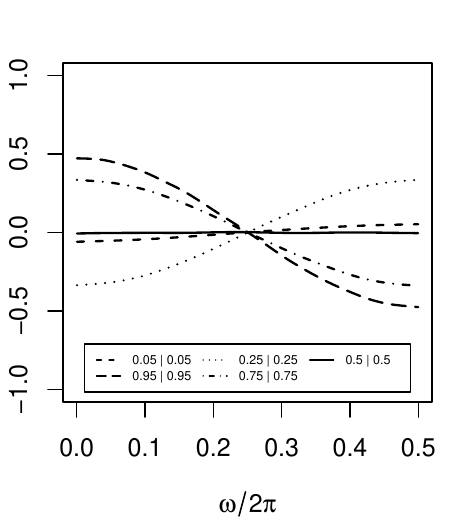}
        \includegraphics[width=0.32\linewidth]{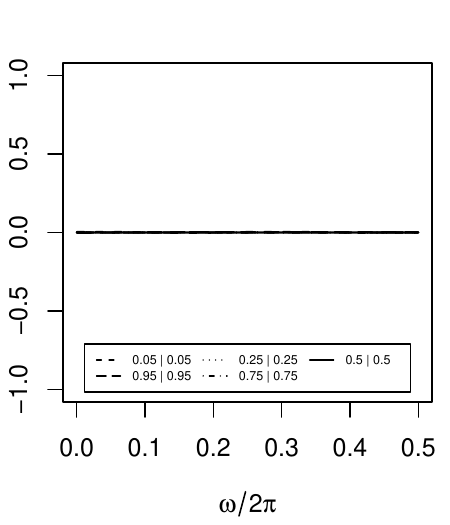} 
	\end{center}
  \vspace*{-0.5cm}\caption{\label{fig:simFreq}Illustration of dependence between processes $x_t$ and $y_t$. 
	}
\end{figure}

\begin{figure}[t!]
  \begin{center}
        \includegraphics[width=0.32\linewidth]{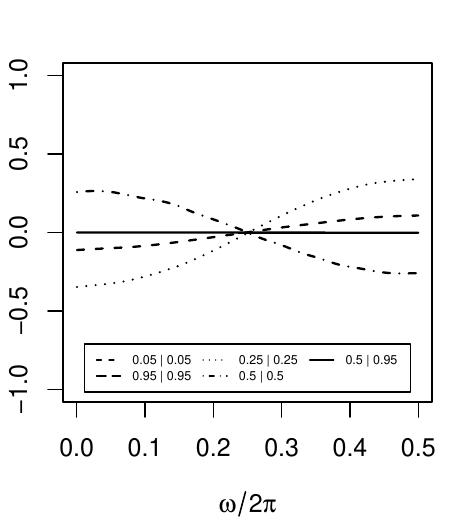} 
        \includegraphics[width=0.32\linewidth]{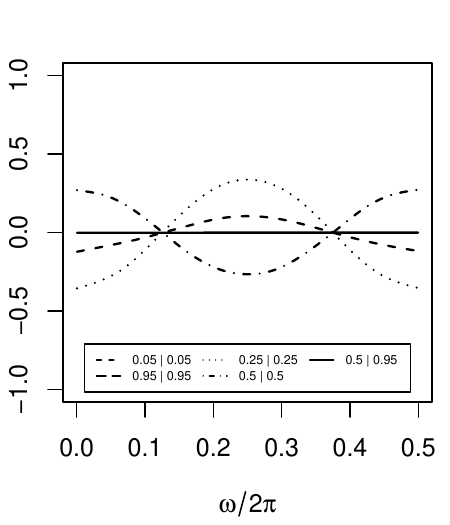}
        \includegraphics[width=0.32\linewidth]{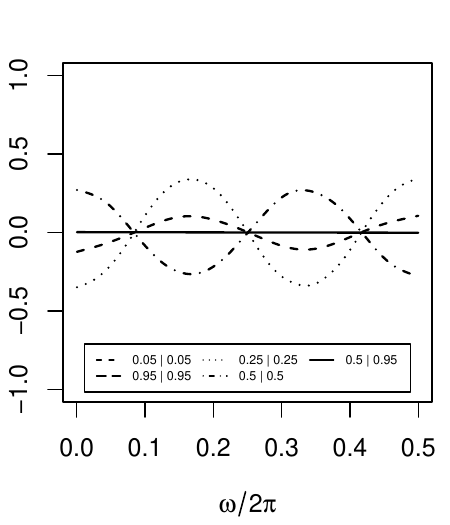} 
	\end{center}
  \vspace*{-0.5cm}\caption{\label{fig:simFreqQVAR}Illustration of dependence between vector quantile autoregressive processes. 
	}
\end{figure}

Three toy examples illustrating the potential offered by quantile cross-spectral analysis are depicted in Figure~\ref{fig:simFreq}. In each example one distinct type of dependence is considered: cross-sectional dependence (left), serial dependence (centre), and independence (right). We consider bivariate processes $(x_t, y_t)$ that possess the desired dependence structure, but are indistinguishable in terms of traditional coherency. In the examples, $(\epsilon_t)$ is an independent sequence of standard normally distributed random variables. In the left column of Figure~\ref{fig:simFreq} the dependence emerging between $\epsilon_t$ and $\epsilon_t^2$ is depicted. It is important to observe that $\epsilon_t$ and $\epsilon_s^2$ are uncorrelated.
Therefore, traditional coherency for $(\epsilon_t, \epsilon_t^2)$ would read zero across all frequencies, even though it is obvious that $\epsilon_t$ and $\epsilon_t^2$ are dependent. From the newly introduced quantile coherency, this dependence can easily be observed. More precisely, we can distinguish various degrees of dependence for each part of the distribution. For example, there is no dependence in the centre of the distribution (i.\,e., $0.5 | 0.5$), but when the quantile levels are different from $0.5$ the dependence becomes visible.\footnote{All plots show real parts of the complex-valued quantities for illustratory purposes. Further discussion on how to interpret the real part and the imaginary part of quantile coherency are deferred to Section~\ref{sec:ccsdk:estim}.} In this example the quantile coherency is constant across frequencies, which corresponds to the fact that there is no serial dependence. In the centre column of Figure~\ref{fig:simFreq} the process $(\epsilon_t, \epsilon_{t-1}^2)$ is studied, where we have introduced a time lag. Intuitively, the dependence in quantiles of this bivariate process will be the same as in the previous example (left column) in the long run, referring to frequencies close to zero. With increasing frequency, dependence will decline or incline gradually to values with opposite signs, as high frequency movements are in opposition due to the lag shift. This is clearly captured by quantile coherency, while the dependence structure would stay hidden away from traditional coherency, again, as it averages the dependence across quantiles. We can think about these processes as being ``spuriously independent''. To demonstrate the behaviour of the quantile coherency when the processes under consideration are \emph{truly independent}, we observe in the right column of Figure~\ref{fig:simFreq} the quantities for independent bivariate Gaussian white noise, where quantile coherency displays zero dependence at all quantiles and frequencies, as expected. These illustrations strongly support our claim that there is need for more general measures that can provide a better understanding of the dependence between variables. These very simple, yet illuminating motivating examples focus on uncovering dependence in uncorrelated variables. Later in the text (Section~\ref{sec:modexerc}), we further discuss a data generating process based on quantile vector autoregression (QVAR), which is able to generate even richer dependence structures, revealing once more the limitations of the traditional approach. In Figure~\ref{fig:simFreqQVAR}, the real part of the quantile coherencies of the QVAR(1), QVAR(2) and QVAR(3) example processes are shown. Further, in Section~\ref{sec:interpretationgauss}, we discuss how to interpret quantile coherency in the special cases of bivariate Gaussian VAR(1).

This paper is organised as follows. In Section~\ref{sec:ccsdk} we introduce notation, define quantile coherency and an estimator for it. In Section~\ref{sec:ccsdk:estim} we discuss the proposed methodology and related literature. In Section~\ref{sec:asymp} we provide a rigorous asymptotic analysis of the estimator's statistical properties. In Section~\ref{sec:emp}, to support our theoretical discussions empirically, we employ the new methodology to inspect bivariate stock market returns, one of the most prominent time series in economics, and reveal dependencies in cycles of quantile-based features. We continue our empirical study in Section~\ref{sec:modexerc} by using quantile coherency to compare time series models with respect to their capabilities to capture the revealed dependencies. In the supplementary material to this paper (available from the publisher's homepage), we discuss additional quantile-based cross-spectral quantities (Section~\ref{sec:additionalquantities}), discuss quantile vector autoregressive processes as examples with rich dynamics (Section~\ref{sec:vqar}), discuss how the new, quantile-based spectral quantities and their traditional counterparts are related (Section~\ref{sec:interpretationgauss}), state additional theoretical results (Section~\ref{sec:asymp2}), comment on the construction of the interval estimators (Section~\ref{sec:CI}), and provide rigorous proofs for all theoretical results (Section~\ref{sec:proofs}).

\section{Quantile cross-spectral quantities and their estimators}
\label{sec:ccsdk}
Throughout the paper $(\bX_t)_{t \in \IZ}$ denotes a $d$-variate, strictly stationary process, with components $X_{t,j}$, $j=1,\ldots,d$; i.\,e. $\bX_t = (X_{t,1}, \ldots, X_{t,d})'$. The marginal distribution function of $X_{t,j}$ will be denoted by $F_j$, and by $q_j(\tau) := F_j^{-1}(\tau) := \inf\{ q \in \IR : \tau \leq F_j(q)\}$, where $\tau \in [0,1]$, we denote the corresponding quantile function. We use the convention $\inf \emptyset = +\infty$, such that, if $\tau = 0$ or $\tau = 1$, then $-\infty$ and $+\infty$ are possible values for $q_j(\tau)$, respectively. We will write $\overline{z}$ for the complex conjugate, $\Re z$ for the real part and $\Im z$ for the imaginary part of $z \in \IC$, respectively. The transpose of a matrix $\bm{A}$ will be denoted by~$\bm{A}'$, the inverse of a regular matrix $\bm{B}$ will be denoted by~$\bm{B}^{-1}$.

As a measure for the serial and cross-dependency structure of $(\bX_t)_{t \in \IZ}$, we define the \emph{matrix of quantile cross-covariance kernels}, $\Gamma_k(\tau_1, \tau_2) := (\gamma_k^{j_1, j_2}(\tau_1, \tau_2))_{j_1, j_2 = 1,\ldots,d}$, where
\begin{equation}
\label{eqn:ccck}
	\gamma_k^{j_1, j_2}(\tau_1, \tau_2) := \Cov \Big( I\{X_{t+k, j_1} \leq q_{j_1}(\tau_1)\}, I\{X_{t, j_2} \leq q_{j_2}(\tau_2)\} \Big),
\end{equation}
$j_1, j_2 \in \{1,\ldots,d\}, \, k \in \IZ, \, \tau_1, \tau_2 \in [0,1]$, and $I\{A\}$ denotes the indicator function of the event~$A$.
In the frequency domain this yields (under appropriate mixing conditions) the \emph{matrix of quantile cross-spectral density kernels}
	$\bfrakf(\omega; \tau_1, \tau_2) := ( \mathfrak{f}^{j_1, j_2}(\omega; \tau_1, \tau_2) )_{j_1, j_2 = 1,\ldots, d}$,
where
\begin{equation}
\label{eqn:ccsdk}
	\mathfrak{f}^{j_1, j_2}(\omega; \tau_1, \tau_2) := (2\pi)^{-1} \sum_{k=-\infty}^{\infty} \gamma_k^{j_1, j_2}(\tau_1, \tau_2) \ee^{-\ii k \omega},
\end{equation}
$j_1, j_2 \in \{1,\ldots,d\}, \, \omega \in \IR, \, \tau_1, \tau_2 \in [0,1]$. 
A closely related quantity that can be used as a measure for the dynamic dependence of the two processes $(X_{t,j_1})_{t \in \IZ}$ and $(X_{t,j_2})_{t \in \IZ}$ is the \emph{quantile coherency kernel} of $(X_{t,j_1})_{t \in \IZ}$ and $(X_{t,j_2})_{t \in \IZ}$, which we define as
\begin{equation}
\label{eqn:qcoherency}
\mathfrak{R}^{j_1, j_2}(\omega; \tau_1, \tau_2)
:= \frac{\mathfrak{f}^{j_1, j_2}(\omega; \tau_1, \tau_2)}{\Big(\mathfrak{f}^{j_1, j_1}(\omega; \tau_1, \tau_1) \mathfrak{f}^{j_2, j_2}(\omega; \tau_2, \tau_2) \Big)^{1/2}},
\end{equation}
$(\tau_1, \tau_2) \in (0,1)^2$.
We define the estimator for the quantile cross-spectral density as the collection
\begin{equation}
	\label{eqn:inr}
I_{n,R}^{j_1, j_2}(\omega; \tau_1, \tau_2) := \frac{1}{2\pi n}
  d^{j_1}_{n,R}(\omega; \tau_1) d^{j_2}_{n,R}(-\omega; \tau_2),
\end{equation}
$j_1, j_2 = 1,\ldots, d$, $\omega\in  \IR$, $(\tau_1,\tau_2) \in [0,1]^2$, 
and call it the \emph{rank-based copula cross-periodograms}, shortly, the CCR-periodograms, where
\[d_{n,R}^{j}(\omega; \tau) := \sum_{t=0}^{n-1} I\{\hat F_{n,j}(X_{t,j}) \leq \tau\} \ee^{- \ii \omega t} = \sum_{t=0}^{n-1} I\{ R_{n;t,j}  \leq  n\tau\} \ee^{- \ii \omega t},\]
$j = 1,\ldots, d$, $\omega\in  \IR$, $\tau \in [0,1]$, and $\hat F_{n,j}(x) := n^{-1} \sum_{t=0}^{n-1} I\{X_{t,j} \leq x\}$ denotes the empirical distribution function of $X_{t,j}$ and $R_{n;t,j}$ denotes the (maximum) rank of $X_{t,j}$ among $X_{0,j}, \ldots, X_{n-1,j}$.
We will denote the \emph{matrix of CCR-periodograms} by
\begin{equation}
	\label{eqn:inrMatr}
	\bI_{n,R}(\omega; \tau_1, \tau_2) := ( I^{j_1, j_2}_{n, R}(\omega; \tau_1, \tau_2) )_{j_1, j_2 = 1,\ldots, d}.
\end{equation}
From the univariate case it is already known (cf.~Proposition~3.4 in \cite{kley2014}) that the CCR-periodograms fail to estimate $\mathfrak{f}^{j_1, j_2}(\omega; \tau_1, \tau_2)$ consistently. Consistency can be achieved by smoothing $I_{n,R}^{j_1, j_2}(\omega; \tau_1, \tau_2)$ across frequencies. More precisely, we consider
\begin{equation}
	\label{eqn:DefRankEstimator}
	\hat G^{j_1, j_2}_{n,R}(\omega; \tau_1, \tau_2) := \frac{2\pi}{n} \sum_{s=1}^{n-1} W_n\big( \omega - 2\pi s / n \big) I_{n,R}^{j_1, j_2}(2 \pi s / n, \tau_1, \tau_2),
\end{equation}
where $W_n$ denotes a sequence of weight functions, precisely to be defined in Section~\ref{sec:asymp}.

We will denote the \emph{matrix of smoothed CCR-periodograms} by
\begin{equation}
	\label{eqn:DefRankEstimatorMatr}
	\hat\bG_{n,R}(\omega; \tau_1, \tau_2) := ( \hat G^{j_1, j_2}_{n, R}(\omega; \tau_1, \tau_2) )_{j_1, j_2 = 1,\ldots, d}.
\end{equation}
The estimators for the quantile coherency is then given by
\begin{equation}\label{def:Rhat}
	\hat{\mathfrak{R}}^{j_1, j_2}_{n,R}(\omega; \tau_1, \tau_2)
:= \frac{\hat G^{j_1, j_2}_{n,R}(\omega; \tau_1, \tau_2)}{\Big(\hat G^{j_1, j_1}_{n,R}(\omega; \tau_1, \tau_1) \hat G^{j_2, j_2}_{n,R}(\omega; \tau_2, \tau_2) \Big)^{1/2}}.
\end{equation}
In Section~\ref{sec:asymp} we will prove that
\[\hat\bfrakR_{n,R}(\omega;\tau_1,\tau_2) := \big( \hat{\mathfrak{R}}^{j_1, j_2}_{n,R}(\omega; \tau_1, \tau_2) \big)_{j_1, j_2 = 1,\ldots,d}\] is a legitimate estimator for
$\bfrakR(\omega;\tau_1,\tau_2) := \big( \mathfrak{R}^{j_1, j_2}(\omega; \tau_1, \tau_2) \big)_{j_1, j_2 = 1,\ldots,d}$, the \emph{matrix of quantile coherencies}.

\section{Discussion of the introduced quantities and estimators}
\label{sec:ccsdk:estim}

The quantile-based quantities defined in Section~\ref{sec:ccsdk} are functions of the two variables $\tau_1$ and $\tau_2$. They are thus richer in information than the traditional counterparts. We have added the term \emph{kernel} to the name for the quantities to stress this fact, but will frequently omit it in the rest of the paper, for the sake of brevity.
For continuous $F_{j_1}$ and $F_{j_2}$, the quantile cross-covariances defined in~\eqref{eqn:ccck} coincide with the difference of the copula of $(X_{t+k,j_1}, X_{t,j_2})$ and the independence copula. Thus, they provide important information about both the serial dependence (by letting $k$ vary) and the cross-section-dependence (by choosing $j_1 \neq j_2$).
For the quantile cross-spectral density we have
\begin{equation}
\label{eqn:relCop}
  \int_{-\pi}^{\pi} \mathfrak{f}^{j_1, j_2}(\omega; \tau_1, \tau_2) \ee^{\ii k \omega} {\rm d}\omega + \tau_1 \tau_2 = \IP \Big( X_{t+k,j_1} \leq q_{j_1}(\tau_1), X_{t,j_2} \leq q_{j_2}(\tau_2) \Big),
\end{equation}
where the quantity on the right hand side, as a function of $(\tau_1, \tau_2)$, is again the copula of the pair $(X_{t+k,j_1}, X_{t,j_2})$. The equality~\eqref{eqn:relCop} thus shows how any of the pair copulas can be derived from the quantile cross-spectral density kernel defined in~\eqref{eqn:ccsdk}. Thus, the quantile cross-spectral density kernel provides a full description of all copulas of pairs in the process. Comparing these new quantities with their traditional counterparts, it can be observed that covariances and means are essentially replaced by copulas and quantiles. Similar to the regression setting, where this approach provides valuable extra information (cf. \cite{Koenker2005}), the quantile-based approach to spectral analysis supplements the traditional $L^2$-spectral analysis.

Observe that $\bfrakR$ takes values in $\IC^{d \times d}$ (the set of all complex-valued $d \times d$ matrices). Further, note that, as a function of $\omega$, but for fixed $\tau_1, \tau_2$, it coincides with the traditional coherency of the bivariate, binary process
\begin{equation}
\label{eqn:clipped}
	\Big( I\{X_{t,j_1} \leq q_{j_1}(\tau_1)\}, I\{X_{t,j_2} \leq q_{j_2}(\tau_2)\} \Big)_{t \in \IZ}.
\end{equation}
The time series in~\eqref{eqn:clipped} has the bivariate time series $(X_{t,j_1}, X_{t,j_2})_{t \in \IZ}$ as a ``latent driver'' and indicates whether the values of the components $j_1$ and $j_2$ are below the respective marginal distribution's $\tau_1$- and $\tau_2$-quantile.

Note the important fact that $\mathfrak{R}^{j_1, j_2}(\omega; \tau_1, \tau_2)$ is undefined when $(\tau_1, \tau_2)$ is on the boundary of $[0,1]^2$. By Cauchy-Schwarz inequality, we further observe that the range of possible values is limited to $\mathfrak{R}^{j_1, j_2}(\omega; \tau_1, \tau_2) \in \{ z \in \IC : |z| \leq 1\}$. Note that, as $(\tau_1, \tau_2)$ approaches the border of the unit square, the quantile cross-spectral density vanishes. Therefore, quantile coherency is better suited to measure dependence of extremes than the quantile cross-spectral density (which is not standardised). Implicitly, we take advantage of the fact that the quantile cross-spectral density and quantile spectral densities vanish at the same rate and therefore the quotient yields a meaningful quantity when the quantile levels $(\tau_1, \tau_2)$ approaches the border of the unit square.

The quantile coherency kernel contains very valuable information about the joint dynamics of the time series $(X_{t,j_1})_{t \in \IZ}$ and $(X_{t,j_2})_{t \in \IZ}$. In contrast to the traditional case, where coherency will always equal one if $j_1 = j_2 =: j$, the quantile-based versions of these quantities are capable of delivering valuable information about one single component of~$(\bX_t)_{t \in \IZ}$ as well. Quantile coherency then quantifies the joint dynamics of $(I\{X_{t,j} \leq q_{j}(\tau_1)\})_{t \in \IZ}$ and $(I\{X_{t,j} \leq q_{j}(\tau_2)\})_{t \in \IZ}$.

Note that quantile coherency is a complex-valued, $2\pi$-periodic function of the variable $\omega$, and Hermitian in the sense that we have
\[\overline{\mathfrak{R}^{j_1, j_2}(\omega; \tau_1, \tau_2)} = \mathfrak{R}^{j_1, j_2}(-\omega; \tau_1, \tau_2) = \mathfrak{R}^{j_2, j_1}(\omega; \tau_2, \tau_1) = \mathfrak{R}^{j_2, j_1}(2\pi + \omega; \tau_2, \tau_1).\]

Following similar arguments as in Section~2.1 of \cite{BirrEtAl2018}, it can be shown that $\Re \mathfrak{R}^{j_1, j_2}(\omega; \tau_1, \tau_2)$ describes the dynamics of the process switching between the $j_1$st component being below the $\tau_1$-quantile and the $j_2$nd component being above the $\tau_2$-quantile. Consequently, for $\tau_1$ close to 0 and for $\tau_2$ close to 1 it describes the dynamics of changing from an extreme in one component to an extreme in another component. Further, it can be shown that $\Im \mathfrak{R}^{j_1, j_2}(\omega; \tau_1, \tau_2)$ contains information about asymmetry.

A discussion of related quantities, how to interpret, how not to interpret them and how they are related to their traditional counterparts in the Gaussian case can be found in Sections~\ref{sec:additionalquantities}, \ref{sec:vqar}, and~\ref{sec:interpretationgauss} of the supplementary material.

Recently, important contributions that aim at accounting for more general dynamics emerged in the literature. Measures as, for example, distance correlation~\citep{SzekelyEtAl2007} and martingale difference correlation~\citep{ShaoZhang2014} go beyond traditional correlation and instead can indicate whether random quantities are independent or martingale differences, respectively.
For time series, in the time domain, \cite{Zhou2012} introduced auto distance correlations that are zero if and only if the measured time series components are independent. \cite{linton2007quantilogram}, and \cite{davis2009extremogram} introduced the (univariate) concepts of quantilograms and extremograms, respectively. More recently, quantile correlation \citep{SchmittEtAl2015}, and quantile autocorrelation functions \citep{li2014quantile} together with cross-quantilograms \citep{han2014cross} have been proposed as a fundamental tool for analysing dependence in quantiles of the distribution.

In the frequency domain, \cite{Hong1999} introduced a generalised spectral density. In the generalised spectral density covariances are replaced by quantities that are closely related to empirical characteristic functions. In \cite{Hong2000} the Fourier transform of empirical copulas at different lags is considered for testing the hypothesis of pairwise independence. Recently, under the names of \mbox{Laplace-,} quantile and copula spectral density and spectral density kernels, various quantile-related spectral concepts have been proposed, for the frequency domain. The approaches by~\cite{hagemann2013robust} and \cite{li2008laplace,li2012quantile} are designed to consider cyclical dependence in the distribution at user-specified quantiles. \cite{MikoschZhao2014, MikoschZhao2015} define and analyse a periodogram (and its integrated version) of extreme events. As noted by~\cite{hagemann2013robust} other approaches aim at discovering ``the presence of \emph{any} type of dependence structure in time series data'', referring to work of \cite{dette2014copulas} and \cite{lee2012quantile}. This comment also applies to~\cite{kley2014}. In the present paper our aim is to generalise the existing approaches to multivariate time series. The extensions to the terminology that we provide, in particular the introduction of the standardised quantile coherency, is very important for economic applications, because it enables the analyst to perform a more detailed joint analysis of the serial and cross sectional dependence in multiple time series.

For the univariate case different approaches to consistent estimation were considered. \cite{li2008laplace} proposed an estimator for a weighted version of the quantile spectra, based on least absolute deviation regression, for the special case where $\tau_1 = \tau_2 = 0.5$. \cite{li2012quantile} generalised the estimator, using quantile regression, to the case where $\tau_1 = \tau_2 \in (0,1)$. The general case, in which the quantities can be related to the copulas of pairs, was first considered by \cite{dette2014copulas}. These authors also were the first to consider a rank-based version of the quantile regression-type estimator which eliminates the need to estimates the weights in~\cite{li2008laplace,li2012quantile}. For the case where $\tau_1 = \tau_2 \in (0,1)$, \cite{hagemann2013robust} proposed a version of the traditional $L^2$-periodogram where the observations are replaced with $I\{\hat F_{n,j}(X_{t, j}) \leq \tau\} = I\{R_{n;t,j} \leq n \tau\}$. \cite{kley2014} generalised this estimator, in the spirit of~\cite{dette2014copulas}, by considering cross-periodograms for arbitrary couples $(\tau_1, \tau_2) \in [0,1]^2$, and proved that it converges, as a stochastic process, to a complex-valued Gaussian limit. An estimator defined in analogy to the traditional lag-window estimator was analysed by~\cite{birr2015quantile} in the context of non-stationary time series.

\section{Asymptotic properties of the proposed estimators}
\label{sec:asymp}


To derive the asymptotic properties of the estimators defined in Section~\ref{sec:ccsdk:estim} some assumptions need to be made. Recall (cf.~\citet{Brillinger1975}, p.~19) that the {$r$th order joint cumulant} $\cum(Z_1, \ldots, Z_r)$ of the random vector $(Z_1, \ldots, Z_r)$ is defined as
\[\cum(Z_1, \ldots, Z_r) := \sum_{\{\nu_1, \ldots, \nu_p\}} (-1)^{p-1} (p-1)! E \Big[ \prod_{j \in \nu_1} Z_j \Big] \cdots E \Big[ \prod_{j \in \nu_p} Z_j\Big],\]
with summation extending over all partitions $\{\nu_1, \ldots, \nu_p\}$, $p=1,\ldots,r$, of~$\{1,\ldots,r\}$.

\noindent Regarding the range of dependence of $(\bX_t)_{t \in \IZ}$ we make the following assumption,
\begin{assumption}\label{ass:exp_alpha_mix}
	The process $(\bX_t)_{t \in \IZ}$ is strictly stationary and exponentially $\alpha$-mixing, that is, there exists constants $K < \infty$ and $\rho \in (0,1)$, such that
\begin{equation}\label{eq:exp_alpha_mix}
\alpha(n) := \sup_{\substack{A \in \sigma(X_0, X_{-1}, \ldots)\\ B \in \sigma(X_n, X_{n+1}, \ldots)}} \big| \IP(A \cap B) - \IP(A) \IP(B) \big| \leq K \rho^n, \quad n \in \IN.
\end{equation}
\end{assumption}

\noindent Further, to establish consistency of the estimates we consider sequences of weights that asymptotically concentrate around multiples of $2\pi$,
\begin{assumption}\label{ass:W}
The weights are defined as
$
W_n(u) := \sum_{j=-\infty}^{\infty} b_n^{-1} W(b_n^{-1} [u + 2\pi j])
$, where $b_n > 0$, $n=1,2,\ldots$, is a sequence of scaling parameters satisfying~$b_n \rightarrow 0$ and $n b_n \rightarrow \infty$, as $n \rightarrow \infty$.
The weight function $W$ is real-valued, even, has support $[-\pi,\pi]$, bound\-ed variation, and satisfies
$
\int_{-\pi}^{\pi} W(u) \text{d}u = 1.
$
\end{assumption}
Comments on the assumptions will follow in the end of this section.
The main result of this section (Theorem~\ref{thm:AsympCohRankEstimator}) will legitimise $\hat\bfrakR_{n,R}(\omega;\tau_1,\tau_2)$ as an estimator of the quantile coherency $\bfrakR(\omega;\tau_1,\tau_2)$. Results that legitimise $\bI_{n,R}(\omega; \tau_1, \tau_2)$ and $\hat\bG_{n,R}(\omega; \tau_1, \tau_2)$ as estimators of the quantile cross-spectral density $\bfrakf(\omega; \tau_1, \tau_2)$ are deferred to the supplementary material to not impair the flow of the paper. The legitimacy of the estimates follows from the fact that the estimators converge weakly in the sense of \emph{Hoffman-J\o{}rgensen} (cf. Chapter 1 of~\cite{vanderVaartWellner1996}). We denote this mode of convergence by $\Rightarrow$ . The estimators under consideration take values in the space of (element-wise) bounded functions $[0,1]^2 \rightarrow \IC^{d \times d}$,  which we denote by $\ell_{\mathbb{C}^{d \times d}}^{\infty}([0,1]^2)$. While results in empirical process theory are typically stated for spaces of real-valued, bounded functions, these results transfer immediately by identifying $\ell_{\mathbb{C}^{d \times d}}^{\infty}([0,1]^2)$ with the product space $\ell^{\infty}([0,1]^2)^{2 d^2}$. Note that the space $\ell_{\mathbb{C}^{d \times d}}^{\infty}([0,1]^2)$ is constructed along the same lines as the space $\ell_{\mathbb{C}}^{\infty}([0,1]^2)$ in \cite{kley2014}.

\noindent We are now ready to state the main result of this section.
\begin{theorem}\label{thm:AsympCohRankEstimator}
Let Assumptions~\ref{ass:exp_alpha_mix} and~\ref{ass:W} hold. Assume that the marginal distribution functions $F_j$, $j=1,\ldots,d$ are continuous and that constants $\kappa > 0$ and $k \in\IN$ exist, such that
$b_n = o(n^{-1/(2k+1)})$ and $b_n n^{1-\kappa} \rightarrow \infty$. Assume that for some $\varepsilon \in (0,1/2)$ we have
$\inf_{\tau \in [\varepsilon, 1-\varepsilon]} \mathfrak{f}^{j, j}(\omega; \tau, \tau) > 0$,
for all $j=1,\ldots,d$. Then, for any fixed $\omega \in  \IR$,
\begin{equation}\label{eqn:convR}
\sqrt{n b_n} \Big( \hat\bfrakR_{n,R}(\omega; \tau_1, \tau_2) - \bfrakR(\omega;\tau_1,\tau_2) - \bfrakB_n^{(k)}(\omega; \tau_1, \tau_2) \Big)_{(\tau_1, \tau_2) \in [\varepsilon, 1-\varepsilon]^2}
\Rightarrow \IL(\omega; \cdot, \cdot),
\end{equation}
in $\ell_{\mathbb{C}^{d \times d}}^{\infty}([\varepsilon,1-\varepsilon]^2)$, where
\begin{equation}\label{eqn:defL}
	\Big\{ \IL(\omega; \tau_1, \tau_2) \Big\}_{j_1, j_2} :=
	\frac{1}{\sqrt{\ff{1}{1} \ff{2}{2}}} \Big( \HH{1}{2}
		- \frac{1}{2} \frac{\ff{1}{2}}{\ff{1}{1}} \HH{1}{1} - \frac{1}{2} \frac{\ff{1}{2}}{\ff{2}{2}} \HH{2}{2} \Big),
\end{equation}
\begin{equation} \label{def:biasR}
\Big\{ \bfrakB_n^{(k)}(\omega; \tau_1, \tau_2) \Big\}_{j_1, j_2} := 
\frac{1}{\sqrt{\ff{1}{1} \ff{2}{2}}} \Big( \BB{1}{2}
		- \frac{1}{2} \frac{\ff{1}{2}}{\ff{1}{1}} \BB{1}{1} - \frac{1}{2} \frac{\ff{1}{2}}{\ff{2}{2}} \BB{2}{2} \Big)
\end{equation}
and we have written $\ff{a}{b}$ for the quantile cross-spectral density $\mathfrak{f}^{j_a, j_b}(\omega; \tau_a, \tau_b)$ as defined in \eqref{eqn:ccsdk},
$\BB{a}{b} := \sum_{\ell=2}^k \frac{b_n^{\ell}}{\ell!} \int_{ -\pi}^{\pi} v^{\ell} W(v) dv \frac{{\rm d}^{\ell}}{{\rm d}\omega^{\ell}}\mathfrak{f}^{j_a, j_b}(\omega; \tau_a, \tau_b)$, and $\HH{a}{b}$ for $\IH^{j_a, j_b}(\omega; \tau_a, \tau_b\big)$; a component of $\IH(\omega;\cdot,\cdot) := (\IH^{j_1, j_2}(\omega;\cdot,\cdot))_{j_1, j_2 = 1,\ldots,d}$ defined as a centred, $\IC^{d \times d}$-valued Gaussian process characterised by
\begin{multline} \label{eqn:CovH}
\Cov\big(\IH^{j_1, j_2}(\omega; u_1, v_1\big), \IH^{k_1, k_2}(\lambda; u_2, v_2)\big) \\ = 2\pi \Big(\int_{-\pi}^\pi W^2(\alpha){\rm d}\alpha \Big)
\Big( \mathfrak{f}^{j_1, k_1}(\omega; u_1, u_2) \mathfrak{f}^{j_2, k_2}(-\omega; v_1, v_2) \eta(\omega - \lambda) \\
+ \mathfrak{f}^{j_1, k_2}(\omega; u_1, v_2) \mathfrak{f}^{j_2, k_1}(-\omega; v_1, u_2)
\eta(\omega + \lambda) \Big),
\end{multline}
where $\eta(x) := I\{x = 0 (\mod 2\pi)\}$ [cf.~\cite[p.\,148]{Brillinger1975}] is the $2\pi$-periodic extension of Kronecker's delta function.
The family~$\{\IH(\omega; \, \cdot, \cdot),$ $ \omega \in  [0,\pi] \}$ is a collection of independent processes and $\IH(\omega; \tau_1, \tau_2) = \overline{\IH(-\omega; \tau_1, \tau_2)} = \IH(\omega+2\pi; \tau_1, \tau_2)$.
\end{theorem}
The proof of Theorem~\ref{thm:AsympCohRankEstimator} is lengthy and technical and therefore delegated to the online supplement (Section~\ref{sec:proofs}).
Comparing Theorem~\ref{thm:AsympCohRankEstimator} with results for the traditional coherency (see, for example, Theorem~7.6.2 in \cite{Brillinger1975}) we observe that the distribution of $\hat\bfrakR_{n,R}(\omega; \tau_1, \tau_2)$ is asymptotically equivalent to that of the traditional estimator [cf.~(7.6.14) in~\cite{Brillinger1975}] computed from the unobserved time series
\begin{equation}\label{eqn:clippedTS}
\big(I\{F_{j_1}(X_{t,j_1}) \leq \tau_1\},I\{F_{j_1}(X_{t,j_2}) \leq \tau_2\}\big), \quad t = 0, \ldots, n-1.
\end{equation}

The convergence to a Gaussian process in~\eqref{eqn:convR} can be employed to obtain asymptotically valid pointwise confidence bands. To this end, the covariance kernel of $\IL$ can easily be determined from~\eqref{eqn:defL} and~\eqref{eqn:CovH}, yielding an expression similar to (7.6.16) in~\cite{Brillinger1975}. A more detailed account on how to conduct inference is given in Section~\ref{sec:CI} of the supplementary material. Note that the bound to the order of the bias given in (7.6.15) in~\cite{Brillinger1975} applies to the expansion given in~\eqref{def:biasR}.

If $W$ is a kernel of order $p\geq 1$ we have that the bias is of order $b_n^p$. Thus, if we choose the mean square error minimising bandwidth $b_n \asymp n^{-1/(2p + 1)}$ the bias will be of order $n^{-p/(2p + 1)}$.
Regarding the restriction $\varepsilon > 0$, note that the convergence~\eqref{eqn:convR} can not hold if $(\tau_1, \tau_2)$ is on the border of the unit square, as the quantile coherency $\bfrakR(\omega;\tau_1,\tau_2)$ is not defined if $\tau_j \in \{0,1\}$, as this implies that $\Var(I\{F_{j}(X_{t,j}) \leq \tau_j\}) = 0$.

We now comment on the assumptions:
Assumption~\ref{ass:exp_alpha_mix} holds for a wide range of popular, linear and nonlinear processes. Examples (possibly, under mild additional assumptions) include the traditional VARMA or vector-ARCH models as well as many others. It is important to observe that Assumption~\ref{ass:exp_alpha_mix} does not require the existence of any moments, which is in sharp contrast to the classical assumptions, where moments up to the order of the respective cumulants have to exist. Assumption~\ref{ass:W} is quite standard in classical time series analysis [cf., for example, \citet{Brillinger1975}, p.\,147].

\section{Quantile cross-spectral analysis of stock market returns: A route to more accurate risk measures?}
\label{sec:emp}

Stock market returns belong to one of the prominent datasets in economics and finance. Although many important stylised facts about their behaviour have been established in the past decades, it remains a very active area of research. Despite the efforts, an important direction, which has not been fully addressed is stylised facts about the joint distribution of returns. Especially during the last turbulent decade, understanding the behaviour of joint quantiles in return distributions became particularly important, as it is essential for understanding systemic risk; {``the risk that the intermediation capacity of the entire system can be impaired''}; cf. \cite{adrian2011covar}. Several authors focus on explaining tails of the bivariate market distributions in different ways. \cite{adrian2011covar} proposed to classify institutions according to the sensitivity of their quantiles to shocks to the market. Most closely related to the notion of how we view the dependence structures is the multivariate regression quantile model of \cite{white2015var}, which studies the degree of tail interdependence among different random variables directly.

Quantile cross-spectral analysis, as designed in this paper, allows to analyse the fundamental dependence quantities in the tails (but also in any other part) of the joint distribution and across frequencies. An application to stock market returns may therefore provide deeper insight about dependence in stock markets, and lead to a more powerful analysis securing us against financial collapses.

One of the important features of stock market returns is time variation in its volatility. Time-varying volatility processes can cross almost every quantile of their distribution (cf. \cite{hagemann2013robust}), and create peaks in quantile spectral densities as shown by \cite{li2014JTSA}. These notions have recently been documented by \cite{engle2004caviar} and \cite{vzikevs2014semi} who propose models for the conditional quantiles of the return distribution based on the past volatility. In the multivariate setting, strong common factors in volatility are found by \cite{barigozzi2014disentangling} who conclude that common volatility is an important risk factor. Hence, common volatility should be viewed as a possible source of dependence. Because we aim to find the common structures in the joint distribution of returns, we study returns standardised by its volatility that we estimate by a GARCH(1,1) model; cf. \cite{bollerslev1986generalized}. This first step is commonly taken in the literature of modelling the joint market distribution using copulas; cf. \cite{granger2006common} and \cite{patton2012review}. In these approaches the volatility in the marginal distributions is modelled first, and the common factors are then considered in the second step. Consequently, this will allow us to discover other possible common factors in the joint distribution of market returns across frequencies, that result in spurious dependence, but which will not be overshadowed by the strong volatility process.

We choose to study the joint distribution of portfolio returns and excess returns on the broad market, hence looking at one of the most commonly studied factor structures in the literature as dictated by asset pricing theories; cf. \cite{sharpe1964capital} and \cite{lintner1965valuation}. As an excess return on the market, we use value-weighted returns of all firms listed on the NYSE, AMEX, or NASDAQ from the Center for Research in Security Price (CRSP) database. For the benchmark portfolio, we use an industry portfolio formed from consumer non-durables.\footnote{Note to choice of the data: we use the publicly available data available and maintained by Fama and French at \url{http://mba.tuck.dartmouth.edu/pages/faculty/ken.french/data_library.html}. This data set is popular among researchers, and while many types of portfolios can be chosen, we chose consumer non-durables randomly for this application. Although very interesting and attractive, it is far beyond the scope of this work to present and discuss results for wider portfolios formed on distinct criteria.} We used $n=23385$ daily observations (from 1 July 1926 through to 30 June 2015). The data includes several crisis periods and therefore might not be suitable to be viewed as a strictly stationary time series. Nevertheless, we choose to study this long period of data as we believe that longer than yearly cycles might constitute an important possible source of dependence, and we believe the empirical results are practically interesting. Moreover, by standardising the returns by their volatility we removed what we believe is the most important source of time-variation in data.\footnote{As a robustness check, we have sliced the time series into decades and found that our results on non-overlapping windows do not materially change.}


\begin{figure}[t!]
  \begin{center}
      \includegraphics[width=0.64\textwidth]{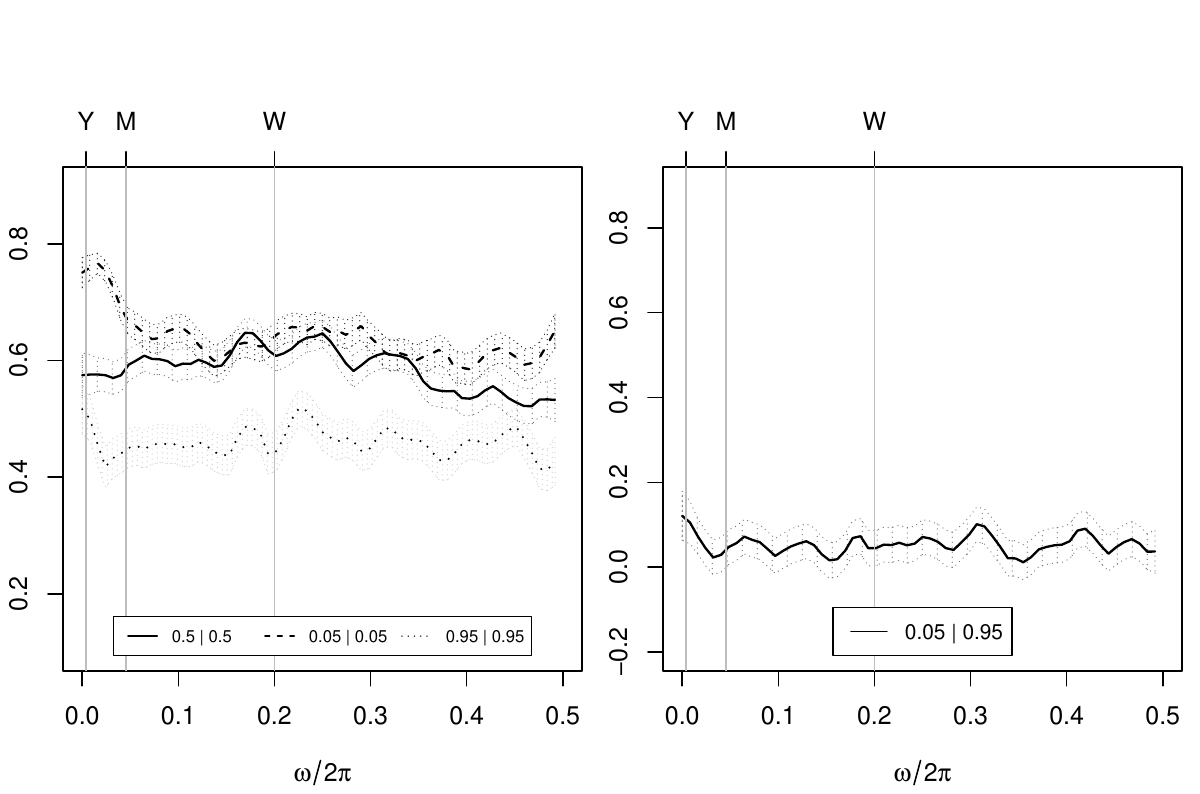}
      \includegraphics[width=0.32\textwidth]{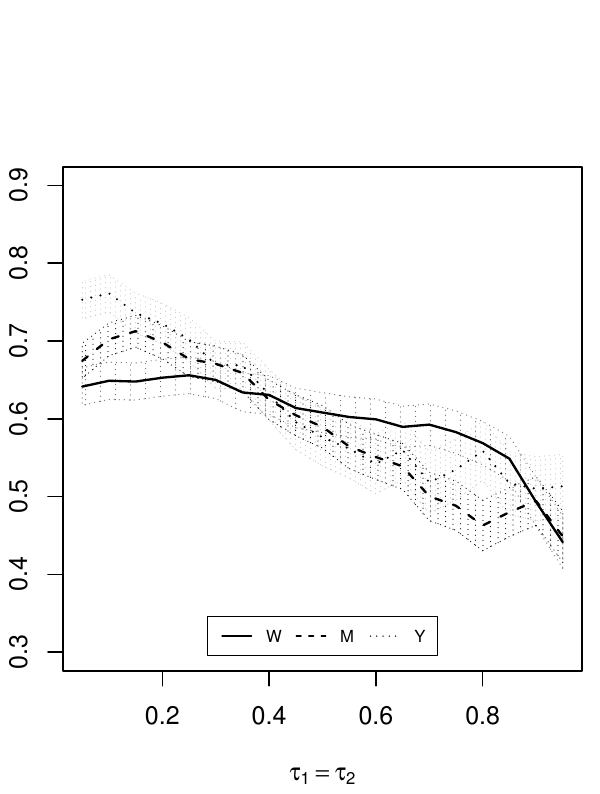} 
   \end{center}
	\vspace*{-0.5cm}
	\caption{\label{fig:emp}Quantile coherency estimates for the portfolio.}
\end{figure}
In the left panel of Figure~\ref{fig:emp}, quantile coherency estimates for the $0.05|0.05$, $0.5|0.5$, and $0.95|0.95$ combinations of quantile levels of the joint distribution are shown for the industry portfolio and excess market returns over frequencies. The centre panel in Figure~\ref{fig:emp}, on which we comment later, shows the $0.05|0.95$ combination. We have used the Epanechnikov kernel and a bandwidth of $b_n = 0.5 n^{1/4}$ for the computation of the estimates (cf. \eqref{def:Rhat}). The confidence intervals, shown as dotted regions, are at the 95\% level and were constructed according to the procedure described in Section~\ref{sec:CI} of the supplementary material.
For clarity, we plot the $x$-axis in daily cycles and also indicate the frequencies that correspond to yearly, monthly, and weekly periods. While we use daily data the highest possible frequency of 0.5 indicates 0.5 cycles per day (i.\,e., a 2-day period). While precise frequencies do not have an economic meaning, one needs to understand the interpretation with respect to the time domain. For example, a sampling frequency of 0.2 corresponds to 0.2 cycles per day translating to a 5 days period (equivalent to one week), but the frequency of 0.3 translates to a hardly interpretable $3.\overline{3}$ period. Hence, the upper label of the $x$-axis is of particular interest to an economist, as one can study how weekly, monthly, or yearly cycles are connected across quantiles of the joint distribution. For the clarity of presentation, we focus on the real part of the quantities, which relates to the dynamics of the process switching between the $j_1$st component being below the $\tau_1$-quantile and the $j_2$nd component being above the $\tau_2$-quantile (cf. Section~\ref{sec:ccsdk}).

The real parts of the quantile coherency estimates reveal frequency dynamics in quantiles of the joint distribution of the returns under study. Generally, cycles at the lower quantiles appear to be more strongly dependent than at the upper quantiles, which is a well documented stylised fact about stock market returns. It points us to the fact that returns are more dependent during business cycle downturns, than upturn; cf. \cite{erb1994forecasting}, \cite{longin2001extreme}, \cite{ang2002asymmetric} and \cite{patton2012review}. More importantly, lower quantiles are more strongly related in periods longer than one week on average in comparison to shorter than weekly periods, and are even more connected at longer than monthly cycles. This suggests that infrequent clusters of large negative portfolio returns are better explained by excess market returns than small daily fluctuations. Returns in upper quantiles of the joint distribution seem to be connected similarly across all frequencies. The same result holds also for the median. For a better exposure, we also present quantile coherency estimates for three fixed weekly, monthly, and yearly periods (corresponding to $\omega \in 2\pi \{1/5,1/22,1/250\}$, respectively) at all quantile levels $\tau_1 = \tau_2 \in \{0.05, 0.1, \ldots, 0.95\}$ in the right panel of Figure~\ref{fig:emp}. This alternative plot highlights the previous discussion.

We now compare our findings to a corresponding analysis with the cross-quantilogram, a related quantile-based measure for serial dependence in the time domain. Considering a strictly stationary, $\IR \times \IR \times \IR^{d_1} \times \IR^{d_2}$-valued time series $(y_{1t}, y_{2t}, x_{1t}, x_{2t})$, with $t \in \IZ$ and $d_{1}, d_{2} \in \IN$, denoting the conditional distribution of the series $y_{it}$ given $x_{it}$ by $F_{y_i|x_i}(\cdot | x_{it})$, and the quantile function as $q_{i,t}(\tau_i) = \inf\{v : F_{y_i|x_i}(\cdot | x_{it}) \geq \tau_i\}$, $\tau_i \in (0,1)$, $i=1,2$; \cite{han2014cross} define the cross-quantilogram as
\[\rho_{(\tau_1, \tau_2)}(k) := \frac{\E\big[ ( I\{y_{1t} < q_{1,t}(\tau_1)\} - \tau_1 ) ( I\{y_{2,t-k} < q_{2,t-k}(\tau_2)\} - \tau_2 ) \big]}{
\Big(\E\big[ ( I\{y_{1t} < q_{1,t}(\tau_1)\} - \tau_1 )^2 \big]
\E\big[ ( I\{y_{2,t-k} < q_{2,t-k}(\tau_2)\} - \tau_2 )^2 \big] \Big)^{1/2} }.\]
With no covariate information in our data example, this reduces to $x_{1t} = x_{2t} = 1$ and $q_{i,t}$ being the quantile of the marginal distribution of $y_{it}$. It is important to note that the cross-quantilogram is defined as a standardised measures of serial dependencies between the events $\{y_{1t} \leq q_{1,t}(\tau_1)\}$ and $\{y_{2t} \leq q_{2,t}(\tau_2)\}$ in the time domain, while quantile coherency is defined similarly, but in the frequency domain.

\begin{figure}[t!]
  \begin{center}
        \includegraphics[width=0.32\linewidth]{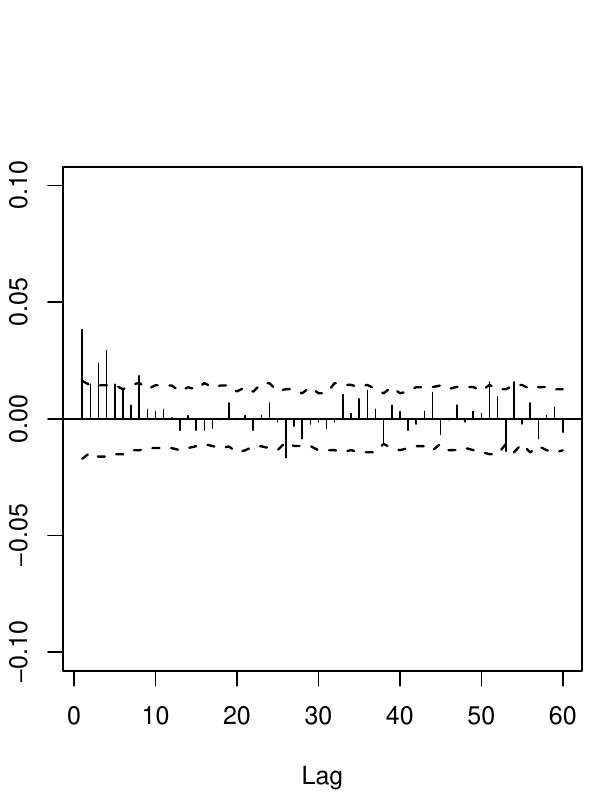} 
        \includegraphics[width=0.32\linewidth]{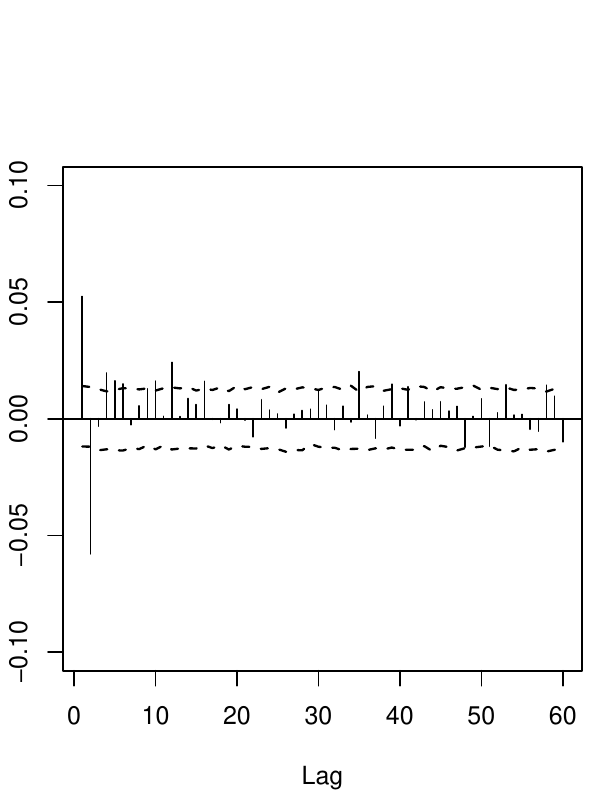}
        \includegraphics[width=0.32\linewidth]{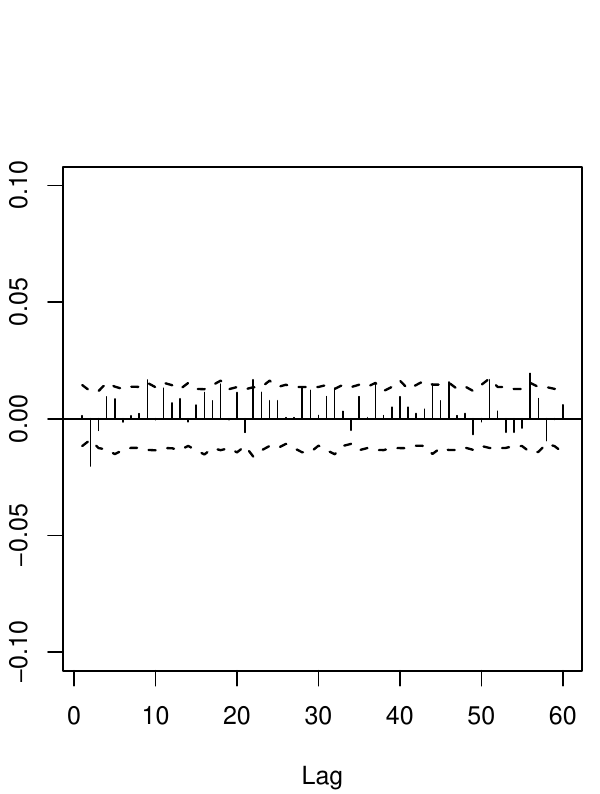} 
   \end{center}
  \vspace*{-0.5cm}
  \caption{\label{fig:emp2}Cross-quantilogram estimates for the portfolio.}
\end{figure}
In Figure~\ref{fig:emp2} we present the cross-quantilograms that we estimated from our data example. For the computation we have used the estimator and stationary bootstrap procedure defined in \cite{han2014cross}. More precisely, we used the implementation that is available in the R package \pkg{quantilogram}; cf. \cite{quantilogram_Rpkg}. Inspecting the plots, it can be seen that there are lags~$k$, typically short, where significant dependence is present. Further, it is possible to guess that there is periodic variation of positive and negative dependence at the 0.05 quantile level, while at the 0.95 quantile level the dependence seems to be largely positive. Yet, taking into account the confidence intervals, it is uncertain if this is a significant pattern.
Further, comparing the discussion of these periodic patterns shown by cross-quantilogram with what we were able to read from quantile coherency in Figure~\ref{fig:emp}, it is difficult to read specific weekly, monthly and yearly periodic components and whether or not they are significant. Thus, at least in the specific case where a researcher is interested in the dependence of cycles, we believe that quantile coherency can provide a perspective that is unavailable in the time domain analysis.

To summarise the result of our empirical analysis: while asymmetry is commonly found by researchers, we document frequency dependent asymmetry of stock market returns (i.\,e., asymmetry with respect to cycles in the joint distribution). In case this behaviour would be common across larger classes of assets, our results may have large implications for one of the cornerstones of asset pricing theory assuming normal distribution of returns. It leads us to the call for more general models, and more importantly to the need of restating the asset pricing theory in a way that allows to distinguish between short run and long run behaviour of investors.

Our results are also crucial for systemic risk measurement, as an investor wishing to optimise a portfolio should focus on stocks which will not be connected at lower quantiles, in a situation of distress, but will be connected at upper quantiles, in a situation of market upturns in a given investment period. We document behaviour which is not favourable to such an investor using traditional pricing theories, as we show that broad stock market returns contain a common factor more frequently during downturns than during upturns. This suggests that the portfolio at hand might be much riskier than it were implied by common measures. Further, our results suggest that this effect becomes even worse for long-run investors.

An important feature of our quantile cross-spectral measures is that they enable us to measure dependence also between $\tau_1 \ne \tau_2$ quantiles of the joint distribution. In the central panel of Figure \ref{fig:emp} we document that the dependence between the $0.05|0.95$ quantiles of the return distribution is not very strong. Generally speaking, no intense dependence can be seen between large negative returns of the stock market, and large positive returns of the portfolio under study. This kind of analysis may be even more interesting in the case where dependence between individual assets is studied. There, negative news may have strong opposite impact on the assets under study.

Finally, some words of caution to the reader, about the interpretation of the quantities which we have estimated, are in order. In Section~\ref{sec:interpretationgauss} of the supplementary material we provide a link between quantile coherency and traditional measures of dependence under the assumption of normally distributed data. The quantile-based measures are designed to capture general dependence types without restrictive assumptions on the underlying distribution of the process. Hence, here we have intentionally not relate it to traditional correlation which, ideally, should only be interpreted when the process is known to be Gaussian. The financial returns under study in this section are known to depart from normality. Therefore, quantile coherency is not directly comparable to traditional correlation measures. What we can see is generally strong dependence between the portfolio returns and excess market returns at all quantiles confirming the fact that excess returns are a strong common factor for the studied portfolio returns. The details that the quantile-based analysis in this section revealed would have remained hidden in an analysis based on the traditional coherency.

\section{Quantile coherency in a model assessing exercise}
\label{sec:modexerc}

In the previous section we demonstrated how quantile coherency can be used by applied researchers to reveal cyclical features of the data that might remain invisible if the data is analysed solely with covariance-based dependency measures. In this section we illustrate how quantile coherency can be used to assess the capability of time series models to capture such cycles documented in the data.

More precisely, we fit several bivariate time series models and then compare the quantile coherencies implied by estimated parameters with those obtained from a non-parametric estimation (cf. Figure~\ref{fig:emp}).
The graphical approach of assessing the models is similar to the one proposed in~\cite{BirrEtAl2018}.
For the sake of clarity, we focus on two classes of models: (a) vector autoregressive (VAR) models, and (b) vector versions of the quantile autoregressive (QVAR) model introduced by \cite{koenker2006quantile}. Classical VAR used by many applied researchers assumes the same autoregressive structure at all quantiles. To model asymmetry, one can employ more flexible copulas allowing for asymmetric dependence. In addition, QVAR allows different autoregressive structure at different quantiles. Hence different quantiles can be driven by processes with different cyclical properties. 

\begin{figure}[t!]
  \begin{center}
      \includegraphics[width=0.32\textwidth]{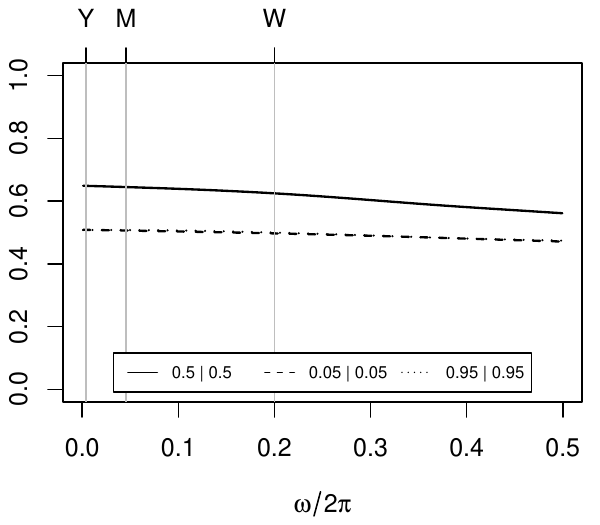}
			\includegraphics[width=0.32\textwidth]{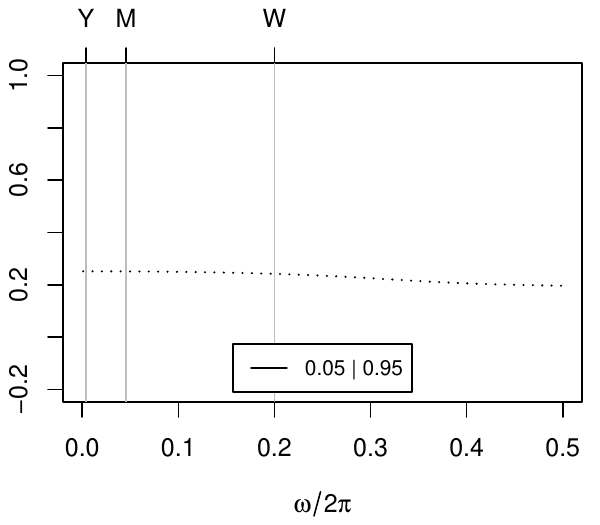}
			\includegraphics[width=0.32\textwidth]{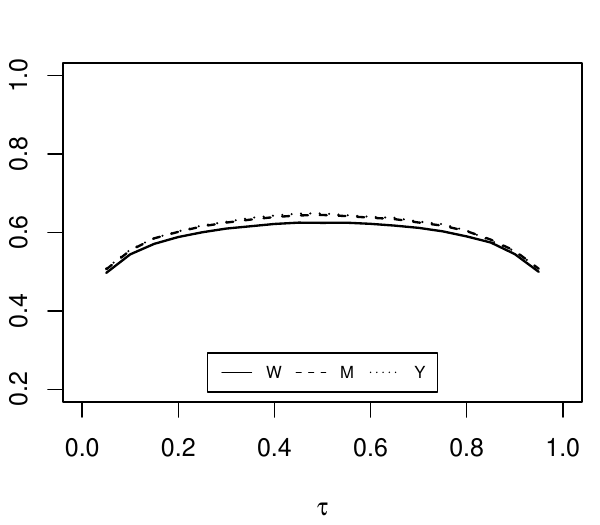}
			
			\includegraphics[width=0.32\textwidth]{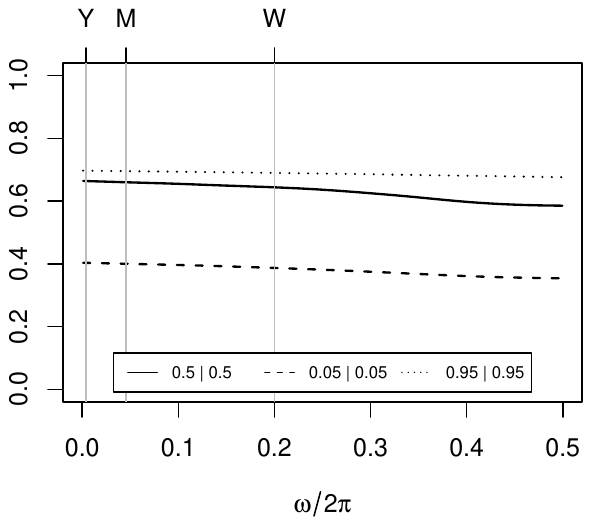}
			\includegraphics[width=0.32\textwidth]{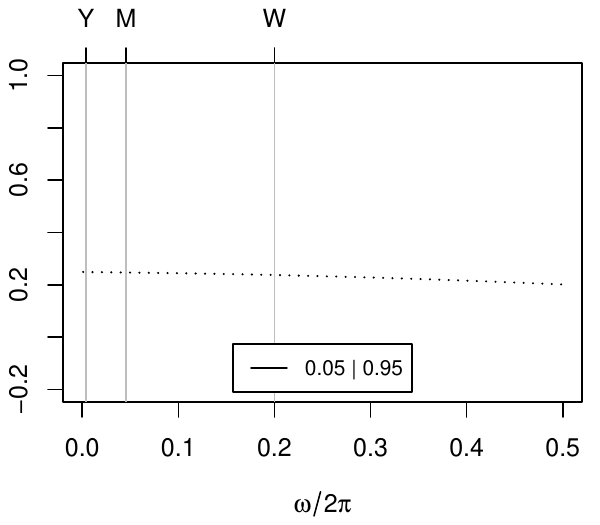}
			\includegraphics[width=0.32\textwidth]{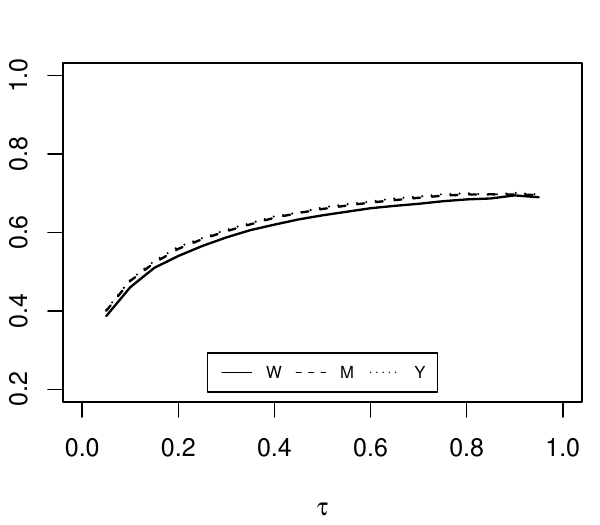}
			
			\includegraphics[width=0.32\textwidth]{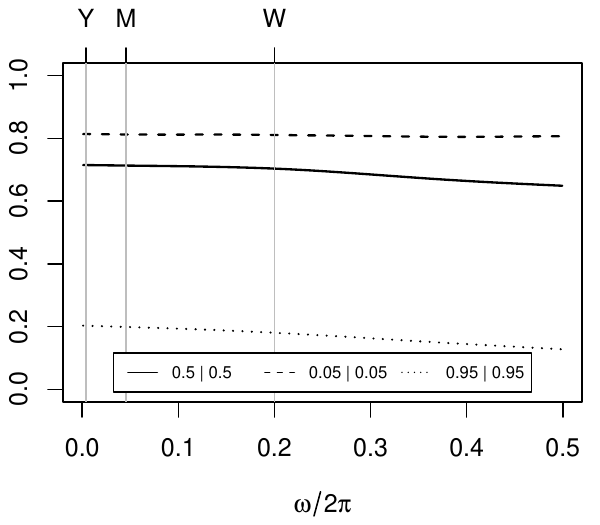}
			\includegraphics[width=0.32\textwidth]{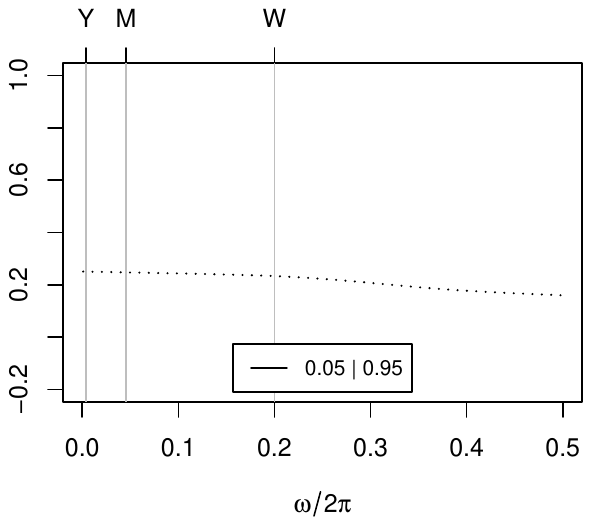}
			\includegraphics[width=0.32\textwidth]{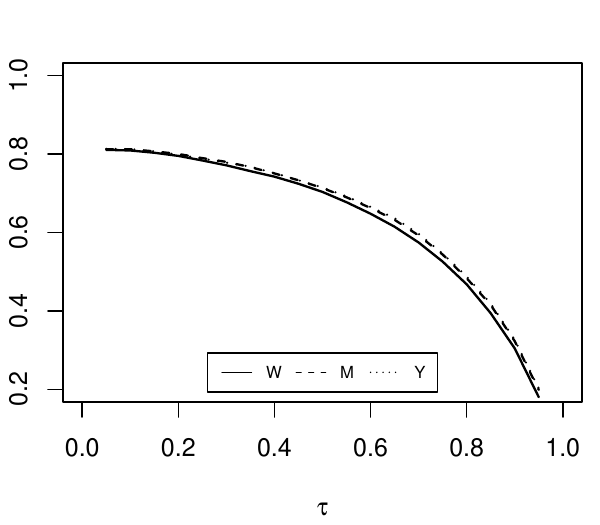}
   \end{center}
	\vspace*{-0.5cm}
	\caption{\label{fig:qc_var}Quantile coherency simulated from the VAR models.}
\end{figure}

We discuss the models in order, from simple to more complex, and evaluate if the more complex models are better suited to capture the weekly, monthly and yearly cycles of quantile-related features which were discovered in the stock market returns analysis of Section~\ref{sec:emp}.

We begin by fitting a VAR(1) to the stock market returns. The fitted model is
\begin{equation}\label{eqn:var_model}
\begin{split}
	Y_{t,1} & = 0.0987 + 0.056 Y_{t-1,1} + 0.186 Y_{t-1,2} + \varepsilon_{t,1},\\
	Y_{t,2} & = 0.0369 - 0.056 Y_{t-1,1} + 0.175 Y_{t-1,2} + \varepsilon_{t,2},\\
\end{split}
\end{equation}
where $(\varepsilon_{t,1}, \varepsilon_{t,2})$ is white noise with an estimated ${\rm Corr}(\varepsilon_{t,1}, \varepsilon_{t,2}) \approx 0.822$. Adding the common assumption that the $(\varepsilon_{t,1}, \varepsilon_{t,2})$ are independent and jointly Gaussian, the corresponding quantile coherencies can be determined. Quantile coherencies implied by the model \eqref{eqn:var_model} are depicted in the top row of Figure~\ref{fig:qc_var}. For easier comparison, we consider the same combinations of frequencies and quantile levels as in Figure~\ref{fig:emp}. In the picture it is clearly visible that dependencies of cycles implied by this Gaussian models are symmetric. For example, the dependence at the $0.05 | 0.05$ and at the $0.95 | 0.95$ level are equally strong for all frequencies. In contrast, the nonparametric estimate obtained from the data (cf. Figure~\ref{fig:emp}) shows strong asymmetry. Further, we can see that for the weekly, monthly and yearly frequencies, which might be of particular interest for applied researchers, the dependencies at the $\tau | \tau$ and at the $1-\tau | 1-\tau$ level coincide as well. If an applied researcher seeks to model dependencies as the ones revealed in Section~\ref{sec:emp}, the Gaussian VAR model might therefore be too restrictive.

Next, we consider non-Gaussian versions of the fitted VAR. To obtain these models, note that the innovations in~\eqref{eqn:var_model} are assumed to be white noise, but are not required to be i.\,i.\,d. Gaussian. Another plausible model is therefore obtained by specifying any joint distribution for $(\varepsilon_{t,1}, \varepsilon_{t,2})$ that has first and second moment as implied by the fitted VAR model. For illustration we now consider the following two cases. In both cases we assume the marginal distributions to be standard normal. In the first case we assume that the dependence is according to a Clayton copula with parameter $\theta = 4$. In the second case we assume that it is according to a Gumbel copula with parameter $\theta = 2.7$. As one might expect, the dependence in the tails of the VAR(1) process is now remarkably different. As it can be seen from the middle-left plot in Figure~\ref{fig:qc_var}, for the case of the Clayton copula there is stronger dependence in the lower tail ($0.05 | 0.05$) and weaker dependence in the upper tail ($0.95 | 0.95$). The dependence is slightly stronger for low frequencies, which is expected from the temporal dependence in the VAR model. In the bottom-left plot of Figure~\ref{fig:qc_var}, on the other hand, we see stronger dependence in the upper and weaker dependence in the lower tail. Interestingly, as can be seen from the centre plots, the dependence of cycles in changing from being below the 0.05-quantile in the first component to being below the 0.95-quantile in the second component does not depend much on the choice of the copula. Finally, in the right plots of Figure~\ref{fig:qc_var}, we see how the dependence changes according to the quantile level when cycles at the weekly, monthly and yearly frequencies, which we think might be most relevant to some practitioners, are considered. As expected, we see that for the case of the Clayton copula the dependence decreases as the quantile level $\tau$ increases, where for the case of the Gumbel copula the dependence increases if $\tau$ increases. Although the models with the Gumbel and Clayton copula capture asymmetric dependence better than the one with the Gaussian copula, we can still see that they depart from the data in terms of quantile coherency.

\begin{figure}
  \begin{center}
		\includegraphics[width=0.75\textwidth]{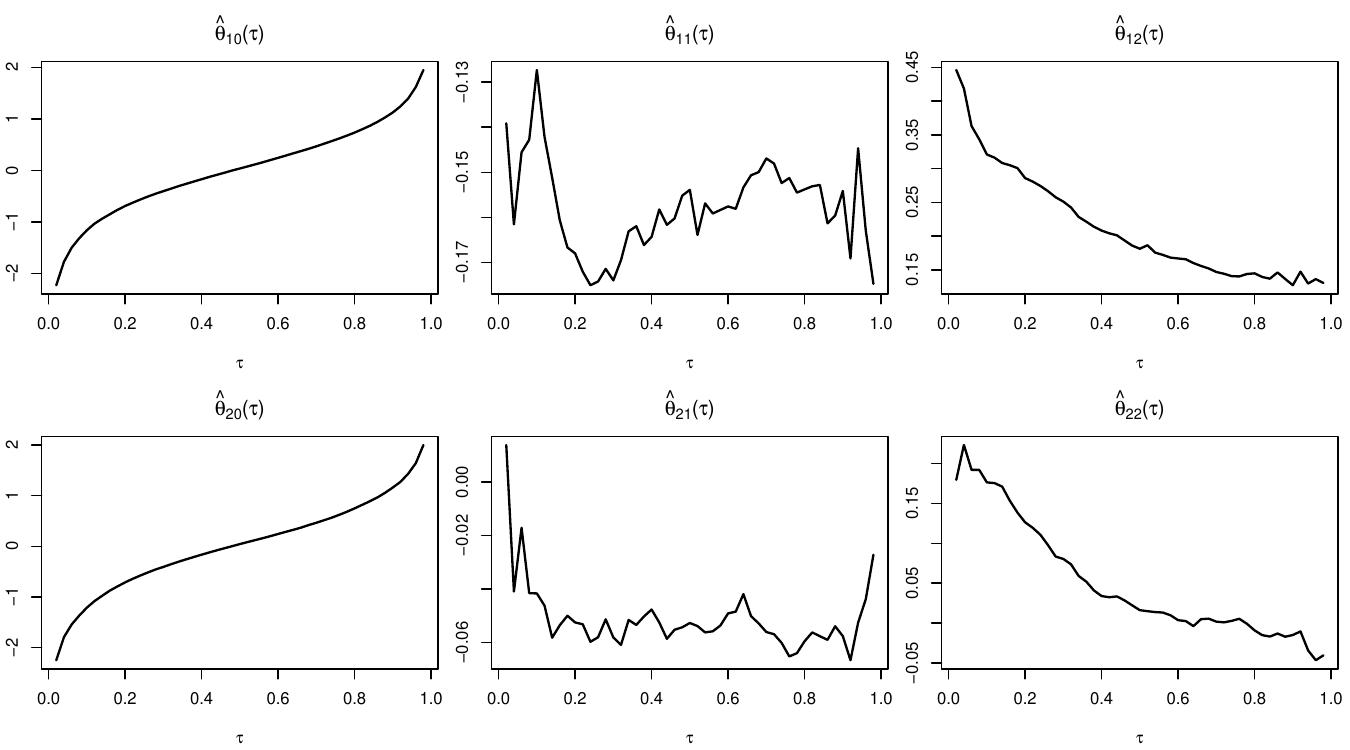}
  \end{center}
	\vspace*{-0.5cm}
	\caption{\label{fig:qvar_a_coef}Estimated parameter functions for model~\eqref{eqn:qvar_a}.}
\end{figure}

In the discussion before, we have seen three versions of a VAR(1) model, neither of which was particularly well suited to capture the type of dependence of cycles at quantile level which we observed in Section~\ref{sec:emp}. In the second part of our modelling exercise we now turn our attention to a more flexible class of time series models. Motivated by the quantile autoregression model that was introduced by \cite{koenker2006quantile}, we consider quantile vector autoregression, QVAR, a VAR model with random coefficients:
\begin{equation}\label{eqn:qvar_a}
	Y_{t,j} = \theta_{j0}(U_{t,j}) + \theta_{j1}(U_{t,j}) Y_{t-1,1} + \theta_{j2}(U_{t,j}) Y_{t-1,2}, \quad j=1,2,
\end{equation}
where the $\theta_{ji}$ are coefficient functions and the $U_{t,j}$ are assumed to be independent and uniformly distributed on $[0,1]$.
\cite{ZhuEtAl2018} discuss a model similar to~\eqref{eqn:qvar_a}. Our aim here is to assess whether the time series model~\eqref{eqn:qvar_a} is flexible enough to capture cyclical features in quantiles that were identified in Section~\ref{sec:emp}. To this end, we choose the parameter functions in a data-driven way and then simulate the corresponding quantile coherency to compare with the the nonparameteric estimate. Motivated by the estimation method in \cite{ZhuEtAl2018}, we compute
\begin{equation}\label{eqn:qvar_a_est}
	\hat\theta(\tau) = \arg\min_{\theta(\tau)} \sum_{j=1}^2 \sum_{t = 2}^n \rho_{\tau} \big( Y_{t,j} - \theta_{j0}(\tau) - \theta_{j1}(\tau) Y_{t-1,1} - \theta_{j2}(\tau) Y_{t-1,2} \big),
\end{equation}
$\tau \in \mathcal{T} := \{1/50, 2/50, \ldots, 48/50, 49/50\}$, where $\rho_{\tau}(u) := u(\tau-I\{u < \tau\})$ is the check function (cf. \cite{Koenker2005}). For $\tau \notin \mathcal{T}$ we define $\hat\theta(\tau) := \hat\theta(\eta)$, $\eta := \arg\min_{\eta \in \mathcal{T}} |\tau-\eta|$ (choose the smaller $\eta$ if there are two). The functions $\hat\theta(\tau) = (\hat\theta_{ji}(\tau))$, obtained from the stock market returns, are shown in Figure~\ref{fig:qvar_a_coef}. It is interesting to observe that the functions $\hat\theta_{j1}$ and $\hat\theta_{j2}$, are not constant across quantile levels. This possibly indicates that a VAR model is too simple to capture the complicated dynamics present in the stock markets returns. The ``shock'' at time $t$ to the $j$th equation is delivered by $\hat\theta_{j0}(U_{tj})$.

\cite{koenker2006quantile} and \cite{ZhuEtAl2018} establish conditions that ensure that quantile regressions, similar to~\eqref{eqn:qvar_a_est}, can be used to consistently estimate the parameter functions of the models in their papers. In particular, their model-defining equations (corresponding to \eqref{eqn:qvar_a} in our model) are assumed to be monotonically increasing in~$U_{t,j}$.
The monotonicity condition further implies a particularly convenient form for the conditional quantile function of $Y_{t,j}$ given $Y_{t-1,1}, Y_{t-1,2}$. \cite{FanFan2006} argue that the quantile regression estimate considered by \cite{koenker2006quantile} will be a consistent estimate for
the argument of the minimum of a population version of the loss function, under some mild conditions. For $\hat\theta(\tau)$, defined in~\eqref{eqn:qvar_a_est}, this corresponds to being a consistent estimator for
\begin{equation*}
	\theta^*(\tau) = \arg\min_{\theta(\tau)} \sum_{j=1}^2 \E \rho_{\tau} \big( Y_{t,j} - \theta_{j0}(\tau) - \theta_{j1}(\tau) Y_{t-1,1} - \theta_{j2}(\tau) Y_{t-1,2} \big).
\end{equation*}
\cite{FanFan2006} point out that additional conditions, such as the monotonicity condition, are necessary for $\theta^*(\tau)$ and $\theta(\tau)$ to coincide. These important arguments have to be taken into account when interpreting $\hat\theta(\tau)$ as an estimator for $\theta(\tau)$. Of course, data can always be generated according to equation~\eqref{eqn:qvar_a} where we substitute $\hat\theta(\tau)$ for~$\theta(\tau)$. To assess whether the class of QVAR models is rich enough to reflect cyclical features in the quantiles as we have seen in the data in Section~\ref{sec:emp} it is sufficient to consider individual models from the class. For the purpose of this section, we select a QVAR model of the kind defined in~\eqref{eqn:qvar_a}, in a data-driven way, to then compare the implied quantile coherency with the one estimated non-parametrically in Section~\ref{sec:emp}.

\begin{figure}
  \begin{center}
			\includegraphics[width=0.32\textwidth]{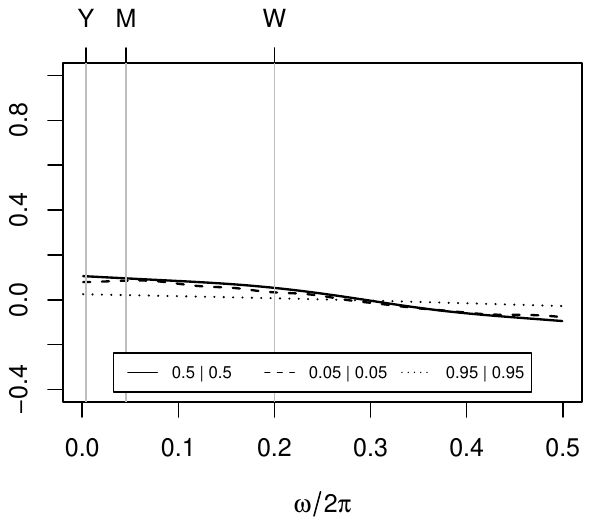}
			\includegraphics[width=0.32\textwidth]{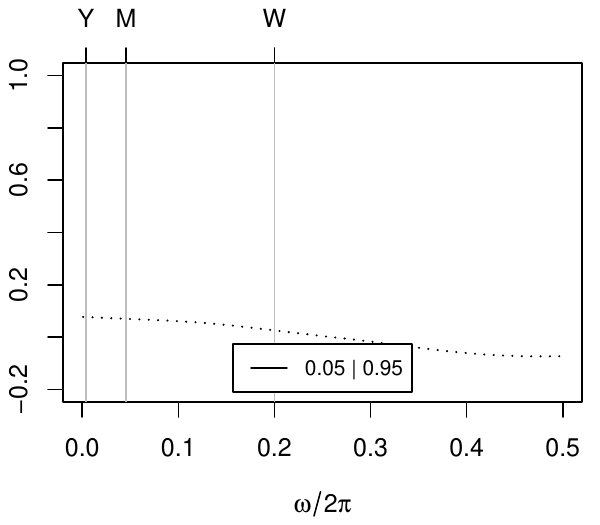}
			\includegraphics[width=0.32\textwidth]{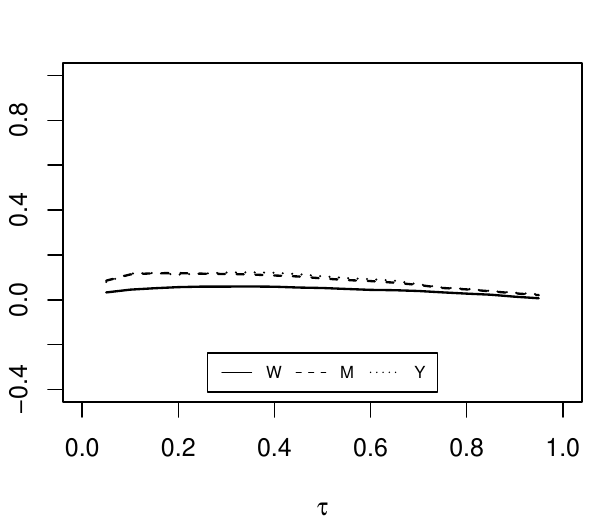}		

			\includegraphics[width=0.32\textwidth]{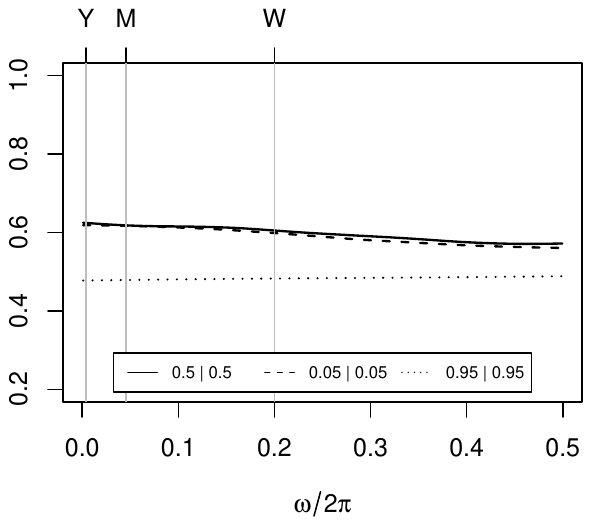}
			\includegraphics[width=0.32\textwidth]{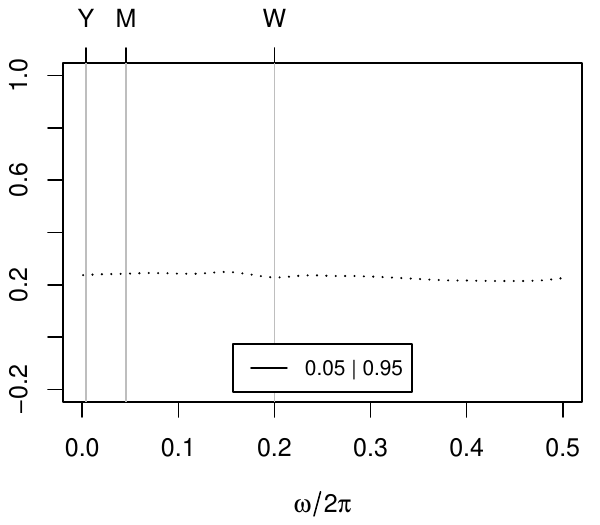}
			\includegraphics[width=0.32\textwidth]{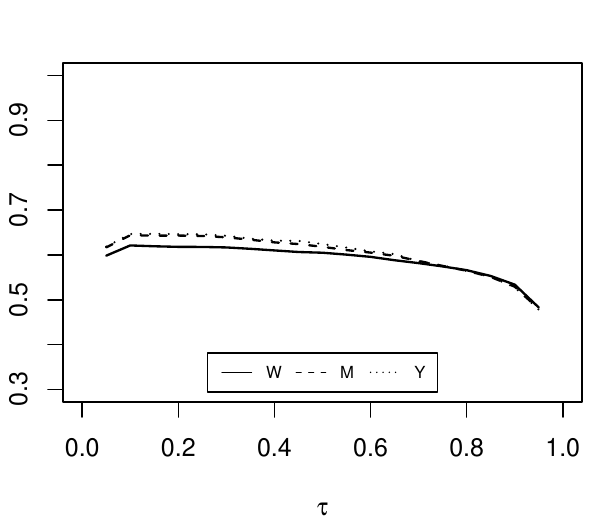}
   \end{center}
  \vspace*{-0.5cm}
	\caption{\label{fig:qc_models_2}Quantile coherency simulated from several QVAR models.}
\end{figure}

In the top row of Figure~\ref{fig:qc_models_2} the quantile coherencies associated with model~\eqref{eqn:qvar_a} where $\hat\theta(\tau)$ was substituted for $\theta(\tau)$ are shown. The plots are of the same format as the ones we had considered before. Strikingly, we observe that the quantile coherency of the fitted model is substantially lower than what we see via the nonparametric estimate in Figure~\ref{fig:emp}. Besides this, in the top row of Figure~\ref{fig:qc_models_2}, we see that the general shape, decreasing lines with frequency, and ordering ($0.95 | 0.95$ shows less dependence than $0.05 | 0.05$) resembles the nonparametric estimate more closely.

\begin{figure}
  \begin{center}
		\includegraphics[width=\textwidth]{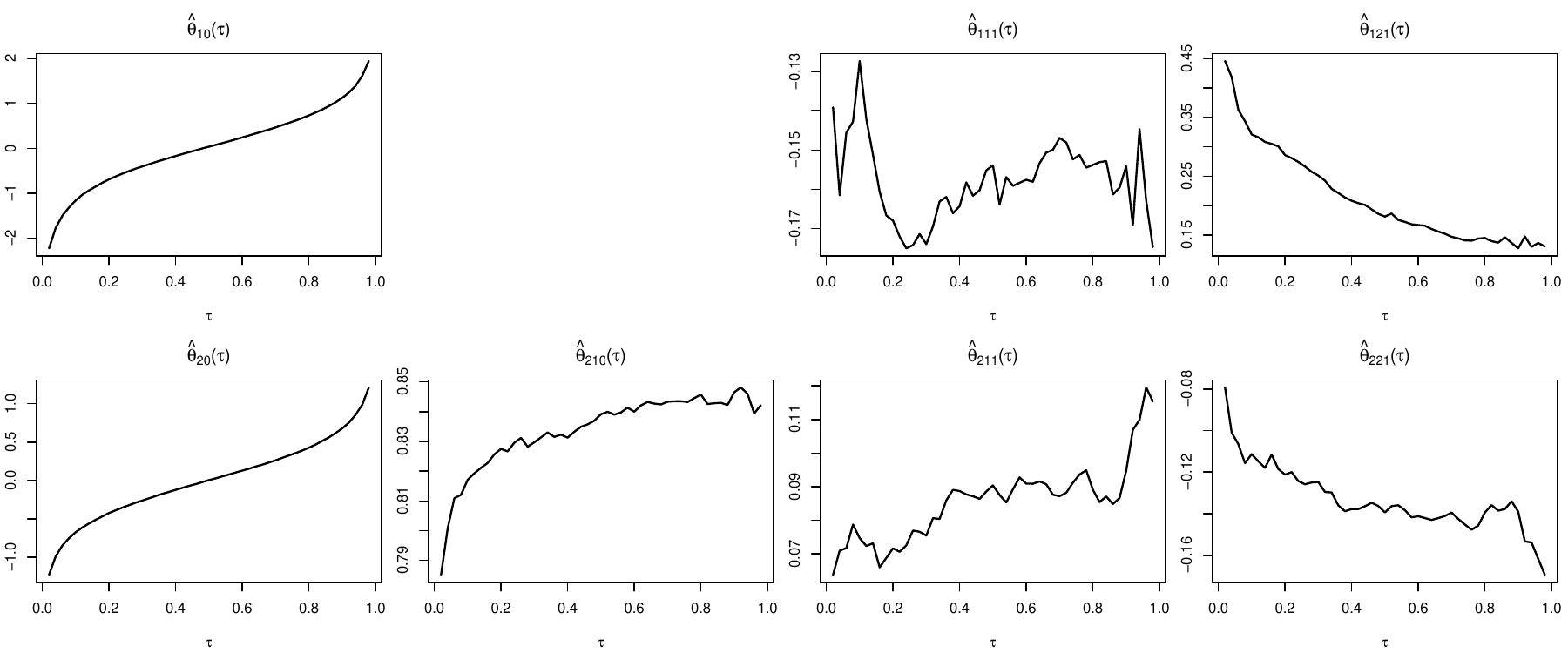}
   \end{center}
	\vspace*{-0.5cm}
  \caption{\label{fig:qvar_coef_2}Parameter functions for model~\eqref{eqn:qvar_b}.}
\end{figure}

Finally, we propose to extend the QVAR(1) stated in~\eqref{eqn:qvar_a}, by adding spatial dependence. More precisely, the model we now consider is
\begin{equation}\label{eqn:qvar_b}
\begin{split}
	Y_{t,1} & = \theta_{10}(U_{t,1}) + \theta_{111}(U_{t,1}) Y_{t-1,1} + \theta_{121}(U_{t,1}) Y_{t-1,2}, \\
	Y_{t,2} & = \theta_{20}(U_{t,2}) + \theta_{211}(U_{t,2}) Y_{t-1,1} + \theta_{221}(U_{t,2}) Y_{t-1,2} + \theta_{210}(U_{t,2}) Y_{t,1}.
\end{split}
\end{equation}
For this model, we compute quantile regression estimates
\begin{multline*}
	\hat\theta(\tau) = \arg\min_{\theta(\tau)} \Big( \sum_{t = 2}^n \rho_{\tau} \big( Y_{t,1} - \theta_{10}(\tau) - \theta_{111}(\tau) Y_{t-1,1} - \theta_{121}(\tau) Y_{t-1,2} \big) \\
	\qquad + \sum_{t = 2}^n \rho_{\tau} \big( Y_{t,2} - \theta_{20}(\tau) - \theta_{210}(\tau) Y_{t,1} - \theta_{211}(\tau) Y_{t-1,1} - \theta_{221}(\tau) Y_{t-1,2} \big) \Big).
\end{multline*}
The estimates obtained from the stock returns data, that also should be cautiously interpreted, are depicted in Figure~\ref{fig:qvar_coef_2}. Note that, if we substitute $Y_{1,t}$ in the second equation of~\eqref{eqn:qvar_b} by the expression given in the first equation, then we see that the ``shocks'' in this model are now dependent, as they are of the form
$(\hat\theta_{10}(U_{t,1}), \hat\theta_{20}(U_{t,2}) + \hat\theta_{210}(U_{t,2}) \hat\theta_{10}(U_{t,1}))$. The parameter function $\hat\theta_{210}$ moderates the strength of dependence. We now again look at the quantile coherency, depicted in the bottom row of Figure~\ref{fig:qc_models_2} and see that the quantile coherencies resemble the nonparameter estimates more closely (in shape, order and magnitude). This is true in particular for the right plot, where the frequency corresponding to the weekly, monthly, and yearly cycles are shown, which could be especially interesting for applied researchers.

In this section we illustrated how quantile coherency can be used by applied researchers to assess time series models regarding their capabilities to capture dependence between general cycles of stock market returns. We have seen that Gaussian VAR models are completely incapable of capturing asymmetries in the dependence of cycles. Our modelling exercise showed how non-Gaussian VAR models can possibly remedy this by allowing more general copulas for the errors in the model. Going further, we have also inspected bivariate quantile autoregression models and seen that their flexibility does better in capturing the general dependence between cycles that we have discovered using quantile coherency in Section~\ref{sec:emp}.

\section{Conclusion}

In this paper we introduced quantile cross-spectral analysis of economic time series providing an entirely model-free, nonparametric theory for the estimation of general cross-dependence structures emerging from quantiles of the joint distribution in the frequency domain. We argue that complex dynamics in time series often arise naturally in many macroeconomic and financial time series, as infrequent periods of large negative values (lower quantiles of the joint distribution) 
may be more dependent than infrequent periods of large positive values (upper quantiles of the joint distribution). Moreover, the dependence may differ in the long-, medium, or short-run. Quantile cross-spectral analysis hence may fundamentally change the way how we view the dependence between economic time series, and may be viewed as a precursor to the subsequent developments in economic research underlying many new modelling strategies.

While connecting two branches of the literature which focus on the dependence between variables in quantiles of their joint distribution and across frequencies separately, the proposed methods may be viewed as an important step in robustifying the traditional cross-spectral analysis as well. Quantile-based spectral quantities are very attractive as they do not require the existence of moments, an important relaxation to the classical assumptions, where moments up to the order of the cumulants involved are typically assumed to exist. The proposed quantities are robust to many common violations of traditional assumptions found in data, including outliers, heavy tails, and changes in higher moments of the distribution. By considering quantiles instead of moments the proposed methods are able to reveal the dependence that remained invisible to the traditional toolsets. As an essential ingredient for a successful applications we have provided a rigorous analysis of the asymptotic properties of the introduced estimators and showed that for a general class of nonlinear processes, properly centred and smoothed versions of the quantile-based estimators converge to centred Gaussian processes.

In an empirical application, we have shown that classical asset pricing theories may not suit the data well, as commonly documented by researchers, because rich dependence structures exists varying across quantiles and frequencies in the joint distribution of returns. We document strong dependence of the bivariate returns series in periods of large negative returns, while positive returns display less dependence over all frequencies. This result is not favourable for an investor, as exactly the opposite would be desired: choosing to invest to stocks with independent negative returns, but dependent positive returns. Our tool reveals that systematic risk originates more strongly from lower quantiles of the joint distribution in the long-, and medium-run investment horizons in comparison to the upper quantiles. In a modelling exercise, we have illustrated how quantile coherency can be employed in the inspection of time series models and might help to find  a model that is capable of capturing the dependencies of cycles of quantile-related features which we had previously revealed in our empirical application.

We believe that our work might open up many exciting new routes for future theoretical as well as empirical research. From the perspective of applications, exploratory analysis based on the quantile cross-spectral estimators can reveal new implications for improvement or even restating of many economic problems. Dependence in many economic time series is of a non-Gaussian nature, calling for an escape from covariance-based methods and allowing for a detailed analysis of the dependence in the quantiles of the joint distribution.

\section*{Acknowledgements}
			
			Authors are listed in alphabetical order, as they have equally contributed to the project. The authors are grateful to Piotr Fryzlewicz, Roger Koenker, Oliver Linton, Stanislav Volgushev, and participants in various seminars and conferences for their comments. In particular, we would like to thank the editor Dennis Kristensen and two anonymous referees, whose reports helped to improve the paper.
			
			Jozef Barun\'{i}k gratefully acknowledges support from the Czech Science Foundation under the GA16-14179S project. Tobias Kley is grateful for being partially supported by the EPSRC fellowship ``New challenges in time series analysis'' (EP/L014246/1) and by the Collaborative Research Center ``Statistical modelling of non-linear dynamic processes'' (SFB~823, Teilprojekt~C1) of the German Research Foundation (DFG).
       
			 For estimation and inference of the quantile cross-spectral measures introduced in this paper, the \textsf{R} package \texttt{quantspec} is provided; cf. \cite{Kley2014JSS}. The \textsf{R} package is available on \url{https://cran.r-project.org/web/packages/quantspec/index.html}

\bibliography{quantile_cross}

\end{document}


\bigskip

\section{Further quantities related to the quantile cross-spectral density kernel}
\label{sec:additionalquantities}

In the situation described in this paper, there exists a right continuous orthogonal increment process $\{Z_j^{\tau}(\omega) : -\pi \leq \omega \leq \pi\}$, for every $j \in \{1,\ldots,d\}$ and $\tau \in [0,1]$, such that the Cram\'{e}r representation
\[I\{X_{t,j} \leq q_{j}(\tau)\} = \int_{-\pi}^{\pi} \ee^{\ii t \omega} {\rm d} Z_j^{\tau}(\omega)\]
holds [cf., e.\,g., ~Theorem~1.2.15 in \cite{TaniguchiKakizawa2000}]. Note the fact that $(X_{t,j})_{t \in \IZ}$ is strictly stationary and therefore $(I\{X_{t,j} \leq q_{j}(\tau)\})_{t \in \IZ}$ is second-order stationary, as the boundedness of the indicator functions implies existence of their second moments.

The quantile cross-spectral density kernels are closely related to these orthogonal increment processes [cf.~\cite[p.\,101]{Brillinger1975} and~\cite[p.\,436]{BrockwellDavis1987}]. More specifically, for $-\pi \leq \omega_1 \leq \omega_2 \leq \pi$, the following relation holds:
\[\int_{\omega_1}^{\omega_2} \mathfrak{f}^{j_1, j_2}(\omega; \tau_1, \tau_2) {\rm d}\omega = \Cov \big( Z_{j_1}^{\tau_1}(\omega_2) - Z_{j_1}^{\tau_1}(\omega_1), Z_{j_2}^{\tau_2}(\omega_2) - Z_{j_2}^{\tau_2}(\omega_1) \big),\]
or shortly: $\mathfrak{f}^{j_1, j_2}(\omega; \tau_1, \tau_2) = \Cov({\rm d}Z_{j_1}^{\tau_1}(\omega), {\rm d}Z_{j_2}^{\tau_2}(\omega))$.
It is important to observe that $\mathfrak{f}^{j_1, j_2}(\omega; \tau_1, \tau_2)$ is complex-valued. One way to represent $\mathfrak{f}^{j_1, j_2}(\omega; \tau_1, \tau_2)$ is to decompose it into its real and imaginary part. The real part is known as the \emph{cospectrum} (of the processes $(I\{X_{t,j_1} \leq q_{j_1}(\tau_1)\})_{t \in \IZ}$ and $(I\{X_{t,j_2} \leq q_{j_2}(\tau_2)\})_{t \in \IZ}$). The negative of the imaginary part is commonly referred to as the \emph{quadrature spectrum}. We will refer to these quantities as the \emph{quantile cospectrum} and \emph{quantile quadrature spectrum} of $(X_{t,j_1})_{t \in \IZ}$ and $(X_{t,j_2})_{t \in \IZ}$. Occasionally, to emphasise that these spectra are functions of $(\tau_1, \tau_2)$, we will refer to them as the \emph{quantile cospectrum kernel} and \emph{quantile quadrature spectrum kernel}, respectively. The quantile quadrature spectrum vanishes if $j_1 = j_2$ and $\tau_1 = \tau_2$. More generally, as described in \cite{kley2014}, for any fixed $j_1, j_2$, the quadrature spectrum will vanish, for all $\tau_1, \tau_2$, if and only if $(X_{t-k,j_1}, X_{t, j_2})$ and $(X_{t+k,j_1}, X_{t, j_2})$ possess the same copula, for all~$k$.

An alternative way to look at $\mathfrak{f}^{j_1, j_2}(\omega; \tau_1, \tau_2)$ is by representing it in polar coordinates. The radius $|\mathfrak{f}^{j_1, j_2}(\omega; \tau_1, \tau_2)|$ is then referred to as the \emph{amplitude spectrum} (of the two processes $(I\{X_{t,j_1} \leq q_{j_1}(\tau_1)\})_{t \in \IZ}$ and $(I\{X_{t,j_2} \leq q_{j_2}(\tau_2)\})_{t \in \IZ}$), while the angle $\arg(\mathfrak{f}^{j_1, j_2}(\omega; \tau_1, \tau_2))$ is the so called \emph{phase spectrum}, respectively. We refer to these quantities as the \emph{quantile amplitude spectrum} and the \emph{quantile phase spectrum} of $(X_{t,j_1})_{t \in \IZ}$ and $(X_{t,j_2})_{t \in \IZ}$.
We note that the quantile spectral distribution function $\int_0^{\omega} \mathfrak{f}^{j_1, j_2}(\lambda; \tau_1, \tau_2)) {\rm d}\lambda$ is clearly another way to represent the quantile-based dependence in the frequency domain. Its properties and estimation procedures are currently investigated in a separate research project and therefore not further discussed here.

Note that quantile coherency $\mathfrak{R}^{j_1, j_2}(\omega; \tau_1, \tau_2)$ which we defined in Section~\ref{sec:ccsdk} as a measure for the dynamic dependence of the two processes $(X_{t,j_1})_{t \in \IZ}$ and $(X_{t,j_2})_{t \in \IZ}$ is the correlation between ${\rm d} Z_{j_1}^{\tau_1}(\omega)$ and ${\rm d} Z_{j_2}^{\tau_2}(\omega)$. Its modulus squared $|\mathfrak{R}^{j_1, j_2}(\omega; \tau_1, \tau_2)|^2$ is referred to as the \emph{quantile coherence kernel} of $(X_{t,j_1})_{t \in \IZ}$ and $(X_{t,j_2})_{t \in \IZ}$. A value of $|\mathfrak{R}^{j_1, j_2}(\omega; \tau_1, \tau_2)|$ close to $1$ indicates a strong (linear) relationship between ${\rm d} Z_{j_1}^{\tau_1}(\omega)$ and ${\rm d} Z_{j_2}^{\tau_2}(\omega)$.

\begin{table}[t]
\caption{\label{tab:newQuantities}
Spectral quantities related to $\mathfrak{f}^{j_1, j_2}(\omega; \tau_1, \tau_2)$ .}
\begin{center}
\begin{tabular}{@{}ll@{}} \toprule
Name & Symbol \\ \midrule
quantile cospectrum of $(X_{t,j_1})_{t \in \IZ}$ and $(X_{t,j_2})_{t \in \IZ}$& $\Re \mathfrak{f}^{j_1, j_2}(\omega; \tau_1, \tau_2)$  \\
quantile quadrature spectrum of $(X_{t,j_1})_{t \in \IZ}$ and $(X_{t,j_2})_{t \in \IZ}$& -$\Im \mathfrak{f}^{j_1, j_2}(\omega; \tau_1, \tau_2)$  \\
quantile amplitude spectrum of $(X_{t,j_1})_{t \in \IZ}$ and $(X_{t,j_2})_{t \in \IZ}$ & $|\mathfrak{f}^{j_1, j_2}(\omega; \tau_1, \tau_2)|$     \\
quantile phase spectrum of $(X_{t,j_1})_{t \in \IZ}$ and $(X_{t,j_2})_{t \in \IZ}$ & $\arg(\mathfrak{f}^{j_1, j_2}(\omega; \tau_1, \tau_2))$     \\
quantile coherency of $(X_{t,j_1})_{t \in \IZ}$ and $(X_{t,j_2})_{t \in \IZ}$ & $\mathfrak{R}^{j_1, j_2}(\omega; \tau_1, \tau_2)$ \\
quantile coherence of $(X_{t,j_1})_{t \in \IZ}$ and $(X_{t,j_2})_{t \in \IZ}$ & $|\mathfrak{R}^{j_1, j_2}(\omega; \tau_1, \tau_2)|^2$   \\ \bottomrule
\end{tabular}
\end{center}
	\footnotesize
	\renewcommand{\baselineskip}{11pm}
	\textbf{Note:} The quantile cross-spectral density kernel $\mathfrak{f}^{j_1, j_2}(\omega; \tau_1, \tau_2)$ of $(X_{t,j_1})_{t \in \IZ}$ and $(X_{t,j_2})_{t \in \IZ}$ is defined in~\eqref{eqn:ccsdk}.
\end{table}

For the readers convenience, a list of the quantities and symbols introduced in this section is provided in Table~\ref{tab:newQuantities}.

Estimators for the quantile cospectrum, quantile quadrature spectrum, quantile amplitude spectrum, quantile phase spectrum, and quantile coherence are then naturally given by $\Re \hat G^{j_1, j_2}_{n,R}(\omega; \tau_1, \tau_2)$, $-\Im \hat G^{j_1, j_2}_{n,R}(\omega; \tau_1, \tau_2)$, $|\hat G^{j_1, j_2}_{n,R}(\omega; \tau_1, \tau_2)|$, $\arg(\hat G^{j_1, j_2}_{n,R}(\omega; \tau_1, \tau_2))$,
and $|\hat{\mathfrak{R}}^{j_1, j_2}_{n,R}(\omega; \tau_1, \tau_2)|^2$, respectively.

\section{An example of a process generating quantile dependence across frequencies: QVAR$(p)$}
\label{sec:vqar}

For a better understanding of the dependence structures that we study in this paper, it is illustrative to introduce a process capable of generating them. We focus on generating dependence at different points of the joint distribution, which will vary across frequencies, but stays hidden from classical measures. In other words, we illustrate the intuition of {spuriously independent} variables, a situation when two variables seem to be independent when traditional cross-spectral analysis is used, while they are indeed clearly dependent at different parts of their joint distribution.

We base our example on a multivariate generalisation of the popular quantile autoregression process (QAR) introduced by \cite{koenker2006quantile}. Inspired by vector autoregression processes (VAR), we link multiple QAR processes through their lag structure and refer to the resulting process as a quantile vector autoregression process (QVAR). This provides a natural way of generating rich dependence structure between two random variables in points of their joint distribution and over different frequencies. The autocovariance function of a stationary QVAR($p$) process is that of a fixed parameter VAR($p$) process. This follows from the argument by \cite{knight2006qar}, who concludes that the exclusive use of autocorrelations may thus ``fail to identify structure in the data that is potentially very informative''. We will show how quantile spectral analysis reveals what otherwise may remain invisible.

Let $\bX_t = (X_{t,1}, \ldots, X_{t,d})'$, $t \in \IZ$, be a sequence of random vectors that fulfills
\begin{equation}
\label{eq:QVAR}
\bm{X}_t = \sum_{j=1}^p\bm{\Theta}^{(j)}(\bm{U}_t)\bm{X}_{t-j}+\bm{\theta}^{(0)}(\bm{U}_t),
\end{equation}
where $\bm{\Theta}^{(1)}, \ldots, \bm{\Theta}^{(p)} $ are $d\times d$ matrices of functions, $\bm{\theta}^{(0)}$ is a $d \times 1$ column vector of functions, and $\bm{U}_t = (U_{t,1}, \ldots, U_{t,d})'$, $t \in \IZ$, is a sequence of independent vectors, with components $U_{t,k}$ that are $\mathcal{U}[0,1]$-distributed. We will assume that the elements of the $\ell$th row $\bm{\theta}_{\ell}^{(j)}(u_\ell) = \big( \theta_{\ell,1}^{(j)}(u_\ell), \ldots, \theta_{\ell,d}^{(j)}(u_\ell) \big)$ of $\bm{\Theta}^{(j)}(u_1,\ldots,u_d)=\big(\bm{\theta}_1^{(j)}(u_1)',\ldots, \bm{\theta}_d^{(j)}(u_d)' \big)'$ and that the $\ell$th element $\theta_{\ell}^{(0)}(u_\ell)$ of
$\bm{\theta}^{(0)} = \big(\theta_{1}^{(0)}(u_1),\ldots, \theta_{d}^{(0)}(u_d) \big)'$
only depend on the $\ell$th variable, respectively. Under this assumption we can rewrite~\eqref{eq:QVAR} as
\begin{equation}
\label{eq:QVARb}
X_{t,i} = \sum_{j=1}^p\bm{\theta}_i^{(j)}(U_{t,i})\bm{X}_{t-j}+\theta_i^{(0)}(U_{t,i}), \quad i = 1, \ldots, d.
\end{equation}
If the right hand side of~\eqref{eq:QVARb} is monotonically increasing, then the conditional quantile function of $X_{t,i}$ given $(\bm{X}_{t-1}, \ldots, \bm{X}_{t-p})$ can be represented as
\[Q_{X_{t,i}}(\tau | \bm{X}_{t-1}, \ldots, \bm{X}_{t-p}) = \sum_{j=1}^p\bm{\theta}_i^{(j)}(\tau)\bm{X}_{t-j}+\theta_i^{(0)}(\tau).\]

Note that in this design the $\ell$th component of $\bm{U}_t$ determines the coefficients for the autoregression equation of the $\ell$th component of $\bX_t$.
We refer to the process as a quantile vector autoregression process of order $p$, hence QVAR$(p)$. The class of processes~\eqref{eq:QVAR} without assumptions regarding the parameters $\bm{\Theta}^{(j)}$ is naturally richer. Yet, the interpretation of the parameters in terms of the conditional quantile functions is possibly lost.

In the bivariate case ($d=2$) of order $p=1$, i.e. QVAR$(1)$, \eqref{eq:QVAR} takes the following form:
\[ \Bigg( \begin{matrix} X_{t,1} \\ X_{t,2} \end{matrix} \Bigg)
=
\Bigg( \begin{matrix}
\theta_{11}^{(1)}(U_{t,1}) & \theta_{12}^{(1)}(U_{t,1}) \\
\theta_{21}^{(1)}(U_{t,2}) & \theta_{22}^{(1)}(U_{t,2})
\end{matrix} \Bigg)
\Bigg( \begin{matrix} X_{t-1,1} \\ X_{t-1,2} \end{matrix} \Bigg)
+
\Bigg( \begin{matrix} \theta_{1}^{(0)}(U_{t,1}) \\ \theta_{2}^{(0)}(U_{t,2}) \end{matrix} \Bigg). \]
For the examples we assume that the components $U_{t,1}$ and $U_{t,2}$ are independent and set the components of $\bm{\theta}^{(0)}$ to $\theta_{1}^{(0)}(u) = \theta_{2}^{(0)}(u) = \Phi^{-1}(u)$, $u\in[0,1]$, where $\Phi^{-1}(u)$ denotes the $u$-quantile of the standard normal distribution. Further, we set the diagonal elements of of $\bm{\Theta}^{(1)}$ to zero (i.\,e., $\theta_{11}^{(1)}(u) = \theta_{22}^{(1)}(u) = 0$, $u\in[0,1]$) and the off-diagonal elements to $\theta_{12}^{(1)}(u) =\theta_{21}^{(1)}(u)= 1.2(u-0.5)$, $u\in[0,1]$.
We thus create cross-dependence by linking the two processes with each other through the other ones lagged contributions. Note that this particular choice of parameter functions leads to the existence of a unique, strictly stationary solution; cf. \cite{BougerolPicard1992}. $(X_{t,1})_{t \in \IZ}$ and $(X_{t,2})_{t \in \IZ}$ are uncorrelated. Note that \cite{hafnerlinton2006} discuss that univariate quantile autoregression nests the popular autoregressive conditional heteroskedasticity (ARCH) models in terms of second order properties. Analogously, our QVAR(1) can be seen to nest a multivariate versions of ARCH.

\begin{figure}
  \begin{center}
      \includegraphics[width=\textwidth]{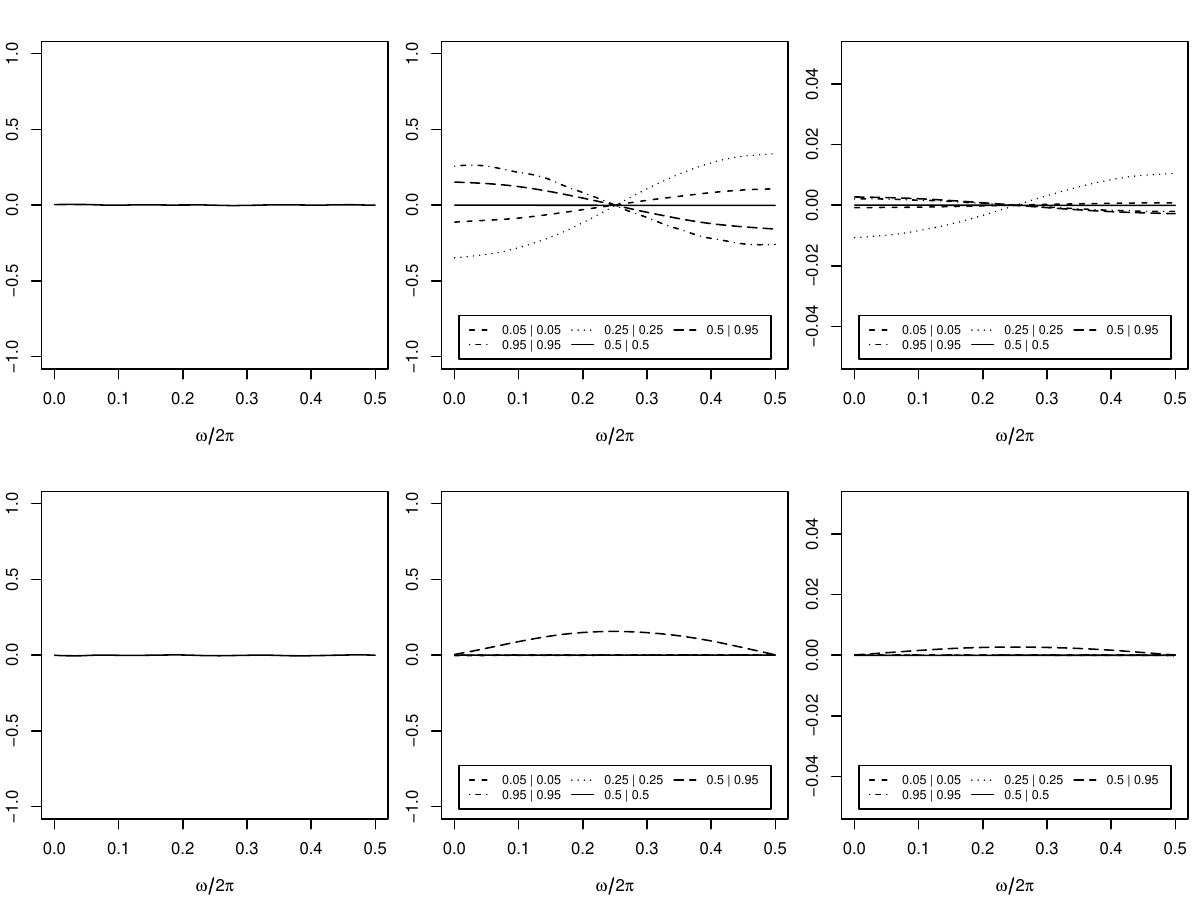}
	\end{center}
	\caption{\label{fig:simQVAR1}Example of dependence structures generated by QVAR(1).
	}
\end{figure}

In Figure~\ref{fig:simQVAR1} the dynamics of the described QVAR(1) process are depicted. In terms of traditional coherency there appears to be no dependence across all frequencies. In terms of quantile coherency, on the other hand, rich dynamics are revealed in the different parts of the joint distribution. While, in the centre of the distribution (at the $0.5|0.5$ level) the dependence is zero across frequencies, we see that the dependence increases if at least one of the quantile levels $(\tau_1, \tau_2)$ is chosen closer to $0$ or $1$. More precisely, we see that the quantile coherency of this QVAR process resembles the shape of an VAR(1) process with coefficient matrix $\bm{\Theta}^{(1)}(\tau_1, \tau_2)$. The two processes are, for example when $\tau_1=0.05$ and $\tau_2=0.95$, clearly positively connected at lower frequencies with exactly the opposite value of quantile coherency at high frequencies, where the processes are in opposition. This also resembles the dynamics of the simple motivating examples from the introductory section of this paper, and highlights the importance of the quantile cross-spectral analysis as the dependence structure stays hidden if only the traditional measures are used. 

\begin{figure}
  \begin{center}
      \includegraphics[width=\textwidth]{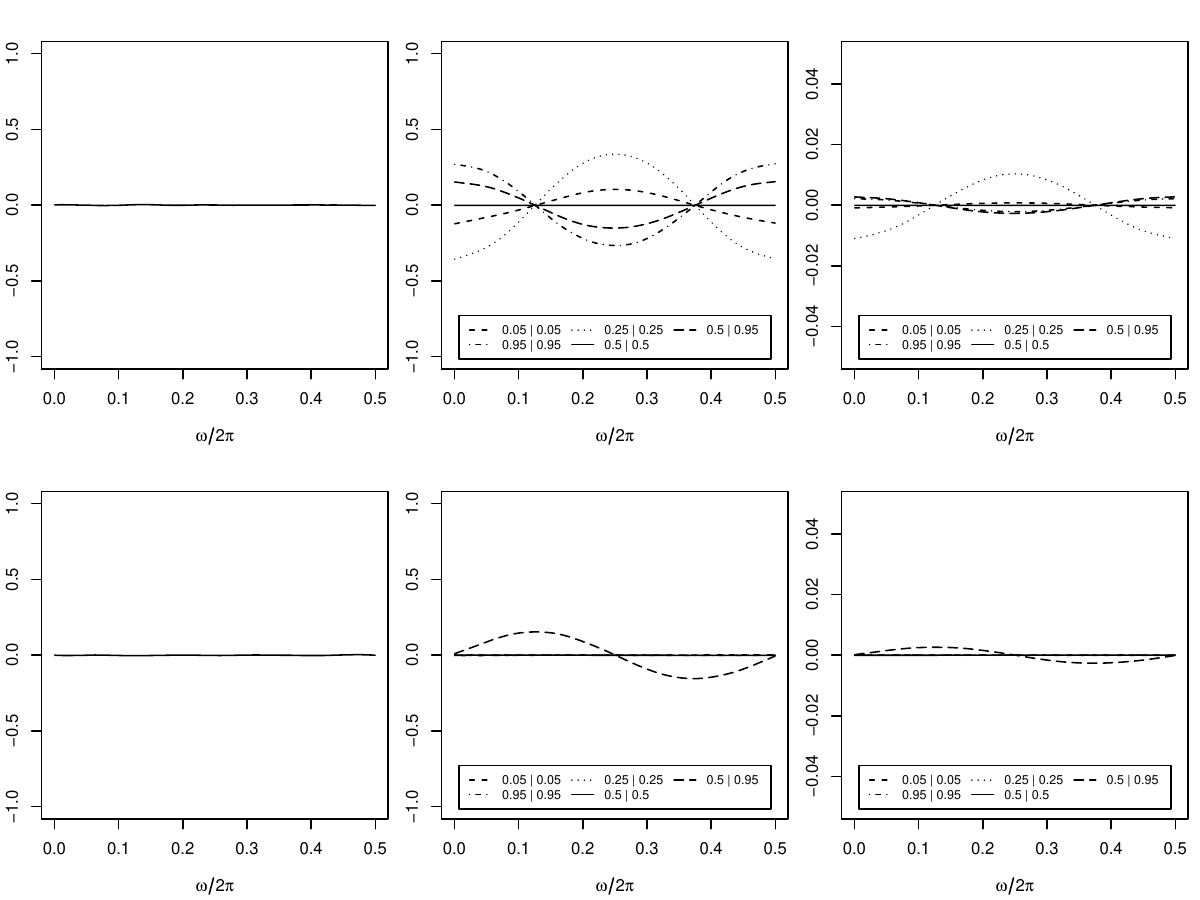} 
\end{center}
\caption{\label{fig:simQVAR2}Example of dependence structures generated by QVAR(2).
}
\end{figure}

\begin{figure}
  \begin{center}
      \includegraphics[width=\textwidth]{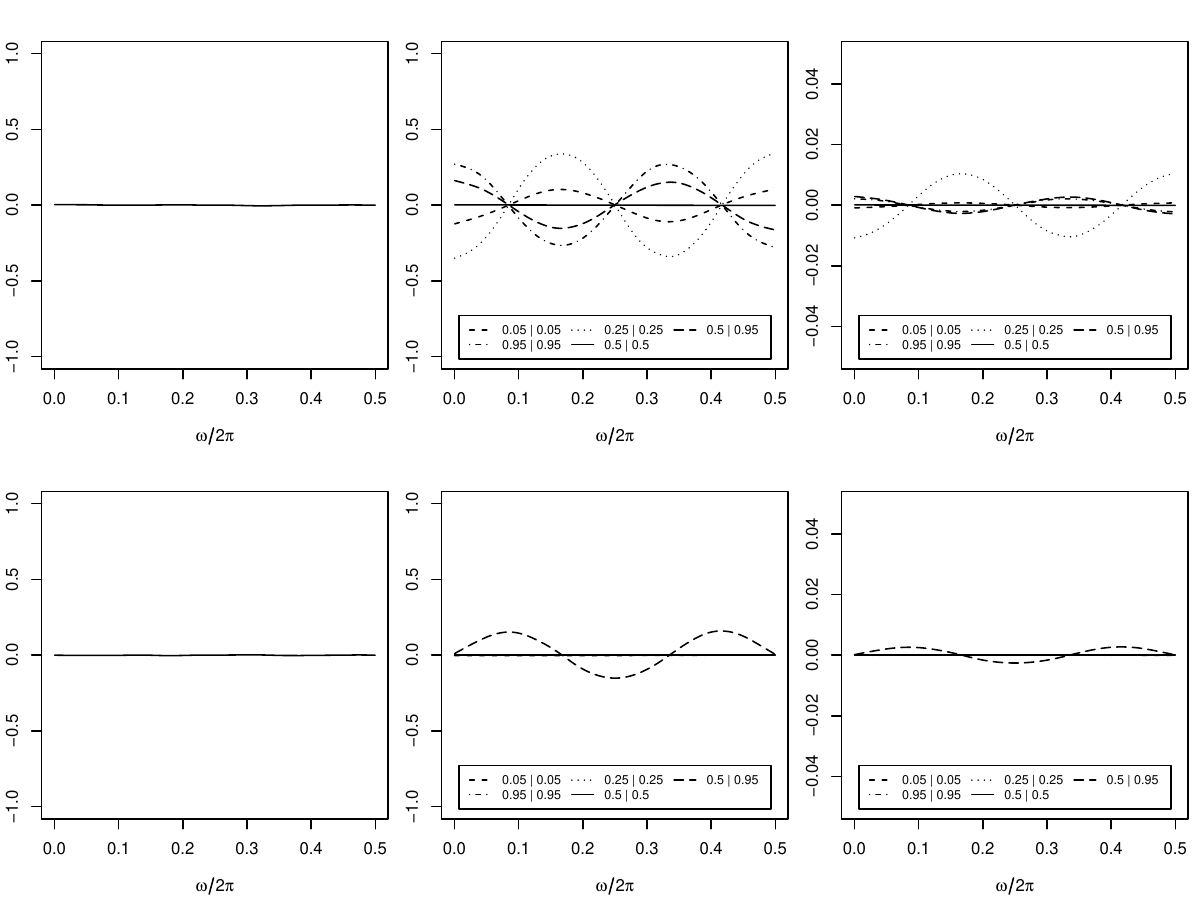}
	\end{center}
	\caption{\label{fig:simQVAR3}Example of dependence structures generated by QVAR(3).
}
\end{figure}

In a second and third example, we consider a similar structure of parameters at the second and third lag. For the QVAR(2) process we let $\theta_{11}^{(j)}(u) = \theta_{22}^{(j)}(u) =0$, for $j=1,2$, $\theta_{12}^{(1)}(u) =\theta_{21}^{(1)}(u)= 0$ and $\theta_{12}^{(2)}(u) =\theta_{21}^{(2)}(u)= 1.2(u-0.5)$. In other words, here, the processes are connected through the second lag of the other one and, again, not directly through their own lagged contributions. In the QVAR(3) process, all coefficients are again set to zero, except for $\theta_{12}^{(3)}(u) =\theta_{21}^{(3)}(u)= 1.2(u-0.5)$, such that the processes are connected only through the third lag of the other component and not through their own contributions.

In Figures~\ref{fig:simQVAR2} and~\ref{fig:simQVAR3} the dynamics of the described QVAR(2) and QVAR(3) processes are shown. Connecting the quantiles of the two processes through the second and third lag gives us richer dependence structures across frequencies. They, again, resemble the shape of the traditional coherencies of VAR(2) and VAR(3) processes.
When traditional coherency is used for the QVAR(2) and QVAR(3) processes, the dependence structure stays completely hidden. 

These examples of the general QVAR$(p)$ specified in~\eqref{eq:QVAR} served to show how rich dependence structures can be created across points of the joint distribution and different frequencies. It is obvious, how more complicated structures for the coefficient functions would lead to even richer dynamics than in the examples shown.

\section{Relation between quantile and traditional spectral quantities in the case of Gaussian processes}
\label{sec:interpretationgauss}

When applying the proposed quantities, it is important to proceed with care when relating them to the traditional correlation and coherency measures. In this section we examine the case of a weakly stationary, multivariate process, where the proposed, quantile-based quantities and their traditional counterparts are directly related. The aim of the discussion is twofold. On one hand it provides assistance in how to interpret the quantile spectral quantities when the model is known to be Gaussian. On the other hand, and more importantly, it provides additional insight in how the traditional quantities break down when the serial dependency structure is not completely specified by the second moments.

We start by the discussion of the general case, where the process under consideration is assumed to be stationary, but needs not to be Gaussian. We will state conditions under which the traditional spectra (i.\,e., the matrix of spectral densities and cross-spectral densities) uniquely determines the quantile spectra (i.\,e., the matrix of quantile spectral densities and cross-spectral densities). In the end of this section we will discuss three examples of bivariate, stationary Gaussian processes and explain how the traditional coherency and the quantile coherency are related.

Denote by $\bm{c} := \{c^{j_1,j_2}_k : j_1,j_2 \in \{1,\ldots,d\}, k \in \IZ\}$. $c^{j_1,j_2}_k := \Cov(X_{t+k,j_1}, X_{t,j_2})$, the family of auto- and cross-covariances. We will also refer to them as the \emph{second moment features} of the process. We assume that $(|c^{j_1,j_2}_k|)_{k \in \IZ}$ is summable, such that the traditional spectra $f^{j_1,j_2}(\omega) := (2 \pi)^{-1} \sum_{k \in \IZ} c^{j_1,j_2}_k {\rm e}^{-{\rm i} k \omega}$
exist. Because of the relation $c^{j_1,j_2}_k = \int_{-\pi}^{\pi} f^{j_1,j_2}(\omega) {\rm e}^{{\rm i} k \omega} {\rm d}\omega$ we will equivalently refer to $\bm{f}(\omega) := (f^{j_1,j_2}(\omega))_{j_1,j_2=1,\ldots,d}$ as the second moment features of the process.

We now state conditions under which the traditional spectra uniquely determine the quantile spectra. Assume that the marginal distribution of $X_{t,j}$ ($j \in \{1,\ldots,d\}$), which we denote by $F_j$, does not depend on $t$ and is continuous. Further, the joint distribution of $\big(F_{j_1}(X_{t+k,j_1}), F_{j_2}(X_{t,j_2}) \big)$, $j_1,j_2 \in \{1,\ldots,d\}$, i.\,e. the copula of the pair $(X_{t+k,j_1}, X_{t,j_2})$, shall depend only on $k$, but not on $t$, and be uniquely specified by the second moment features of the process. More precisely, we assume the existence of functions $C_k^{j_1,j_2}$, such that
\[C_k^{j_1,j_2} \big(\tau_1, \tau_2; \bm{c} \big) = \IP\big( F_{j_1}(X_{t+k,j_1}) \leq \tau_1, F_{j_2}(X_{t,j_2}) \leq \tau_2 \big).\]
Obviously, $\mathfrak{f}^{j_1, j_2}(\omega; \tau_1, \tau_2)$ is then, if it exists, uniquely determined by $\bm{c}$ [note~\eqref{eqn:ccsdk} and the fact that $\gamma_k^{j_1, j_2}(\tau_1, \tau_2) = C_k^{j_1,j_2} \big(\tau_1, \tau_2; \bm{c} \big) - \tau_1 \tau_2$].

In the case of stationary Gaussian processes the assumptions sufficient for the quantile spectra to be uniquely identified by the traditional spectra hold with
\[C_k^{j_1,j_2} \big(\tau_1, \tau_2; \bm{c} \big) := C^{\rm Gauss}(\tau_1, \tau_2; c^{j_1,j_2}_k (c^{j_1,j_1}_0 c^{j_2,j_2}_0)^{-1/2}),\]
where we have denoted the Gaussian copula by $C^{\rm Gauss}(\tau_1, \tau_2; \rho)$.

The converse can be stated under less restrictive conditions. If the marginal distributions are both known and both possess second moments, then the quantile spectra uniquely determine the traditional spectra.

Assume now the previously described situation in which the second moment features $\bm{f}$ uniquely determine the quantile spectra, which we denote by $\mathfrak{f}_{\bm{f}}^{j_1, j_2}(\omega; \tau_1, \tau_2)$ to stress the fact that it is determined by $\bm{f}$.
Thus, the relation between the traditional spectra and the quantile spectra is 1-to-1. Denote the traditional coherency by $R^{j_1,j_2}(\omega) := f^{j_1,j_2}(\omega) / (f^{j_1,j_1}(\omega) f^{j_2,j_2}(\omega))^{1/2}$ and observe that it is also uniquely determined by the second moment features $\bm{f}$. Because the quantile coherency is determined by the quantile spectra which is related to the second moment features $\bm{f}$, as previously explained, we have established the relation of the traditional coherency and the quantile coherency. Obviously, this relation is not necessarily 1-to-1 anymore.

If the stationary process is from a parametric family of time series models the second moment features can be determined for each parameter. We now discuss three examples of Gaussian processes. Each example will have more complex serial dependence than the previous one. Without loss of generality we consider only bivariate examples. The first example is the one of non-degenerate Gaussian white noise. More precisely, we consider a Gaussian process $(X_{t,1}, X_{t,2})_{t \in \IZ}$, where $\Cov(X_{t,i}, X_{s,j}) = 0$ and $\Var(X_{t,i}) > 0$, for all $t \neq s$ and $i, j \in \{1,2\}$.

Observe that, due to the independence of $(X_{t,1}, X_{t,2})$ and $(X_{s,1}, X_{s,2})$, $t \neq s$, we have $\gamma_k^{1,2}(\tau_1, \tau_2) = 0$ for all $k \neq 0$ and $\tau_1, \tau_2 \in [0,1]$. It is easy to see that
\begin{equation}\label{eqn:relQCohTCohIID}
\mathfrak{R}^{1, 2}(\omega; \tau_1, \tau_2)
= \frac{C^{\rm Gauss}(\tau_1, \tau_2; R^{1, 2}(\omega)) - \tau_1 \tau_2}{\sqrt{\tau_1(1-\tau_1)}\sqrt{\tau_2(1-\tau_2))}}
\end{equation}
where $R^{1, 2}(\omega)$ denotes the traditional coherency, which in this case (a bivariate i.\,i.\,d. sequence) equals $c^{1,2}_0 (c^{1,1}_0 c^{2,2}_0)^{-1/2}$ (for all $\omega$).

\begin{figure}
  \begin{center}
    \includegraphics[width=1 \textwidth]{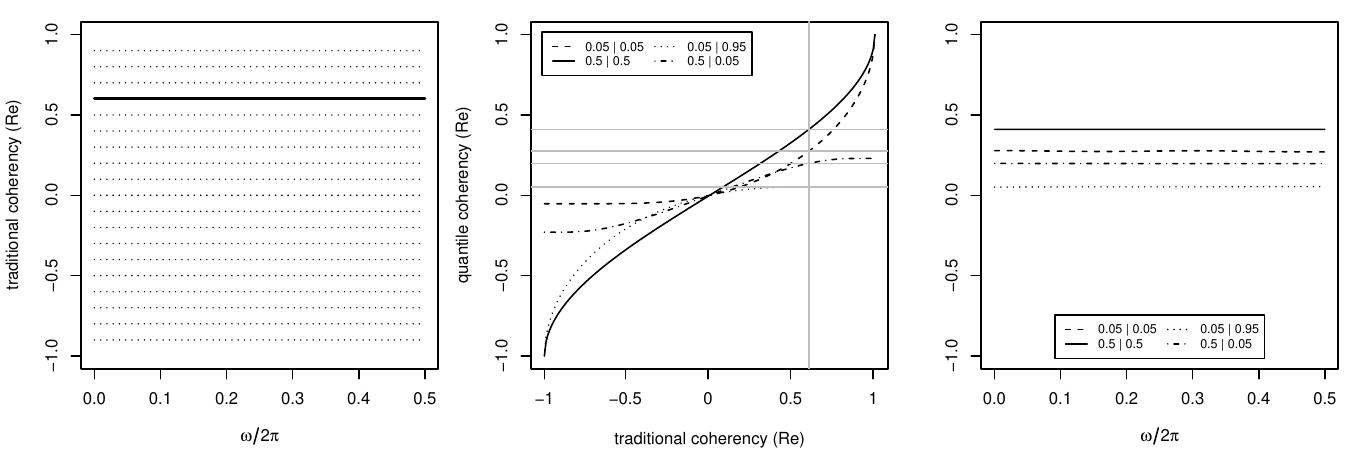} \\
    \vspace{1em}
		\includegraphics[width=1 \textwidth]{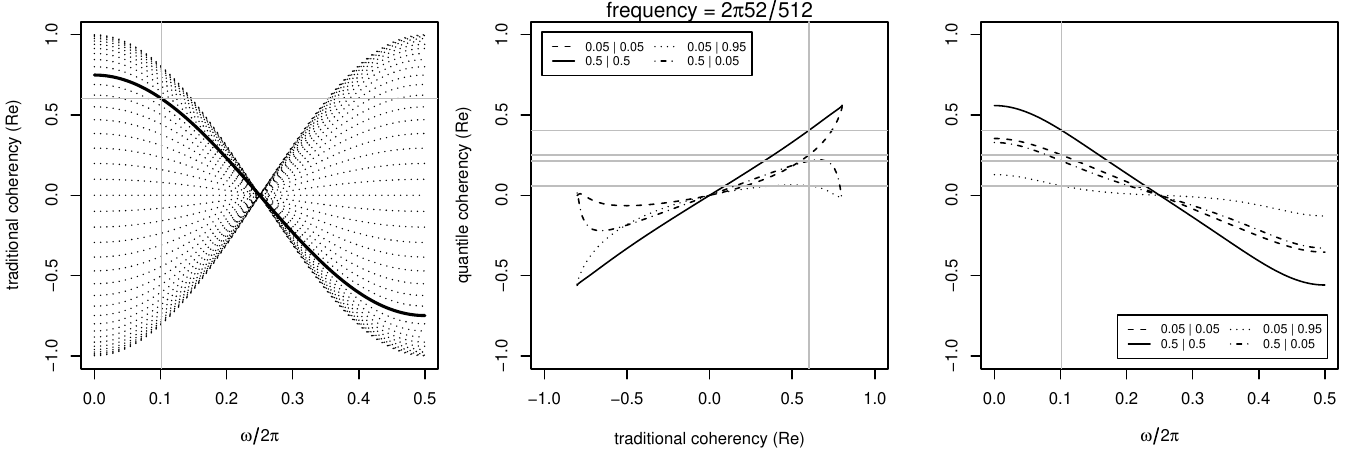}\\
		\vspace{1em}
		\includegraphics[width=1 \textwidth]{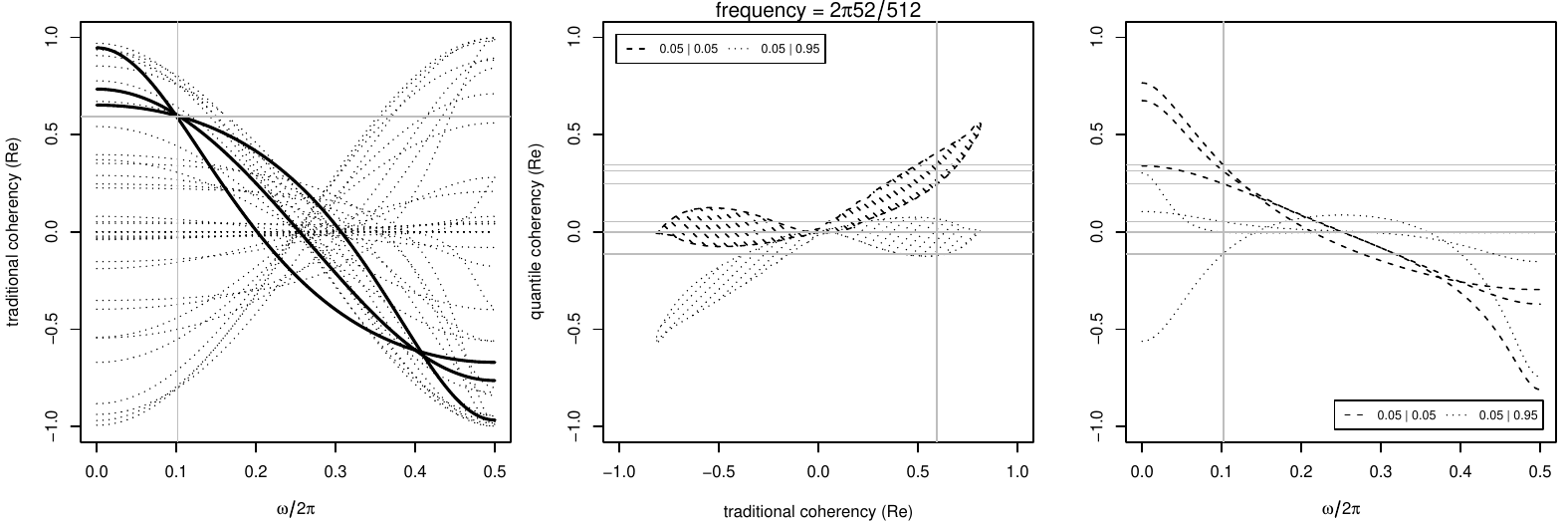}
	\end{center}
\caption{\label{fig:Qcorr}Quantile and traditional coherency for selected Gaussian processes.
}
\end{figure}

By employing~\eqref{eqn:relQCohTCohIID}, we can thus determine the quantile coherency for any given traditional coherency and fixed combination of $\tau_1, \tau_2 \in (0,1)$. In the top-centre part of Figure~\ref{fig:Qcorr} this conversion is visualised for four pairs of quantile levels and any possible traditional coherency. It is important to observe the limited range of the quantile coherency. For example, there never is strong positive dependence between the $\tau_1$-quantile in the first component and the $\tau_2$-quantile in the second component when both $\tau_1$ and $\tau_2$ are close to 0. Similarly, there never is strong negative dependence when one of the quantile levels is chosen close to 0 while the other one is chosen close to 1. This observation is not special for the Gaussian case, but holds for any sequence of pairwise independent bivariate random variables. Bounds that correspond to the case of perfect positive or perfect negative dependence (at the level of quantiles), can be derived from the Fr\'{e}chet/Hoeffding bounds for copulas: in the case of serial independence quantile coherency is bounded by
\[
\frac{\max\{\tau_1 + \tau_2 - 1, 0\} - \tau_1 \tau_2}{\sqrt{\tau_1(1-\tau_1)}\sqrt{\tau_2(1-\tau_2))}}
\leq \mathfrak{R}^{1, 2}(\omega; \tau_1, \tau_2)
\leq \frac{\min\{\tau_1, \tau_2\} - \tau_1 \tau_2}{\sqrt{\tau_1(1-\tau_1)}\sqrt{\tau_2(1-\tau_2))}}.\]
Note that these bounds hold for \emph{any} joint distribution of $(X_{t,i}, X_{t,j})$. In particular, the bound holds independent of the correlation.

In the top-left part of Figure~\ref{fig:Qcorr} traditional coherencies are shown for this example. Because no serial dependence is present, all coherencies are flat lines. Their level is equal to the correlation between the two components. In the top-right part of Figure~\ref{fig:Qcorr} the quantile coherency for the example is shown when the correlation is $0.6$ (the corresponding coherency is marked with a bold line in the top-left figure). Note that for fixed $\tau_1$ and $\tau_2$ the value of the quantile coherency corresponds to the value in the top-centre figure where the vertical grey line and the corresponding graph intersect. The quantile coherency in the right part does not depend on the frequency, because in this example there is no serial dependence.

In the top-centre part of Figure~\ref{fig:Qcorr} it is important to observe that for traditional coherency 0 (i.\,e., when the components are independent, due to $(X_{t,1}, X_{t,2})$ being uncorrelated jointly Gaussian) quantile coherency is zero at all quantile levels.

In the next two examples we stay in the Gaussian framework, but introduce serial dependence. Consider a bivariate, stable VAR(1) process $\bX_t=(X_{t,1}, X_{t,2})'$, $t \in \IZ$, fulfilling the difference equation
\begin{equation}\label{eqn:VAR}
\bX_t = \bA \bX_{t-1} + \bm{\varepsilon}_t,
\end{equation}
with parameter $\bA \in \IR^{2 \times 2}$ and i.\,i.\,d., centred, bivariate, jointly normally distributed innovations $\bm{\varepsilon}_t$ with unit variance $\E (\bm{\varepsilon}_t \bm{\varepsilon}'_t) = \bI_2$.

In our second example serial dependence is introduced, by relating each component to the lagged other component in the regression equation. In other words, we consider model~\eqref{eqn:VAR} where the matrix $\bA$ has diagonal elements equal to 0 and some value $a$ on the off-diagonal. Assuming $|a| < 1$ yields a stable process. As described earlier, the traditional spectral density matrix, which in this example is of the form
\[
	\bm{f}(\omega) := (2\pi)^{-1} \Big(\bI_2 - \begin{pmatrix} 0 & a \\ a & 0 \end{pmatrix} {\rm e}^{-i \omega}\Big)^{-1} \Big(\bI_2 - \begin{pmatrix} 0 & a \\ a & 0 \end{pmatrix} {\rm e}^{i \omega}\Big)^{-1}, \ |a| < 1,
\]
uniquely determines the traditional coherency and, because of the Gaussian innovations, also the quantile coherency.

In the middle-left plot of Figure~\ref{fig:Qcorr} the traditional coherencies for this model are shown when $a$ takes different values. If we now fix a frequency [$\neq \pi/4$], then the value of the traditional coherency for this frequency uniquely determines the value of $a$. In Figure~\ref{fig:Qcorr} we have marked the frequency of $\omega = 2\pi 52/512$ and coherency value of $0.6$ by grey lines and printed the corresponding coherency (as a function of $\omega$) in bold. Note that of the many pictured coherencies [one for each $a \in (-1,1)$] only one has the value of $0.6$ at this frequency. In the centre plot of the middle row we show the relation between the traditional coherency and quantile coherency for the considered model. For four combinations of quantile levels and all values of $a \in (-1,1)$ the corresponding traditional coherencies and quantile coherencies are shown. It is important to observe that the relation is shown only for one frequency [$\omega = 2\pi 52/512$]. We observe that the range of values for the quantile coherency is limited and that the range depends on the combination of quantile levels and on the frequency. While this is quite similar to the first example where \emph{quantile} coherency had to be bounded due to the Fr\'{e}chet/Hoeffding bounds, we here also observe (for this particular model and frequency) that the range of values for the \emph{traditional} coherency is limited. This fact is also apparent in the middle-left plot. To relate the traditional and quantile coherency at this particular frequency, one can, using the centre-middle plot, proceed as in the first example. For a given frequency choose a valid traditional coherency (x-axis of the middle-centre plot) and combination of quantile levels (one of the lines in the plot) and then determine the value for the quantile coherency (depicted in the right plot). Note that (in this example), for a given frequency and combination of quantile levels the relation is still a function of the traditional coherency, but fails to be injective.

In our final example we consider the Gaussian VAR(1) model~\eqref{eqn:VAR} where we now allow for an additional degree of freedom, by letting the matrix $\bA$ be of the form where the diagonal elements both are equal to $b$ and keep the value $a$ on the off-diagonal as before. Thus, compared to the previous example, where $b=0$ was required, each component now may also depend on its own lagged value. It is easy to see that $|a + b| < 1$ yields a stable process. In this case the tradtional spectral density matrix is of the form
\[
	\bm{f}(\omega) := (2\pi)^{-1} \Big(\bI_2 - \begin{pmatrix} b & a \\ a & b \end{pmatrix} {\rm e}^{-i \omega}\Big)^{-1} \Big(\bI_2 - \begin{pmatrix} b & a \\ a & b \end{pmatrix} {\rm e}^{i \omega}\Big)^{-1}, \ |a + b| < 1.
\]

In the bottom-left part of Figure~\ref{fig:Qcorr} a collection of traditional coherencies (as functions of $\omega$) is shown. Due to the extra degree of freedom in the model the variety of shapes increased dramatically. In particular, for a given frequency, the value of the traditional coherency does not uniquely specify the model parameter any more. We have marked three coherencies (as functions of $\omega$) that have value $0.6$ at $\omega = 2\pi 52/512$ in bold to stress this fact. The corresponding processes have (for a fixed combination of quantile levels) different values of quantile coherency at this frequency. This fact can be seen from the bottom-centre part of Figure~\ref{fig:Qcorr}, where the relation between traditional and quantile coherency is depicted for the frequency fixed and two combinations of quantile levels are shown in black and grey. Note the important fact that the relation (for fixed frequency) is \emph{not} a function of the traditional coherency any more. The bottom-right part of the figure shows the quantile coherency curves (as a function of $\omega$) for the three model parameters (shown in bold in the bottom-left part of the figure) and the two combination of quantile levels. It is clearly visible that even though, for the particular fixed frequency, the traditional coherency coincide, the value and shape of the quantile coherency can be very different depending on the underlying process. This third example illustrated how a frequency-by-frequency comparison of the traditional coherency with its quantile-based counterpart may fail, even when the process is quite simple.

We have seen, from the theoretical discussion in the beginning of this section, that for Gaussian processes, when the marginal distributions are fixed, a relation between the traditional spectra and the quantile spectra exists. This relation is a 1-to-1 relation between the quantities \emph{as functions} of frequency (and quantile levels). The three examples have illustrated that a comparison on a frequency-by-frequency basis may be possible in special cases but does not hold in general.

In conclusion we therefore advise to see the quantile cross-spectral density as a measure for dependence on its own, as the quantile-based quantities focus on more general types of dependence. We further point out that quantile coherency may be used in examples where the conditions that make a relation possible are fulfilled, but also, for example, to analyse the dependence in the quantile vector autoregressive (QVAR) processes, described in Section~\ref{sec:vqar}. The QVAR processes possess more complicated dynamics, which cannot be described only by the second order moment features.

\section{Asymptotic properties of the proposed estimators for quantile cross-spectral densities}
\label{sec:asymp2}

We are now going to state a result on the asymptotic properties of the CCR-periodogram $\bI_{n,R}(\omega; \tau_1, \tau_2)$ defined in \eqref{eqn:inr} and~\eqref{eqn:inrMatr}.

\begin{proposition}\label{prop:inr}
Assume that $(\bX_t)_{t\in\IZ}$ is strictly stationary and satisfies Assumption~\ref{ass:exp_alpha_mix}. Further assume that the marginal distribution functions $F_j$, $j=1,\ldots,d$ are continuous. Then, for every fixed $\omega \neq 0 \mod 2\pi$,
\begin{equation}
\label{eq:IntoI}
\Big( \bI_{n, R}(\omega; \tau_1, \tau_2) \Big)_{ (\tau_1,\tau_2) \in [0,1]^2} \Rightarrow
\Big(\mathbb{I}(\omega;\tau_1,\tau_2) \Big)_{ (\tau_1,\tau_2) \in [0,1]^2}
\quad \text{in } \ell_{\mathbb{C}^{d \times d}}^{\infty}([0,1]^2).
\end{equation}
The $\IC^{d \times d}$-valued limiting processes $\mathbb{I}$, indexed by $(\tau_1, \tau_2) \in [0,1]^2$, is of the form
\[
\mathbb{I}(\omega;\tau_1,\tau_2) = \frac{1}{2\pi} \mathbb{D}(\omega;\tau_1) \overline{\mathbb{D}(\omega;\tau_2)'},
\]
where $\mathbb{D}(\omega; \tau) = ( \mathbb{D}^j(\omega; \tau) )_{j = 1,\ldots,d}$, $\tau \in [0,1]$, $\omega \in \IR$ is a centred, $\IC^d$-valued Gaussian processes with covariance structure of the following form
\[
\Cov(\mathbb{D}^{j_1}(\omega; \tau_1), \mathbb{D}^{j_2}(\omega; \tau_2))
= 2 \pi \mathfrak{f}^{j_1, j_2}(\omega; \tau_1, \tau_2).
\]
Moreover,  $\mathbb{D}(\omega; \tau) = \overline{\mathbb{D}(-\omega; \tau)} = \mathbb{D}(\omega+2\pi; \tau)$, and the family~$\{\mathbb{D}(\omega; \, \cdot) \ : \ \omega \in [0,\pi] \}$ is a collection of independent processes. In particular, the weak convergence~(\ref{eq:IntoI}) holds jointly for any finite  fixed collection of frequencies~$\omega$.
\end{proposition}
For $\omega = 0 \mod 2\pi$ the asymptotic behaviour of the CCR-periodogram is as follows: we have $d^j_{n,R}(0; \tau) = n \tau + o_p(n^{1/2})$, where the exact form of the remainder term depends on the number of ties in $X_{j,0}, \ldots, X_{j,n-1}$. Therefore, under the assumptions of Proposition~\ref{prop:inr}, we have $\bI_{n, R}(0; \tau_1, \tau_2) = n(2\pi)^{-1}\tau_1\tau_2 1_d 1'_d + o_p(1)$, where $1_d := (1, \ldots, 1)' \in \IR^d$.

We now state a result that quantifies the uncertainty in estimating $\bfrakf(\omega; \tau_1, \tau_2)$ by $\bG_{n,R}(\omega; \tau_1, \tau_2)$ asymptotically.
\begin{theorem} \label{thm:AsympDensityRankEstimator}
Let Assumptions~\ref{ass:exp_alpha_mix} and~\ref{ass:W} hold. Assume that the marginal distribution functions $F_j$, $j=1,\ldots,d$ are continuous and that constants $\kappa > 0$ and $k \in\IN$ exist, such that
$b_n = o(n^{-1/(2k+1)})$ and $b_n n^{1-\kappa} \rightarrow \infty$.
Then, for any fixed $\omega \in  \IR $, the process
\[
	\mathbb{G}_{n}(\omega; \cdot, \cdot) := \sqrt{n b_n} \Big(\hat \bG_{n,R}(\omega; \tau_1, \tau_2) - \bfrakf(\omega; \tau_1, \tau_2) - \bB_n^{(k)}(\omega; \tau_1, \tau_2)\Big)_{\tau_1,\tau_2 \in [0,1]}
\]
  satisfies
\begin{equation} \label{eq:gntoh}
 \mathbb{G}_{n}(\omega; \cdot,\cdot) \Rightarrow \IH(\omega; \cdot,\cdot) \text{ in $\ell_{\mathbb{C}^{d \times d}}^{\infty}([0,1]^2)$,}
\end{equation}
 where the elements of the bias matrix $\bB_n^{(k)}$ are given by
\begin{equation}\label{def:bias}
\Big\{ \bB_n^{(k)}(\omega; \tau_1, \tau_2) \Big\}_{j_1, j_2}:=
  \sum_{\ell=2}^k \frac{b_n^{\ell}}{\ell!} \int_{ -\pi}^{\pi} v^{\ell} W(v) dv \frac{{\rm d}^{\ell}}{{\rm d}\omega^{\ell}}\mathfrak{f}^{j_1, j_2}(\omega; \tau_1, \tau_2)
\end{equation}
  and $\mathfrak{f}^{j_1, j_2}(\omega; \tau_1, \tau_2)$ is defined in \eqref{eqn:ccsdk}. The process $\IH(\omega;\cdot,\cdot) := (\IH^{j_1, j_2}(\omega;\cdot,\cdot))_{j_1, j_2 = 1,\ldots,d}$ in~\eqref{eq:gntoh}  is a centred, $\IC^{d \times d}$-valued Gaussian process characterised by
\begin{multline} 
\Cov\big(\IH^{j_1, j_2}(\omega; u_1, v_1\big), \IH^{k_1, k_2}(\lambda; u_2, v_2)\big) \\ = 2\pi \Big(\int_{-\pi}^\pi W^2(\alpha){\rm d}\alpha \Big)
\Big( \mathfrak{f}^{j_1, k_1}(\omega; u_1, u_2) \mathfrak{f}^{j_2, k_2}(-\omega; v_1, v_2) \eta(\omega - \lambda) \\
+ \mathfrak{f}^{j_1, k_2}(\omega; u_1, v_2) \mathfrak{f}^{j_2, k_1}(-\omega; v_1, u_2)
\eta(\omega + \lambda) \Big),
\end{multline}
where $\eta(x) := I\{x = 0 (\mod 2\pi)\}$ [cf.~\cite[p.\,148]{Brillinger1975}] is the $2\pi$-periodic extension of Kronecker's delta function.
The family~$\{\IH(\omega; \, \cdot, \cdot),$ $ \omega \in  [0,\pi] \}$ is a collection of independent processes and $\IH(\omega; \tau_1, \tau_2) = \overline{\IH(-\omega; \tau_1, \tau_2)} = \IH(\omega+2\pi; \tau_1, \tau_2)$.
\end{theorem}

A few remarks on the result are in order.
	In sharp contrast to classical spectral analysis, where higher-order moments are required to obtain smoothness of the spectral density [cf.~\cite{Brillinger1975}, p.\,27], Assumption~\ref{ass:exp_alpha_mix} guarantees that the quantile cross-spectral density is an analytical function of $\omega$. Hence, the $k$th derivative of $\omega \mapsto \mathfrak{f}^{j_1, j_2}(\omega; \tau_1, \tau_2)$ in~\eqref{def:bias} exists without further assumptions.

	The case $\omega = 0 \mod 2\pi$ does not require separate treatment as in Proposition~\ref{prop:inr}, because $I_{n,R}^{j_1, j_2}(0, \tau_1, \tau_2)$ is excluded in~\eqref{eqn:DefRankEstimator}: the definition of $\hat G^{j_1, j_2}_{n,R}(\omega; \tau_1, \tau_2)$.
	
	Assume that $W$ is a kernel of order $p$; i.\,e., for some $p$, satisfies $\int_{-\pi}^{\pi} v^j W(v) {\rm d}v = 0$, for all $j < p$, and $0 < \int_{-\pi}^{\pi} v^p W(v) {\rm d}v < \infty$. E.\,g., the Epanechnikov kernel is a kernel of order $p=2$. Then, the bias is of order $b_n^p$. As the variance is of order $(n b_n)^{-1}$, the mean squared error is minimal, if $b_n \asymp n^{-1/(2p+1)}$. This optimal bandwidth fulfills the assumptions of Theorem~\ref{thm:AsympDensityRankEstimator}.
A detailed discussion of how Theorem~\ref{thm:AsympDensityRankEstimator} can be used to construct asymptotically valid confidence intervals is deferred to Section~D.

The independence of the limit $\{\IH(\omega; \, \cdot, \cdot),$ $ \omega \in  [0,\pi] \}$ has two important implications. On one hand, the weak convergence~(\ref{eq:gntoh}) holds jointly for any \emph{finite} fixed collection of frequencies~$\omega$. On the other hand, if one were to consider the smoothed CCR-periodogram as a function of the three arguments $(\omega, \tau_1, \tau_2)$, weak convergence cannot hold any more. This limitation of convergence is due to the fact that there exists no tight element in $\ell_{\mathbb{C}^{d \times d}}^{\infty}([0,\pi] \times [0,1]^2)$ that has the right finite-dimensional distributions, which would be required for process convergence in $\ell_{\mathbb{C}^{d \times d}}^{\infty}([0,\pi] \times [0,1]^2)$.

Fixing $j_1, j_2$ and $\tau_1, \tau_2$ the CCR-periodogram $\hat G^{j_1, j_2}_{n,R}(\omega; \tau_1, \tau_2)$ and traditional smoothed cross-periodogram determined from the unobservable, bivariate time series
\begin{equation} 
\big(I\{F_{j_1}(X_{t,j_1}) \leq \tau_1\},I\{F_{j_1}(X_{t,j_2}) \leq \tau_2\}\big), \quad t = 0, \ldots, n-1,
\end{equation}
are asymptotically equivalent. Theorem~\ref{thm:AsympDensityRankEstimator} thus reveals that in the context of the estimation of the quantile cross-spectral density the estimation of the marginal distribution has no impact on the limit distribution (cf. comment after Remark~3.5 in~\cite{kley2014}).

\section{On the construction of interval Estimators}\label{sec:CI}

In this section we collect details on how to construct pointwise confidence bands.

Sections~\ref{sec:asymp} and~\ref{sec:asymp2} contained asymptotic results on the uncertainty of point estimation of the newly introduced quantile cross-spectral quantities. In this section we describe strategies to estimate the variances (of the real and imaginary parts) that appear in those limit results and describe how asymptotically valid pointwise confidence bands can be constructed.

In all three subsections the following comment is relevant. Assuming that we have determined the weights $W_n$ form a kernel $W$ that is of order $d$. We will choose a bandwidth $b_n = o(n^{-1/(2d+1)})$. This choice implies that compared to the variance the bias (that in some form appears in both limit results) is asymptotically negligible: $\sqrt{n b_n} \bB_n^{(k)}(\omega; \tau_1, \tau_2) = o(1)$.

\subsection[Pointwise confidence bands for f]{Pointwise confidence bands for $\bfrakf$}\label{rem:CI:f}

Utilising Theorem~\ref{thm:AsympDensityRankEstimator} we now construct pointwise asymptotic $(1-\alpha)$-level confidence bands for the real and imaginary parts of $\mathfrak{f}^{j_1, j_2}(\omega_{kn};\tau_1,\tau_2)$, $\omega_{kn} := 2\pi k / n$, as follows:

\[C^{(1)}_{{\rm r},n}(\omega_{kn}; \tau_1, \tau_2) := \Re \tilde{G}^{j_1, j_2}_{n,R}(\omega_{kn}; \tau_1, \tau_2) \pm \Re \sigma_{(1)}^{j_1, j_2}(\omega_{kn}; \tau_1, \tau_2) \Phi^{-1}(1-\alpha/2),\]
for the real part, and 
\[C^{(1)}_{{\rm i},n}(\omega_{kn}; \tau_1, \tau_2) := \Im \tilde{G}^{j_1, j_2}_{n,R}(\omega_{kn}; \tau_1, \tau_2) \pm \Im \sigma_{(1)}^{j_1, j_2}(\omega_{kn}; \tau_1, \tau_2) \Phi^{-1}(1-\alpha/2),\]
for the imaginary part of the quantile cross-spectrum. Here,
\[ \tilde{G}^{j_1, j_2}_{n,R}(\omega_{kn}; \tau_1, \tau_2) :=  \hat{G}^{j_1, j_2}_{n,R}(\omega_{kn}; \tau_1, \tau_2) / W_n^k, \quad
W_n^k := \frac{2\pi}{n} \sum_{s = 1}^{n-1} W_n(\omega_{kn} - \omega_{sn}),\]
and $\Phi$ denotes the cumulative distribution function of the standard normal distribution,%
\footnote{Note that for $k=0,\ldots,n-1$ we have $W_n^k := 2\pi / n \sum_{0 = s \neq k}^{n-1} W_n(2\pi s/n)$. For $k \in \IZ$ with $k < 0$ or $k \geq n$ we can define it as the $n$ periodic extension.}
\begin{equation*}
\big( \Re \sigma^{j_1, j_2}(\omega_{kn}; \tau_1, \tau_2) \big)^2
:= 0 \vee \begin{cases}
			\Cov( \HH{1}{2}, \HH{1}{2} ) & \text{if $j_1 = j_2$ and $\tau_1 = \tau_2$,} \\
      \frac{1}{2} \big( \Cov( \HH{1}{2}, \HH{1}{2} ) + \Re \Cov( \HH{1}{2}, \HH{2}{1} ) \big) & \text{otherwise},
\end{cases}
\end{equation*}
and
\begin{equation*}
\big( \Im \sigma^{j_1, j_2}(\omega_{kn}; \tau_1, \tau_2) \big)^2
:= 0 \vee \begin{cases}
      0 & \text{if $j_1 = j_2$ and $\tau_1 = \tau_2$,}\\
      \frac{1}{2} \big( \Cov( \HH{1}{2}, \HH{1}{2} ) - \Re \Cov( \HH{1}{2}, \HH{2}{1} ) \big) & \text{otherwise},
\end{cases}
\end{equation*}

where $\Cov( \HH{a}{b}, \HH{c}{d} )$ denotes an estimator of
$\Cov\big(\IH^{j_a, j_b}(\omega_{kn}; \tau_a, \tau_b\big), \IH^{j_c, j_d}(\omega_{kn}; \tau_c, \tau_d)\big)$. Here, motivated by Theorem~7.4.3 in~\cite{Brillinger1975}, we use 
\begin{multline}\label{eqn:rem:CI:f:1}
\Big(\frac{2\pi}{n \cdot W_n^k}\Big)
\times \Bigg[ \sum_{s=1}^{n-1} W_n\big(2\pi (k -  s) / n \big) W_n\big(2\pi (k -  s) / n \big) \tilde G^{j_a, j_c}_{n,R}(\tau_a, \tau_c; 2\pi s / n) \tilde G^{j_b, j_d}_{n,R}(\tau_b, \tau_d; - 2\pi s / n) \\
+ \sum_{s=1}^{n-1} W_n\big(2\pi (k -  s) / n \big) W_n\big(2\pi (k +  s) / n \big) \tilde G^{j_a, j_d}_{n,R}(\tau_a, \tau_d; 2\pi s / n)  \tilde G^{j_b, j_c}_{n,R}(\tau_b, \tau_c; -2\pi s / n)   \Bigg]
\end{multline}

The definition of $\sigma_{(1)}^{j_1, j_2}(\omega_{kn}; \tau_1, \tau_2)$ is motivated by the fact that $\Im \hat{G}^{j_1, j_2}_{n,R}(\omega_{kn}; \tau_1, \tau_2) = 0$, if $j_1 = j_2$ and $\tau_1 = \tau_2$. Furthermore, note that, for any complex-valued random variable $Z$, with complex conjugate $\bar{Z}$,
\begin{equation}\label{eqn:rem:CI:f:2}
\Var(\Re Z) = \frac{1}{2} \big( \Var(Z) + \Re \Cov(Z, \bar{Z})\big);\ \
\Var(\Im Z) = \frac{1}{2} \big( \Var(Z) - \Re \Cov(Z, \bar{Z})\big),
\end{equation}
and we have $\overline{\HH{1}{2}} = \HH{2}{1}$.

\subsection[Pointwise confidence bands for R]{Pointwise confidence bands for $\bfrakR$}\label{rem:CI:R}
We utilise Theorem~\ref{thm:AsympCohRankEstimator} to construct pointwise asymptotic $(1-\alpha)$-level confidence bands for the real and imaginary parts of $\mathfrak{R}^{j_1, j_2}(\omega;\tau_1,\tau_2)$ as follows:

\[C^{(2)}_{{\rm r},n}(\omega_{kn}; \tau_1, \tau_2) := \Re \hat{\mathfrak{R}}^{j_1, j_2}_{n,R}(\omega_{kn}; \tau_1, \tau_2) \pm \Re \sigma_{(2)}^{j_1, j_2}(\omega_{kn}; \tau_1, \tau_2) \Phi^{-1}(1-\alpha/2),\]
for the real part, and 
\[C^{(2)}_{{\rm i},n}(\omega_{kn}; \tau_1, \tau_2) := \Im \hat{\mathfrak{R}}^{j_1, j_2}_{n,R}(\omega_{kn}; \tau_1, \tau_2) \pm \Im \sigma_{(2)}^{j_1, j_2}(\omega_{kn}; \tau_1, \tau_2) \Phi^{-1}(1-\alpha/2),\]
for the imaginary part of the quantile coherency. Here, $\Phi$ stands for the cdf of the standard normal distribution,
\begin{equation*}
\big( \Re \sigma_{(2)}^{j_1, j_2}(\omega_{kn}; \tau_1, \tau_2) \big)^2
:= 0 \vee \begin{cases}
			0 & \text{if $j_1 = j_2$} \\
			  & \text{\quad and $\tau_1 = \tau_2$,} \\
      \frac{1}{2} \big(\Cov( \LL{1}{2}, \LL{1}{2} ) + \Re \Cov( \LL{1}{2}, \LL{2}{1} ) \big) & \text{otherwise},
\end{cases}
\end{equation*}
and
\begin{equation*}
\big( \Im \sigma_{(2)}^{j_1, j_2}(\omega_{kn}; \tau_1, \tau_2) \big)^2
:= 0 \vee \begin{cases}
      0 & \text{if $j_1 = j_2$}\\
			  & \text{\quad  and $\tau_1 = \tau_2$,}\\
      \frac{1}{2} \big(\Cov( \LL{1}{2}, \LL{1}{2} ) - \Re \Cov( \LL{1}{2}, \LL{2}{1} ) \big) & \text{otherwise}.
\end{cases}
\end{equation*}

The definition of $\sigma_{(2)}^{j_1, j_2}(\omega_{kn}; \tau_1, \tau_2)$ is motivated by~\eqref{eqn:rem:CI:f:2} and the fact that we have $\overline{\LL{1}{2}} = \LL{2}{1}$. Furthermore, note that $\hat{\mathfrak{R}}^{j_1, j_2}_{n,R}(\omega_{kn}; \tau_1, \tau_2) = 1$, if $j_1 = j_2$ and $\tau_1 = \tau_2$..

\noindent In the definition of $\sigma_{(2)}^{j_1, j_2}(\omega_{kn}; \tau_1, \tau_2)$ we have used $\Cov( \LL{a}{b}, \LL{c}{d} )$ to denote an estimator for
\[\Cov\big(\IL^{j_1, j_2}(\omega_{kn}; \tau_1, \tau_2 \big), \IL^{j_3, j_4}(\omega_{kn}; \tau_3, \tau_4 )\big).\]
Recalling the definition of he limit process in Theorem~\ref{thm:AsympCohRankEstimator} we derive the following expression:
\begin{equation*}
\begin{split}
	& \frac{1}{\sqrt{\ff{1}{1} \ff{2}{2} \ff{3}{3} \ff{4}{4}}}
	\Cov \Big( \HH{1}{2} - \frac{1}{2} \frac{\ff{1}{2}}{\ff{1}{1}} \HH{1}{1} - \frac{1}{2} \frac{\ff{1}{2}}{\ff{2}{2}} \HH{2}{2},
	\HH{3}{4} - \frac{1}{2} \frac{\ff{3}{4}}{\ff{3}{3}} \HH{3}{3} - \frac{1}{2} \frac{\ff{3}{4}}{\ff{4}{4}} \HH{4}{4} \Big) \\
	& = \frac{\Cov( \HH{1}{2}, \HH{3}{4})}{\sqrt{\ff{1}{1} \ff{2}{2} \ff{3}{3} \ff{4}{4}}} 
		- \frac{1}{2} \frac{\overline{\ff{3}{4}} \Cov( \HH{1}{2}, \HH{3}{3})}{\sqrt{\ff{1}{1} \ff{2}{2} \ff{3}{3}^3 \ff{4}{4}}}
		- \frac{1}{2} \frac{\overline{\ff{3}{4}} \Cov( \HH{1}{2}, \HH{4}{4})}{\sqrt{\ff{1}{1} \ff{2}{2} \ff{3}{3} \ff{4}{4}^3}} \\
	& \qquad
	- \frac{1}{2} \frac{\ff{1}{2} \Cov( \HH{1}{1}, \HH{3}{4})}{\sqrt{\ff{1}{1}^3 \ff{2}{2} \ff{3}{3} \ff{4}{4}}} 
		+ \frac{1}{4} \frac{\ff{1}{2} \overline{\ff{3}{4}} \Cov( \HH{1}{1}, \HH{3}{3})}{\sqrt{\ff{1}{1}^3 \ff{2}{2} \ff{3}{3}^3 \ff{4}{4}}}
		+ \frac{1}{4} \frac{\ff{1}{2} \overline{\ff{3}{4}} \Cov( \HH{1}{1}, \HH{4}{4})}{\sqrt{\ff{1}{1}^3 \ff{2}{2} \ff{3}{3} \ff{4}{4}^3}} \\
	& \qquad
	- \frac{1}{2} \frac{\ff{1}{2} \Cov( \HH{2}{2}, \HH{3}{4})}{\sqrt{\ff{1}{1} \ff{2}{2}^3 \ff{3}{3} \ff{4}{4}}} 
		+ \frac{1}{4} \frac{\ff{1}{2} \overline{\ff{3}{4}} \Cov( \HH{2}{2}, \HH{3}{3})}{\sqrt{\ff{1}{1} \ff{2}{2}^3 \ff{3}{3}^3 \ff{4}{4}}}
		+ \frac{1}{4} \frac{\ff{1}{2} \overline{\ff{3}{4}} \Cov( \HH{2}{2}, \HH{4}{4})}{\sqrt{\ff{1}{1} \ff{2}{2}^3 \ff{3}{3} \ff{4}{4}^3}},
	\end{split}
\end{equation*}
where we have written $\ff{a}{b}$ for the quantile spectral density $\mathfrak{f}^{j_a, j_b}(\omega_{kn}; \tau_a, \tau_b)$, and $\HH{a}{b}$ for the limit distribution $\IH^{j_a, j_b}(\omega_{kn}; \tau_a, \tau_b\big)$ for any $a,b = 1,2,3,4$).

Thus, considering the special case where $\tau_3 = \tau_1$ and $\tau_4 = \tau_2$, we have
\begin{equation}\label{eqn:rem:CI:R:1}
\begin{split}
	& \Cov( \LL{1}{2}, \LL{1}{2} ) \\
	& = \frac{1}{\ff{1}{1} \ff{2}{2}} \Big( \Cov( \HH{1}{2}, \HH{1}{2} )
		- \Re \frac{\ff{1}{2} \Cov( \HH{1}{1}, \HH{1}{2} )}{\ff{1}{1}}
		- \Re \frac{\ff{1}{2} \Cov( \HH{2}{2}, \HH{1}{2} )}{\ff{2}{2}} \\
	& \qquad
		+ \frac{1}{4} |\ff{1}{2}|^2 \Big( \frac{ \Cov( \HH{1}{1}, \HH{1}{1})}{\ff{1}{1}^2}
		+ 2 \Re \frac{ \Cov( \HH{1}{1}, \HH{2}{2})}{\ff{1}{1} \ff{2}{2}}
		+ \frac{ \Cov( \HH{2}{2}, \HH{2}{2})}{\ff{2}{2}^2} \Big) \Big)
	\end{split}
\end{equation}
and for the special case where $\tau_3 = \tau_1$ and $\tau_4 = \tau_2$ we have 
\begin{equation*}
\begin{split}
	& \Cov( \LL{1}{2}, \LL{2}{1} ) \\
	& = \frac{1}{\ff{1}{1} \ff{2}{2}} \Big( \Cov( \HH{1}{2}, \HH{2}{1})
		- \frac{\ff{1}{2} \Cov( \HH{1}{2}, \HH{2}{2})}{\ff{2}{2}}
		- \frac{\ff{1}{2} \Cov( \HH{1}{2}, \HH{1}{1})}{\ff{1}{1}} \\
	& \qquad
		+ \frac{1}{4} \ff{1}{2}^2 \Big(
		\frac{ \Cov( \HH{1}{1}, \HH{1}{1})}{ \ff{1}{1}^2 }
		+ 2 \Re \frac{\Cov( \HH{1}{1}, \HH{2}{2})}{\ff{1}{1} \ff{2}{2}}
		+  \frac{ \Cov( \HH{2}{2}, \HH{2}{2})}{\ff{2}{2}^2} \Big) \Big).
	\end{split}
\end{equation*}
We substitute consistent estimators for the unknown quantities. To do so we abuse notation using $\ff{a}{b}$ to denote $\tilde{G}^{j_a, j_b}_{n,R}(\omega_{kn}; \tau_a, \tau_b)$ and write $\Cov( \HH{a}{b}, \HH{c}{d})$ for the quantity defined in~\eqref{eqn:rem:CI:f:1}.

\section{Proofs of the results in Sections~4 and~S4}\label{sec:proofs}

In this section the proofs to the results in Sections~\ref{sec:asymp} and~\ref{sec:asymp2} are given. Before we begin, note that by a trivial generalisation of Proposition~3.1 in \cite{kley2014} we have that Assumption~\ref{ass:exp_alpha_mix} implies that there exist constants $\bar\rho\in(0,1)$ and $K < \infty$ such that, for arbitrary intervals $A_1,...,A_p \subset \mathbb{R}$, arbitrary indices $j_1, \ldots, j_p \in \{1,\ldots,d\}$ and times $t_1,...,t_p \in \IZ$,
\begin{equation}\label{eq:boundcum}
|\cum(I\{X_{t_1, j_1} \in A_1\}, \ldots, I\{X_{t_p, j_p} \in A_p\})| \leq K \bar\rho^{\max_{i,j} |t_i-t_j|}.
\end{equation}
We will use this fact several times throughout the proofs in this section.

\subsection{Proof of Theorem~\ref{thm:AsympCohRankEstimator}}
\label{proof:thm:AsympCohRankEstimator}

By a Taylor expansion we have, for every $x, x_0 > 0$,

\[\frac{1}{\sqrt{x}} = \frac{1}{\sqrt{x_0}} - \frac{1}{2} \frac{1}{\sqrt{x_0^3}} (x-x_0) + \frac{3}{8} \xi_{x,x_0}^{-5/2} (x-x_0)^2,\]
where $\xi_{x,x_0}$ is between $x$ and $x_0$. Let $R_n(x,x_0) := \frac{3}{8} \xi_{x,x_0}^{-5/2} (x-x_0)^2$, then
\begin{equation}
\label{eqn:expansion}
	\frac{x}{\sqrt{y z}} - \frac{x_0}{\sqrt{y_0 z_0}}
	= \frac{1}{\sqrt{y_0 z_0}} \Big( (x - x_0) - \frac{1}{2} \frac{x_0}{y_0} (y-y_0) - \frac{1}{2} \frac{x_0}{z_0} (z-z_0) + r_n \Big),
\end{equation}
where
\begin{equation*}
\begin{split}
	r_n & = (x - x_0) \Big(- \frac{1}{2} \frac{1}{y_0} (y-y_0) - \frac{1}{2} \frac{1}{z_0} (z-z_0) \Big) \\
	& \quad + x \Big( R_n(y,y_0) \sqrt{y_0} \big(1 - \frac{1}{2} \frac{1}{z_0} (z-z_0) \big) +  R_n(z,z_0) \sqrt{z_0} \big(1 - \frac{1}{2} \frac{1}{y_0} (y-y_0) \big)   \\
	& \qquad + \frac{1}{4} \frac{1}{y_0} (y-y_0) \frac{1}{z_0} (z-z_0)  + \sqrt{y_0 z_0} R_n(y,y_0) R_n(z,z_0)\Big)
\end{split}
\end{equation*}

Write $\ff{a}{b}$ for $\mathfrak{f}^{j_a, j_b}(\omega; \tau_a, \tau_b)$, $\GG{a}{b}$ for $\hat G_{n,R}^{j_a, j_b} (\omega; \tau_a, \tau_b)$, and $\BB{a}{b}$ for $\{ \bB_n^{(k)}(\omega; \tau_a, \tau_b)\}_{j_a, j_b}$ ($a, b = 1,2,3,4$). We want to employ~\eqref{eqn:expansion} and to this end let
\begin{align*}
x   & := \GG{a}{b}							&	y   & := \GG{a}{a}							&	z 	& := \GG{b}{b} \\
x_0 & := \ff{a}{b} + \BB{a}{b}	&	y_0 & := \ff{a}{a} + \BB{a}{a}	& z_0 & := \ff{b}{b} + \BB{b}{b}
\end{align*}
By Theorem~\ref{thm:AsympDensityRankEstimator} the differences $x-x_0$, $y-y_0$, and $z-z_0$ are in $O_p ((n b_n)^{-1/2})$, uniformly with respect to $\tau_1, \tau_2$. Under the assumption that $n b_n \rightarrow \infty$, as $n \rightarrow \infty$, this entails $\GG{a}{a} - \BB{a}{a} \rightarrow \ff{a}{a}$, in probability. For $\varepsilon \leq \tau_1, \tau_2 \leq 1-\varepsilon$, we have $\ff{a}{a} > 0$, such that, by the Continuous Mapping Theorem we have $(\GG{a}{a} - \BB{a}{a})^{-5/2} \rightarrow \ff{a}{a}^{-5/2}$, in probability. As $\BB{a}{a} = o(1)$, we have $y^{-5/2} - y_0^{-5/2} = o_p(1)$. Finally, due to
\[\xi_{y,y_0}^{-5/2} \leq y_n^{-5/2} \vee y_0^{-5/2} \leq (y_n^{-5/2} - y_0^{-5/2}) \vee 0 + y_0^{-5/2} = o_p(1) + O(1) = O_p(1),\]
we have that $R_n(y,y_0) = O_p((n b_n)^{-1})$.\\
Analogous arguments yields $R_n(z,z_0) = O_p((n b_n)^{-1})$. Thus we have shown that
\begin{multline*}
		\hat{\mathfrak{R}}^{j_1, j_2}_{n,R}(\omega; \tau_1, \tau_2) - \frac{\ff{a}{b} + \BB{a}{b}}{\sqrt{\ff{a}{a} + \BB{a}{a}}\sqrt{\ff{b}{b} + \BB{b}{b}}} \\
		= \frac{1}{\sqrt{\ff{1}{1} \ff{2}{2}}} \Big([\GG{1}{2} - \ff{1}{2} - \BB{1}{2}]
		- \frac{1}{2} \frac{\ff{1}{2}}{\ff{1}{1}} [\GG{1}{1} - \ff{1}{1} - \BB{1}{1}]
		- \frac{1}{2} \frac{\ff{1}{2}}{\ff{2}{2}} [\GG{2}{2} - \ff{2}{2} - \BB{2}{2}]\Big) \\ + O_p \big(1/ (n b_n) \big),
\end{multline*}
with the $O_p$ holding uniformly with respect to $\tau_1, \tau_2$. Further more, note that
\begin{multline*}
	\frac{\ff{a}{b} + \BB{a}{b}}{\sqrt{\ff{a}{a} + \BB{a}{a}}\sqrt{\ff{b}{b} + \BB{b}{b}}}
	= \frac{\ff{a}{b}}{\sqrt{\ff{a}{a} \ff{b}{b}}}
		+ \frac{1}{\sqrt{\ff{a}{a} \ff{b}{b}}} \Big( \BB{a}{b}
		- \frac{1}{2} \frac{\ff{a}{b}}{\ff{a}{a}} \BB{a}{a} - \frac{1}{2} \frac{\ff{a}{b}}{\ff{b}{b}} \BB{b}{b} \Big) \\
	+ O(|\BB{a}{b}| ( \BB{a}{a} + \BB{b}{b} ) + \BB{a}{a}^2 + \BB{b}{b}^2 + \BB{a}{a} \BB{b}{b}),
\end{multline*}
where we have used~\eqref{eqn:expansion} again.
By Lemma~\ref{lem:LipschitzLaplaceSD} we have that
\begin{equation*}
\sup_{\tau_1, \tau_2 \in [\varepsilon, 1-\varepsilon]} \Big| \frac{{\rm d}^{\ell}}{{\rm d}\omega^{\ell}}\mathfrak{f}^{j_1, j_2}(\omega; \tau_1, \tau_2) \Big| \leq C_{\varepsilon,\ell}.
\end{equation*}
Therefore, $b_n$ satisfies
\begin{equation*}
\sup_{\tau_1, \tau_2 \in [\varepsilon, 1-\varepsilon]} \Big| \sum_{\ell=2}^k \frac{b_n^{\ell}}{\ell!} \int_{ -\pi}^{\pi} v^{\ell} W(v) dv \frac{{\rm d}^{\ell}}{{\rm d}\omega^{\ell}}\mathfrak{f}^{j_1, j_2}(\omega; \tau_1, \tau_2) \Big| = o\big( (n b_n)^{-1/4} \big),
\end{equation*}
for all $j_1, j_2=1,\ldots,d$,
which implies that
\[|\BB{a}{b}| ( \BB{a}{a} + \BB{b}{b} ) + \BB{a}{a}^2 + \BB{b}{b}^2 + \BB{a}{a} \BB{b}{b} = o\big( (n b_n)^{-1/2} \big).\] Therefore,
\begin{multline*}
\sqrt{n b_n}\Big(\hat{\mathfrak{R}}^{j_1, j_2}_{n,R}(\omega; \tau_1, \tau_2)
- \mathfrak{R}^{j_1, j_2}(\omega;\tau_1,\tau_2) \\
- \frac{1}{\sqrt{\ff{a}{a} \ff{b}{b}}} \Big( \BB{a}{b} 
		- \frac{1}{2} \frac{\ff{a}{b}}{\ff{a}{a}} \BB{a}{a} - \frac{1}{2} \frac{\ff{a}{b}}{\ff{b}{b}} \BB{b}{b} \Big)
		\Big)_{\tau_1, \tau_2 \in [0,1]}
\end{multline*}
and 
\begin{equation}
\label{eqn:expExpr}
	\frac{\sqrt{n b_n}}{\sqrt{\ff{1}{1} \ff{2}{2}}} \Big([\GG{1}{2} - \ff{1}{2} - \BB{1}{2}]
		- \frac{1}{2} \frac{\ff{1}{2}}{\ff{1}{1}} [\GG{1}{1} - \ff{1}{1} - \BB{1}{1}]
		- \frac{1}{2} \frac{\ff{1}{2}}{\ff{2}{2}} [\GG{2}{2} - \ff{2}{2} - \BB{2}{2}]\Big)
\end{equation}
are asymptotically equivalent in the sense that if one of the two converges weakly in $\ell_{\mathbb{C}^{d \times d}}^{\infty}([0,1]^2)$, then so does the other.
The assertion then follows by Theorem~\ref{thm:AsympDensityRankEstimator}, Slutzky's lemma and the Continuous Mapping Theorem.\hfill\qed

\subsection{Proof of Proposition~\ref{prop:inr}}\label{proof:prop:inr}

\noindent The proof resembles the proof of Proposition~3.4 in~\cite{kley2014}, where the univariate case was handled.
For $j=1,\ldots,d$ we have, from the continuity of~$F_j$ that the ranks of the random variables $X_{0,j},...,X_{n-1,j}$ and $F_j(X_{0,j}),...,F_j(X_{n-1,j})$ coincide almost surely. Thus, without loss of generality,  we can assume that the CCR-periodogram is computed from the unobservable data $(F_j(X_{0,j}))_{j=1,\ldots,d},...,(F_j(X_{n-1,j}))_{j=1,\ldots,d}$. In particular, we can assume the marginals to be uniform.

Applying the Continuous Mapping Theorem afterward, it suffices to prove
\begin{equation}
  \label{prop:inr:eqn:dfts}
  \Big( n^{-1/2} d^j_{n,R}(\omega; \tau) \Big)_{ \tau \in [0,1], j=1,\ldots,d} \Rightarrow
  \Big( \mathbb{D}^j (\omega; \tau) \Big)_{ \tau \in [0,1], j=1,\ldots,d}
  \quad \text{in } \ell_{\IC^d}^\infty([0,1]),
\end{equation}
where $\ell_{\IC^d}^\infty([0,1])$ is the space of bounded functions $[0,1] \rightarrow \IC^d$ that we identify with the product space $\ell^{\infty}([0,1])^{2d}$.
Let
\[d_{n,U}^{j}(\omega; \tau) := \sum_{t=0}^{n-1} I\{F_{j}(X_{t,j}) \leq \tau\} \ee^{- \ii \omega t},\]
$j = 1,\ldots, d$, $\omega\in  \IR$, $\tau \in [0,1]$, and note that for \eqref{prop:inr:eqn:dfts} to hold, it is sufficient that
\[\big( n^{-1/2} d^j_{n,U}(\omega; \tau) \big)_{ \tau \in [0,1], j=1,\ldots,d}\]
satisfies the following two conditions:
\begin{itemize}
  \item[(i1)] convergence of the finite-dimensional distributions, i.\,e.,
    \begin{equation}
      \label{prop:inr:eqn:fidis}
      \big( n^{-1/2} d^{j_\ell}_{n,U}(\omega_{\ell}; \tau_{\ell}) \big)_{\ell=1,\ldots,k} \xrightarrow{d} \big( \mathbb{D}^{j_\ell}(\omega_{\ell}; \tau_\ell) \big)_{\ell=1,\ldots,k},
    \end{equation}
   for any $(j_{\ell}, \tau_{\ell}) \in \{1,\ldots,d\} \times [0,1]$, $\omega_{\ell} \neq 0 \mod 2\pi$, $\ell =1,\ldots,k$ and $k \in \IN$;
  \item[(i2)]  stochastic equicontinuity:
    for any $x > 0$ and any $\omega \neq 0 \mod 2\pi$,
    \begin{equation}
      \label{prop:inr:eqn:stochequicont}
      \lim_{\delta \downarrow 0} \limsup_{n \rightarrow \infty} \IP \Big( \sup_{\substack{\tau_1, \tau_2 \in [0,1] \\ | \tau_1 - \tau_2 | \leq \delta}} | n^{-1/2} (d^{j}_{n,U}(\omega; \tau_1) - d^{j}_{n,U}(\omega; \tau_2)) | > x \Big) = 0, \quad \forall j = 1,\ldots,d.
    \end{equation}
\end{itemize}

\noindent Under (i1) and (i2), an application of Theorems~1.5.4 and~1.5.7 from \cite{vanderVaartWellner1996} then yields
\begin{equation}
  \label{prop:inr:eqn:dftsU}
  \Big( n^{-1/2} d^j_{n,U}(\omega; \tau) \Big)_{ \tau \in [0,1], j = 1,\ldots,d} \Rightarrow
  \Big( \mathbb{D}^j(\omega; \tau) \Big)_{ \tau \in [0,1], j = 1,\ldots,d}
  \quad \text{in } \ell_{\IC^d}^\infty([0,1]).
\end{equation}
In combination with
\begin{equation}
\label{prop:inr:eqn:equiv-dftsRandU}
  \sup_{\tau \in [0,1]} | n^{-1/2} (d_{n,R}^{j}(\omega; \tau) - d_{n,U}^{j}(\omega; \tau) ) | = o_p(1), \quad \text{for $\omega \neq 0 \mod 2\pi$, $j = 1,\ldots,d$},
\end{equation}
which we will prove below, \eqref{prop:inr:eqn:dftsU} yields the desired result: \eqref{prop:inr:eqn:dfts}.
For the proof of~\eqref{prop:inr:eqn:equiv-dftsRandU}, we denote by $\hat F_{n,j}^{-1}(\tau) := \inf\{x: \hat F_{n,j} (x) \geq \tau\}$ the generalised inverse of $\hat F_{n,j}$ and let $\inf \emptyset := 0$. Then, we have, as in (7.25) of~\cite{kley2014}, that
\begin{equation}\label{eq:dntau}
\sup_{\omega \in \IR}\sup_{\tau\in [0,1]}\Big|d_{n,R}^{j}(\omega; \tau) -  d_{n,U}^{j}(\omega; \hat F_{n,j}^{-1}(\tau))\Big| \leq n \sup_{\tau \in [0,1]} |\hat F_{n,j}(\tau) - \hat F_{n,j}(\tau-)| = O_p(n^{1/2k})
\end{equation}
where $\hat F_{n,j}(\tau-) := \lim_{\xi \uparrow 0} \hat F_{n,j}(\tau-\xi)$. The $O_p$-bound in~\eqref{eq:dntau} follows from Lemma~\ref{lem:BoundCdf}. Therefore, it suffices to bound the terms
\[
\sup_{\tau \in [0,1]} n^{-1/2} |
  d_{n,U}^{j}(\omega; \hat F_{n,j}^{-1}(\tau)) - d_{n,U}^{j}(\omega, \tau) ) |, \text{ for all $j = 1,\ldots, d$.}
\]
To do so, note that, for  any $x>0$ and $\delta_n = o(1)$ satisfying $n^{1/2} \delta_n \rightarrow \infty$, we have
\begin{equation*}
\begin{split}
    &\IP \Big( \sup_{\tau \in [0,1]} n^{-1/2} |
  d_{n,U}^{j}(\omega; \hat F_{n,j}^{-1}(\tau)) - d_{n,U}^{j}(\omega; \tau) ) | > x \Big) \\
  &\leq  \IP \Big( \sup_{\tau \in [0,1]} \sup_{|u - \tau| \leq \delta_n} |  d_{n,U}^{j}(\omega; u) - d_{n,U}^{j}(\omega; \tau) | > x n^{1/2},
  \sup_{\tau \in [0,1]} | \hat F_{n, j}^{-1}(\tau) - \tau | \leq \delta_n \Big) \\
  & \quad + \IP \Big(\sup_{\tau \in [0,1]} | \hat F_{n,j}^{-1}(\tau) - \tau | > \delta_n \Big) = o(1) + o(1).
\end{split}
\end{equation*}
The first $o(1)$ follows from~\eqref{prop:inr:eqn:stochequicont}. The second one is a consequence of Lemma~\ref{lem:quantproc}.

It thus remains to prove~\eqref{prop:inr:eqn:fidis} and~\eqref{prop:inr:eqn:stochequicont}.
For any fixed $j=1,\ldots,d$ the process $\big( d_{n,U}^j(\omega, \tau) \big)_{\tau \in [0,1]}$ is determined by the univariate time series $X_{0,j}, \ldots, X_{n-1,j}$. Under the assumptions made here, \eqref{prop:inr:eqn:stochequicont} therefore follows from (8.7) in~\cite{kley2014}.

Finally, we establish~\eqref{prop:inr:eqn:fidis}, by employing Lemma~\ref{lem:OrderDftYinA} in combination with Lemma~P4.5 and Theorem~4.3.2 from~\cite{Brillinger1975}. More precisely, to apply Lemma~P4.5 from \cite{Brillinger1975}, we have to verify that, for any $j_1, \ldots, j_{\ell} \in \{1,\ldots,d\}$, $\tau_1, \ldots, \tau_{\ell} \in [0,1]$, $\ell \in \IN$, and $\omega_1, \ldots, \omega_{\ell} \neq 0 \mod 2\pi$, all cumulants of the vector
\[n^{-1/2} \big( d_{n,U}^{j_1}(\omega_1; \tau_1), d_{n,U}^{j_1}(-\omega_1; \tau_1), \ldots,
d_{n,U}^{j_{\ell}}(\omega_{\ell}; \tau_{\ell}), d_{n,U}^{j_{\ell}}(-\omega_{\ell}; \tau_{\ell}) \big)\]
converge to the corresponding cumulants of the vector
\[\big( \mathbb{D}^{j_1}(\omega_1; \tau_1), \mathbb{D}^{j_1}(-\omega_1; \tau_1), \ldots,
\mathbb{D}^{j_{\ell}}(\omega_{\ell}; \tau_{\ell}), \mathbb{D}^{j_{\ell}}(-\omega_{\ell}; \tau_{\ell}) \big).\]
For the cumulants of order one the arguments from the univariate case (cf. the proof of Proposition~3.4 in \cite{kley2014}) apply: we have
$|\IE (n^{-1/2} d_{n,U}^{j}(\omega; \tau))| = o(1)$,
for any $j = 1,\ldots, d$, $\tau \in [0,1]$ and fixed $\omega \neq 0 \mod 2\pi$.
Furthermore, for the cumulants of order two, applying Theorem~4.3.1 in \cite{Brillinger1975} to the bivariate process
\[(I\{X_{t,j_1} \leq q_{j_1}(\mu_1)\}, I\{X_{t, j_2} \leq q_{j_2}(\mu_2)\}),\]
we obtain
\[\cum(n^{-1/2} d_{n,U}^{i_1}(\lambda_1; \mu_1), n^{-1/2} d_{n,U}^{i_2}(\lambda_2; \mu_2))
= 2 \pi n^{-1} \Delta_n(\lambda_1 + \lambda_2) \mathfrak{f}^{i_1, i_2}(\lambda_1; \mu_1, \mu_2) + o(1)\]
for any $(i_1, \lambda_1, \mu_1), (i_2, \lambda_2, \mu_2) \in \bigcup_{\ell=1}^{k} \{(i_{\ell}, \omega_{\ell}, \tau_{\ell}), (j_{\ell}, -\omega_{\ell}, \tau_{\ell})\}$, which yields the correct second moment structure. The function $\Delta_n$ is defined in Lemma~\ref{lem:OrderDftYinA}.
Finally, the cumulants of order $J$, with $J \in \IN$ and $J \geq 3$, all tend to zero, as in view of Lemma~\ref{lem:OrderDftYinA}
\begin{multline*}
  \cum(n^{-1/2} d_{n,U}^{i_1}(\lambda_1; \mu_1), \ldots, n^{-1/2} d_{n,U}^{i_J}(\lambda_J; \mu_J)) \\
  \leq C n^{-J/2} (|\Delta_n(\sum_{j=1}^J \lambda_j)| + 1) \varepsilon (|\log \varepsilon|+1)^d = O(n^{-(J-2)/2}) = o(1),
\end{multline*}
for $(i_1, \lambda_1, \mu_1), \ldots ,(i_J, \lambda_J, \mu_J) \in \bigcup_{\ell=1}^k \{(i_{\ell}, \omega_{\ell}, \tau_{\ell}), (i_{\ell}, -\omega_{\ell}, \tau_{\ell})\}$, where $\varepsilon := \min_{j=1}^J \mu_j$.
This implies that the limit $\mathbb{D}^{j}(\tau;\omega)$ is Gaussian, and completes the proof of~\eqref{prop:inr:eqn:fidis}. Proposition~\ref{prop:inr} follows. \qed

\subsection{Proof of Theorem~\ref{thm:AsympDensityRankEstimator}}
\label{proof:thm:AsympDensityRankEstimator}

We proceed in a similar fashion as in the proof of the univariate estimator which was analysed in~\cite{kley2014}.
First, we state an {asymptotic representation} result by which the estimator $\hat \bG_{n,R}$ can be approximated, in a suitable uniform sense, by another process $\hat\bG_{n,U}$  which is not defined as a function of the standardised ranks $\hat F_{n,j}(X_{t,j})$, but as a function of the unobservable quantities $F_j(X_{t,j})$, $t=0,\ldots,n-1$, $j=1,\ldots,d$. More precisely, this process is defined as
\[\hat\bG_{n,U}(\omega; \tau_1, \tau_2) := ( \hat G^{j_1, j_2}_{n, U}(\omega; \tau_1, \tau_2) )_{j_1, j_2 = 1,\ldots, d},\]
where
\begin{align}
  \hat G^{j_1, j_2}_{n,U}(\omega; \tau_1, \tau_2)
    & := \frac{2\pi}{n} \sum_{s=1}^{n-1} W_n\big( \omega - 2\pi s / n \big) I_{n,U}^{j_1, j_2}(2 \pi s / n, \tau_1, \tau_2) \nonumber\\
  I_{n,U}^{j_1, j_2}(\omega; \tau_1, \tau_2)
    & := \frac{1}{2\pi n} d^{j_1}_{n,U}(\omega; \tau_1) d^{j_2}_{n,U}(-\omega; \tau_2) \nonumber\\
  d_{n,U}^{j}(\omega; \tau)
    & := \sum_{t=0}^{n-1} I\{F_{j}(X_{t,j}) \leq \tau\} \ee^{- \ii \omega t}. \label{eq:defdnu}
\end{align}
  
Theorem~\ref{thm:AsympDensityRankEstimator} then follows from the asymptotic representation of $\hat\bG_{n,R}$ by $\hat\bG_{n,U}$ (i.\,e., Theorem~\ref{thm:AsympDensityEstimator}(iii)) and the asymptotic properties of~$\hat\bG_{n,U}$ (i.\,e., Theorem~\ref{thm:AsympDensityEstimator}(i)--(ii)), which we now state:

\begin{theorem}\label{thm:AsympDensityEstimator}
Let Condition~\eqref{eq:boundcum} and Assumption~\ref{ass:W}  hold, and assume that the distribution functions $F_j$ of $X_{0,j}$ are continuous for all $j = 1,\ldots,d$.
Let $b_n$ satisfy the assumptions of Theorem~\ref{thm:AsympDensityRankEstimator}.
Then,
\begin{itemize}
\item[(i)] for any fixed $\omega \in  \IR  $,  as $n\to\infty$,
\[
\sqrt{nb_n}\big( \hat \bG_{n,U}(\omega; \tau_1, \tau_2)  - \IE \hat \bG_{n,U}(\omega; \tau_1, \tau_2) \big)_{\tau_1, \tau_2 \in [0,1]} \Rightarrow \IH(\omega; \cdot,\cdot)
\]
in $\ell_{\IC^{d \times d}}^{\infty}([0,1]^2)$, where the process $\IH(\omega; \cdot,\cdot)$ is defined in Theorem \ref{thm:AsympDensityRankEstimator};
\item[(ii)] still as $n\to\infty$,
\end{itemize}
\begin{multline*}
\sup_{\substack{j_1, j_2 \in \{1,\ldots,d\} \\ \tau_1, \tau_2 \in [0,1] \\ \omega \in \mathbb{R}}} \Big| \IE \hat G_{n,U}^{j_1, j_2}(\tau_1, \tau_2; \omega) - \mathfrak{f}^{j_1, j_2}(\omega; \tau_1, \tau_2) - \big\{ \bB_n^{(k)}(\omega; \tau_1, \tau_2) \big\}_{j_1, j_2} \Big| \\ = O((nb_n)^{-1}) + o(b_n^k),
\end{multline*}
\begin{itemize}
\item[] where $\big\{ \bB_n^{(k)}(\omega; \tau_1, \tau_2) \big\}_{j_1, j_2} $ is defined in~\eqref{def:bias};
\item[(iii)] for any fixed $\omega \in \mathbb{R}$,
\[
\sup_{\substack{j_1, j_2 \in \{1,\ldots,d\} \\ \tau_1, \tau_2 \in [0,1]}} | \hat G_{n,R}^{j_1, j_2}(\tau_1, \tau_2; \omega) - \hat G_{n,U}^{j_1, j_2}(\tau_1, \tau_2; \omega) | = o_p \big( (n b_n)^{-1/2}+ b_n^k \big);
\]
if moreover the kernel $W$ is uniformly Lipschitz-continuous, this bound is uniform with respect to $\omega \in \mathbb{R}$.
\end{itemize}
\end{theorem}

The proof of Theorem~\ref{thm:AsympDensityEstimator} is lengthy, technical and in many places similar to the proof of Theorem~3.6 in~\cite{kley2014}. We provide the proof in Sections~\ref{Parti}--\ref{Partiii}, with technical details deferred to Section~\ref{sec:auxlemmas}. For the reader's convenience we first give a brief description of the necessary steps.

Part (ii) of Theorem~\ref{thm:AsympDensityEstimator} can be proved along the lines of classical results from~\cite{Brillinger1975}, but uniformly with respect to the arguments $\tau_1$ and $\tau_2$. Parts~(i) and (iii) require additional arguments that are different from the classical theory. These additional arguments are due to the fact that the estimator is a stochastic process and stochastic equicontinuity of
\begin{equation}\label{thm:AsympDensityEstimator_i}
 \big( \hat H^{j_1, j_2}_n(a; \omega) \big)_{a \in [0,1]^2} := \sqrt{nb_n}\big( \hat G^{j_1, j_2}_{n,U}(\omega; \tau_1, \tau_2)  - \IE \hat G^{j_1, j_2}_{n,U}(\omega; \tau_1, \tau_2) \big)_{\tau_1, \tau_2 \in [0,1]}
\end{equation}
for all $j_1, j_2 = 1,\ldots,d$ has to be proven to ensure that the convergence holds not only pointwise, but also uniformly. The key to the proof of (i) and (iii) is a uniform bound on the increments $\hat H^{j_1, j_2}_n(a;\omega) - \hat H^{j_1, j_2}_n(b;\omega)$ of the process $\hat H^{j_1, j_2}_{n}$. This bound is  needed to show the stochastic equicontinuity of the process. To employ a restricted chaining technique (cf. Lemma~\ref{lem:Thm224VW}), we require two different bounds. First, we prove a general bound, uniform in $a$ and $b$, on the moments of the increments $\hat H^{j_1, j_2}_n(a;\omega) - \hat H^{j_1, j_2}_n(b;\omega)$ (cf. Lemma~\ref{lem:SixthMomIncrementHn}). Second, we prove a sharper bound on the increments $\hat H^{j_1, j_2}_n(a;\omega) - \hat H^{j_1, j_2}_n(b;\omega)$ when $a$ and $b$ are ``sufficiently close'' (cf. Lemma~\ref{lem:OrderSmallIncrements}).

Condition~\eqref{lem:SixthMomIncrementHn:eqn:Assumption} which we will required for Lemma~\ref{lem:SixthMomIncrementHn} to hold is rather general. In Lemma~\ref{lem:OrderDftYinA} we prove that condition~\eqref{eq:boundcum}, which is implied by Assumption~\ref{ass:exp_alpha_mix}, implies~\eqref{lem:SixthMomIncrementHn:eqn:Assumption}.

\subsubsection{Proof of Theorem~\ref{thm:AsympDensityEstimator}(i)}
\label{Parti}

It is sufficient to prove the following two claims:
\begin{itemize}
  \item[(i1)] convergence of the finite-dimensional distributions of the process~\eqref{thm:AsympDensityEstimator_i}, that is,
    \begin{equation}
      \label{thm:AsympDensityEstimator:eqn:fidis}
      \big(\hat H^{j_{1\ell}, j_{2\ell}}_n\big( (a_{1\ell},a_{2\ell}) ;\omega_j\big)\big)_{j=1,\ldots,k} \xrightarrow{d} \big( \IH^{j_{1\ell}, j_{2\ell}}\big( (a_{1\ell},a_{2\ell}) ;\omega_j\big)\big)_{j=1,\ldots,k}
    \end{equation}
   for any  $(j_{1\ell}, j_{2\ell}, a_{1\ell},a_{2\ell},\omega_\ell) \in \{1,\ldots,d\} \times [0,1]^2\times \IR$, $\ell=1,\ldots,k$ and $k \in \IN$;
  \item[(i2)]  stochastic equicontinuity:
    for any $x > 0$, any $\omega \in \IR$, and any $j_1, j_2 = 1,\ldots,d$,
    \begin{equation}
      \label{thm:AsympDensityEstimator:eqn:stochequicont}
      \lim_{\delta \downarrow 0} \limsup_{n \rightarrow \infty} \IP \Big( \sup_{\substack{a, b \in [0,1]^2 \\ \| a - b \|_1 \leq \delta}} | \hat H^{j_1, j_2}_n(a;\omega) - \hat H^{j_1, j_2}_n(b; \omega) | > x\Big) = 0.
    \end{equation}
\end{itemize}

By~\eqref{thm:AsympDensityEstimator:eqn:stochequicont} we have stochastic equicontinuity of all real parts~$\Re \hat H^{j_1, j_2}_n( \cdot ;\omega)$ and imaginary parts $\Im \hat H_n^{j_1, j_2}( \cdot ;\omega)$. Therefore, in view of Theorems~1.5.4 and~1.5.7 in \cite{vanderVaartWellner1996}, we will have proven part~(i).

First we prove (i1). For fixed $\tau_1, \tau_2$, $\hat G^{j_1, j_2}_{n,U}(\omega; \tau_1,\tau_2)$ is the traditional smoothed periodogram estimator of the cross-spectrum of the {clipped processes}  $(I\{F_{j_1}(X_{t,j_1})\leq \tau_1\})_{t\in\IZ}$ and $(I\{F_{j_2}(X_{t,j_2})\leq \tau_2\})_{t\in\IZ}$ [see Chapter 7.1 in~\citet{Brillinger1975}]. Thus,~\eqref{thm:AsympDensityEstimator:eqn:fidis} follows from Theorem~7.4.4 in~\citet{Brillinger1975}, by which these estimators are asymptotically jointly Gaussian. The first and second moment structures of the limit are given by Theorem~7.4.1 and Corollary~7.4.3 in~\citet{Brillinger1975}. The joint convergence~\eqref{thm:AsympDensityEstimator:eqn:fidis} follows. Note that condition~\eqref{eq:boundcum}, which is implied by Assumption~\ref{ass:exp_alpha_mix}, implies the summability condition [i.\,e., Assumption~2.6.2($\ell$) in~\cite{Brillinger1975}, for every $\ell$] required for the three theorems in~\cite{Brillinger1975} to be applied.

Now to the proof of (i2). The Orlicz norm $\| X \|_{\Psi} = \inf \{C > 0 : \IE \Psi( |X| / C ) \leq 1\}$ with $\Psi(x) := x^6$ coincides with the $L_6$ norm~$\| X \|_6 = ( \IE | X |^6 )^{1/6}$. Therefore,  for any $\kappa>0$ and   sufficiently small $\|a-b\|_1$, we have by Lemma~\ref{lem:SixthMomIncrementHn} and Lemma~\ref{lem:OrderDftYinA} that
\[
\| \hat H^{j_1, j_2}_n(a;\omega) - \hat H^{j_1, j_2}_n(b;\omega) \|_{\Psi} \leq K\Big( \frac{\|a-b\|_1^\kappa}{(n b_n)^2} + \frac{\|a-b\|_1^{2\kappa}}{n b_n} + \|a-b\|_1^{3\kappa} \Big)^{1/6}.
\]
Consequently, for all $a, b$ with $\|a-b\|_1$ sufficiently small and $\|a-b\|_1 \geq (n b_n)^{-1/\gamma}$ and all $\gamma \in (0,1)$ such that $\gamma < \kappa$,
\begin{equation*}
  \| \hat H^{j_1, j_2}_n(a;\omega) - \hat H^{j_1, j_2}_n(b;\omega) \|_{\Psi}
 \leq \bar K \|a-b\|_1^{\gamma/2}.
\end{equation*}
Note that
$\|a-b\|_1\geq (n b_n)^{-1/\gamma}$ if and only if $d(a,b) := \|a-b\|_1^{\gamma/2} \geq (n b_n)^{-1/2} =: \bar \eta_n / 2$.
The \emph{packing number} \cite[p.\,98]{vanderVaartWellner1996}
$D(\varepsilon,d)$ of $([0,1]^2,d)$
satisfies $D(\varepsilon,d) \asymp \varepsilon^{-4/\gamma}$.
By Lemma~\ref{lem:Thm224VW}, we therefore have, for all $x, \delta > 0$ and $\eta \geq \bar \eta_n$,
\begin{equation*}
\begin{split}
  & \IP\Big( \sup_{\| a - b\|_1 \leq \delta^{2/\gamma}} | \hat H^{j_1, j_2}_n(a;\omega) - \hat H^{j_1, j_2}_n(b;\omega)| > x\Big) \\
  & = \IP\Big( \sup_{d(a,b) \leq \delta} | \hat H^{j_1, j_2}_n(a;\omega) - \hat H^{j_1, j_2}_n(b;\omega)| > x\Big) \\
  & \leq \Bigg[ \frac{8 \tilde K}{x} \Bigg(\int_{\bar\eta_n/2}^{\eta} \epsilon^{-2/(3\gamma)} \text{d}\epsilon + (\delta + 2 \bar\eta_n) \eta^{-4/(3\gamma)} \Bigg)\Bigg]^6 \\
	& \quad + \IP\Big( \sup_{d(a,b) \leq \bar\eta_n} | \hat H^{j_1, j_2}_n(a;\omega) - \hat H^{j_1, j_2}_n(b;\omega) | > x/4 \Big).
\end{split}
\end{equation*}
Now, choosing $2/3 < \gamma < 1$ and letting $n$ tend to infinity, the second term tends to zero by Lemma~\ref{lem:OrderSmallIncrements}, because, by construction,  $1/\gamma > 1$ and
$d(a,b) \leq \bar\eta_n$ if and only if $\|a-b\|_1 \leq 2^{2/\gamma} (n b_n)^{-1/\gamma}$.
All together, this yields
\[\lim_{\delta \downarrow 0} \limsup_{n \rightarrow \infty} \IP\Big( \sup_{d(a,b) \leq \delta} | \hat H_n(a;\omega) - \hat H_n(b;\omega)| > x\Big)
  \leq \Bigg[ \frac{8 \tilde K}{x} \int_{0}^{\eta} \epsilon^{-2/(3\gamma)} \text{d}\epsilon \Bigg]^6,\]
for every $x, \eta > 0$. The claim then follows, as the integral on the right-hand side may be arbitrarily small by choosing $\eta$ accordingly.\hfill\qed

\subsubsection{Proof of Theorem~\ref{thm:AsympDensityEstimator}(ii)}
\label{Partii}

Following the arguments which were applied in Section~8.1 of \cite{kley2014} we can derive asymptotic expansions for $\IE[I_{n,U}^{j_1, j_2}(\omega; \tau_1,\tau_2)]$ and $\IE[\hat G_{n,U}^{j_1, j_2}(\omega; \tau_1,\tau_2)]$. In fact, it is easy to see that the proofs can still be applied when the Laplace cumulants
\[\cum\big( I\{X_{k_1} \leq x_1\}, I\{X_{k_2} \leq x_2\}, \ldots, I\{X_{0} \leq x_p\}\big)\]
which were considered in~\cite{kley2014} are replaced by their multivariate counterparts
\[\cum\big( I\{X_{k_1,j_1} \leq x_1\}, I\{X_{k_2,j_2} \leq x_2\}, \ldots, I\{X_{0,j_p} \leq x_p\}\big).\]

More precisely, we now state Lemma~\ref{thm:522unif} and~\ref{thm:562unif} (without proof) that are multivariate counterparts to Lemmas~8.4 and~8.5 in \cite{kley2014}, for which we assume
\begin{assumption}\label{ass:CS}
Let $p\geq 2, \delta>0$. There exists a non-increasing function $a_p: \IN \to\IR^+$ such that
$\sum_{k \in \IN} k^\delta a_p(k) < \infty$
and
\[
\sup_{x_1,...,x_p} | \cum\big( I\{X_{k_1,j_1} \leq x_1\}, I\{X_{k_2,j_2} \leq x_2\}, \ldots, I\{X_{0,j_p} \leq x_p\}\big) | \leq a_p\big(\max_j |k_j|\big),
\]
for all $j_1, \ldots, j_p = 1,\ldots,d$.
\end{assumption}
Note that Assumption~\ref{ass:CS} follows from condition~\eqref{eq:boundcum}, which is in turn implied by Assumption~\ref{ass:exp_alpha_mix}, but that it is in fact somewhat weaker. We now state the first of the two lemmas. It is a generalisation of Theorem~5.2.2 in \cite{Brillinger1975}.

\begin{lemma}
  \label{thm:522unif}
  Under Assumption~\ref{ass:CS} with $K=2, \delta > 3$,
  \begin{equation}\label{also}
    \IE I_{n, U}^{j_1, j_2}(\omega; \tau_1, \tau_2) =
    \begin{cases}
      \mathfrak{f}^{j_1, j_2}(\omega; \tau_1, \tau_2)
      + \frac{1}{2\pi n} \Big[\frac{\sin(n \omega/2)}{\sin(\omega/2)}\Big]^2 \tau_1 \tau_2
      + \varepsilon_n^{\tau_1, \tau_2}(\omega) & \omega \neq 0 \mod 2\pi \\
      \mathfrak{f}^{j_1, j_2}(\omega; \tau_1, \tau_2)
      + \frac{n}{2\pi} \tau_1 \tau_2
      + \varepsilon_n^{\tau_1, \tau_2}(\omega) & \omega = 0 \mod 2\pi
    \end{cases}
  \end{equation}
  with $\sup_{\tau_1, \tau_2 \in [0,1], \omega \in \IR} | \varepsilon_n^{\tau_1, \tau_2}(\omega) | = O(1/n)$.
\end{lemma}

\noindent The second of the two lemmas is a generalisation of Theorem~5.6.1 in~\cite{Brillinger1975}.

\begin{lemma}\label{thm:562unif}
Assume that Assumption~\ref{ass:CS}, with $p=2$ and $ \delta > k+1$, and Assumption~\ref{ass:W} hold. Then, with the notation of Theorem~\ref{thm:AsympDensityRankEstimator},
\begin{multline*}
\sup_{\tau_1, \tau_2 \in [0,1], \omega \in \IR} \Big| \IE \hat G_n^{j_1, j_2} (\omega; \tau_1, \tau_2) - \mathfrak{f}^{j_1, j_2}(\omega; \tau_1, \tau_2) - \big\{ \bB_n^{(k)}(\omega; \tau_1, \tau_2) \big\}_{j_1, j_2} \Big| \\ = O((n b_n)^{-1}) + o(b_n^k).
\end{multline*}
\end{lemma}

\noindent Because condition~\eqref{eq:boundcum}, which is implied by Assumption~\ref{ass:exp_alpha_mix}, implies Assumption~\ref{ass:CS}, Lemma~\ref{thm:562unif} implies Theorem~\ref{thm:AsympDensityEstimator}(ii).\hfill\qed

\subsubsection{Proof of Theorem~\ref{thm:AsympDensityEstimator}(iii)}
\label{Partiii}

Using~\eqref{prop:inr:eqn:equiv-dftsRandU} and argument similar to the ones in the proof of Lemma~\ref{lem:OrderSmallIncrements} it follows that

\[
\sup_{\omega\in\IR}\sup_{\tau_1, \tau_2 \in [0,1]} | \hat G^{j_1, j_2}_{n,R} (\omega; \tau_1, \tau_2) - \hat G^{j_1, j_2}_{n,U} (\omega; \hat F_{n,j_1}^{-1}(\tau_1), \hat F_{n,j_2}^{-1}(\tau_2)) | = o_p(1).
\]
It therefore suffices to bound the differences
\[
\sup_{\tau_1, \tau_2 \in [0,1]} | \hat G^{j_1, j_2}_{n,U} (\omega; \tau_1, \tau_2) - \hat G^{j_1, j_2}_{n,U} (\omega; \hat F_{n,j_1}^{-1}(\tau_1), \hat F_{n,j_2}^{-1}(\tau_2)) |
\]
for $j_1, j_2 = 1,\ldots,d$, pointwise and uniformly in $\omega$.

We first prove the statement for fixed $\omega \in \IR$ in full details and will later sketch the additional arguments needed for the proof of the uniform result.
For any $x>0$ and sequence $\delta_n$ we have,
\begin{align}
& P^n(\omega) := \nonumber \\
&  \IP \Big( \sup_{\tau_1, \tau_2 \in [0,1]} | \hat G^{j_1, j_2}_{n,U} (\omega; \hat F_{n,j_1}^{-1}(\tau_1), \hat F_{n,j_2}^{-1}(\tau_2)) - \hat G^{j_1, j_2}_{n,U} (\omega; \tau_1, \tau_2) | > x ((nb_n)^{-1/2}+b_n^k) \Big) \nonumber \\
& \leq  \IP \Big( \sup_{\tau_1, \tau_2 \in [0,1]} \sup_{\substack{\|(u,v) - (\tau_1, \tau_2)\|_{\infty} \\ \leq \sup_{i=1,2; \tau \in [0,1]} | \hat F_{n, j_i}^{-1}(\tau) - \tau | }} | \hat G^{j_1, j_2}_{n, U} (\omega; u, v) - \hat G^{j_1, j_2}_{n,U} (\omega; \tau_1, \tau_2) | \nonumber \\
& \hspace*{9cm} > x((nb_n)^{-1/2}+b_n^k) \Big) \nonumber
\end{align}
\begin{align}
  & \leq  \IP \Big( \sup_{\tau_1, \tau_2 \in [0,1]} \sup_{\substack{|u - \tau_1| \leq \delta_n \\ |v - \tau_2| \leq \delta_n}} | \hat G^{j_1, j_2}_{n, U} ( \omega; u, v) - \hat G^{j_1, j_2}_{n,U} (\omega; \tau_1, \tau_2) | > x ((nb_n)^{-1/2}+b_n^k), \nonumber \\
  & \qquad \sup_{i=1,2; \tau \in [0,1]} | \hat F_{n,j_i}^{-1}(\tau) - \tau | \leq \delta_n \Big) + \sum_{i=1}^2 \IP \Big(\sup_{\tau \in [0,1]} | \hat F_{n,j_i}^{-1}(\tau) - \tau | > \delta_n \Big) \nonumber\\
  & = P^n_1 + P^n_2,\quad\text{say}.\nonumber
  \end{align}
We choose $\delta_n$ such that
$n^{-1/2} \ll \delta_n = o( n^{-1/2} b_n^{-1/2} (\log n)^{-D} )$,
where $D$ denotes the constant from Lemma~\ref{lem:LipschitzLaplaceSD}. It then follows from Lemma~\ref{lem:quantproc} that $P^n_2$ is $o(1)$. For $P^n_1$, on the other hand, we have the following bound:
  \begin{align}
  &  \IP \Big( \sup_{\tau_1, \tau_2 \in [0,1]} \sup_{\substack{|u - \tau_1| \leq \delta_n \\ |v - \tau_2| \leq \delta_n}} | \hat H^{j_1, j_2}_{n, U} ( \omega; u, v) - \hat H^{j_1, j_2}_{n,U} (\omega; \tau_1, \tau_2) | > (1+(nb_n)^{1/2}b_n^k)x/2  \Big) \nonumber\\
  & \quad + I \Big\{ \sup_{\tau_1, \tau_2 \in [0,1]} \sup_{\substack{|u - \tau_1| \leq \delta_n \\ |v - \tau_2| \leq \delta_n}} | \IE \hat G^{j_1, j_2}_{n, U} ( \omega; u, v) - \IE \hat G^{j_1, j_2}_{n,U} (\omega; \tau_1, \tau_2) | > ((nb_n)^{-1/2}+b_n^k) x/2 \Big\}.\nonumber
\end{align}
The first term tends to zero because of~\eqref{thm:AsympDensityEstimator:eqn:stochequicont}. The indicator vanishes for $n$ large enough, because we have
\begin{equation*}
\begin{split}
  & \sup_{\tau_1, \tau_2 \in [0,1]} \sup_{\substack{|u - \tau_1| \leq \delta_n \\ |v - \tau_2| \leq \delta_n}} | \IE \hat G^{j_1, j_2}_{n, U} ( \omega; u, v) - \IE \hat G^{j_1, j_2}_{n,U} (\omega; \tau_1, \tau_2) | \\
  \leq& \sup_{\tau_1, \tau_2 \in [0,1]} \sup_{\substack{|u - \tau_1| \leq \delta_n \\ |v - \tau_2| \leq \delta_n}} | \IE \hat G^{j_1, j_2}_{n, U} ( \omega; u, v) - \mathfrak{f}^{j_1, j_2}(\omega; u, v) - \big\{ B_n^{(k)}(\omega; u, v) \big\}_{j_1, j_2}| \\
  & + \sup_{\tau_1, \tau_2 \in [0,1]} \sup_{\substack{|u - \tau_1| \leq \delta_n \\ |v - \tau_2| \leq \delta_n}} | \big\{ B_n^{(k)}(\omega; \tau_1, \tau_2) \big\}_{j_1, j_2} +  \mathfrak{f}^{j_1, j_2}(\omega; \tau_1, \tau_2) - \IE \hat G^{j_1, j_2}_{n,U} (\omega; \tau_1, \tau_2) |\\
  & + \sup_{\tau_1, \tau_2 \in [0,1]} \sup_{\substack{|u - \tau_1| \leq \delta_n \\ |v - \tau_2| \leq \delta_n}} | \mathfrak{f}^{j_1, j_2}(\omega; u, v) + \big\{ B_n^{(k)}(\omega; u, v) \big\}_{j_1, j_2} \\
  & \qquad\qquad\qquad\qquad\qquad\qquad - \mathfrak{f}^{j_1, j_2}(\omega; \tau_1, \tau_2) - \big\{ B_n^{(k)}(\omega; \tau_1, \tau_2) \big\}_{j_1, j_2} |\\
  =& o( n^{-1/2} b_n^{-1/2} + b_n^k) + O(\delta_n(1+|\log \delta_n|)^{D}),
\end{split}
\end{equation*}
where $D$ is still the constant from Lemma~\ref{lem:LipschitzLaplaceSD}. To bound the first two terms we have applied part (ii) of Theorem~\ref{thm:AsympDensityEstimator} and Lemma~\ref{lem:LipschitzLaplaceSD} for the third one. Thus, for any fixed~$\omega$, we have shown $P^n(\omega )=o(1)$, which is the pointwise version of the claim.

Next, we outline the proof of the uniform (with respect to $\omega$) convergence. For any~$y_n > 0$, by similar arguments as above, using the same $\delta_n$, we have
\begin{equation*}
\begin{split}
&\IP \Big( \sup_{\omega\in\IR}\sup_{\tau_1, \tau_2 \in [0,1]} | \hat G^{j_1, j_2}_{n,R} (\omega; \tau_1, \tau_2) - \hat G^{j_1, j_2}_{n,U} (\omega; \tau_1, \tau_2) | > y_n\Big)
\\
&\leq \IP \Big( \sup_{\omega\in\IR}\sup_{\tau_1, \tau_2 \in [0,1]} \sup_{\substack{|u - \tau_1| \leq \delta_n \\ |v - \tau_2| \leq \delta_n}} | \hat H^{j_1, j_2}_{n, U} ( \omega; u, v) - \hat H^{j_1, j_2}_{n,U} (\omega; \tau_1, \tau_2) | > (nb_n)^{1/2}y_n/2  \Big)
\\
& + I \Big\{\sup_{\omega \in \IR} \sup_{\tau_1, \tau_2 \in [0,1]} \sup_{\substack{|u - \tau_1| \leq \delta_n \\ |v - \tau_2| \leq \delta_n}} | \IE \hat G^{j_1, j_2}_{n, U} ( \omega; u, v) - \IE \hat G^{j_1, j_2}_{n,U} (\omega; \tau_1, \tau_2) | > y_n/2 \Big\}
+ o(1).
\end{split}
\end{equation*}
The indicator in the latter expression is $o(1)$ by the same arguments as above  [note that Lemma~\ref{lem:LipschitzLaplaceSD} and the statement of part (ii) both hold uniformly with respect to $\omega \in \IR$]. For the bound of the probability, note that by Lemma~\ref{lem:BoundDFT},
\[\sup_{\tau_1,\tau_2}\sup_{k=1,...,n} |I_{n,U}^{j_1, j_2}(2\pi k/n; \tau_1,\tau_2)| = O_p(n^{2/K}), \text{ for any $K>0$.}\]
Moreover, by the uniform Lipschitz continuity of~$W$ the function $W_n$ is also uniformly Lipschitz continuous with constant of order $O(b_n^{-2})$. Combining those facts with Lemma~\ref{lem:LipschitzLaplaceSD} and the assumptions on~$b_n$, we obtain
\[
\sup_{\substack{\omega_1,\omega_2 \in\IR \\|\omega_1-\omega_2|\leq n^{-3}}} \sup_{\tau_1, \tau_2 \in [0,1]} | \hat H^{j_1, j_2}_{n, U} ( \omega_1; \tau_1,\tau_2) - \hat H^{j_1, j_2}_{n,U} (\omega_2; \tau_1, \tau_2) | = o_p(1).
\]
By the periodicity of $\hat H^{j_1, j_2}_{n, U}$ (with respect to $\omega$), it suffices to show that
\[
\max_{\omega =0,2\pi n^{-3},\ldots ,2\pi}\sup_{\tau_1, \tau_2 \in [0,1]}\sup_{\substack{|u - \tau_1| \leq \delta_n \\ |v - \tau_2| \leq \delta_n}} | \hat H^{j_1, j_2}_{n, U} ( \omega; u, v) - \hat H^{j_1, j_2}_{n,U} (\omega; \tau_1, \tau_2) | = o_p(1).
\]
By Lemmas~\ref{lem:Thm224VW} and~\ref{lem:OrderSmallIncrements} there exists a random variable $S(\omega)$ such that
\begin{eqnarray*}
\sup_{\tau_1, \tau_2 \in [0,1]}\sup_{\substack{|u - \tau_1| \leq \delta_n \\ |v - \tau_2| \leq \delta_n}} | \hat H^{j_1, j_2}_{n, U} ( \omega; u, v) - \hat H^{j_1, j_2}_{n,U} (\omega; \tau_1, \tau_2) |
\leq |S(\omega)| + R_n(\omega),
\end{eqnarray*}
for any fixed  $\omega \in \IR$, with $\sup_{\omega\in\IR} |R_n(\omega)| = o_p(1)$ and
\[
\max_{\omega =0,2\pi n^{-3}...,2\pi}\IE[|S^{2L}(\omega)|]  \leq K_L^{2L} \Bigg(\int_{0}^{\eta} \epsilon^{-4/(2L\gamma)} \text{d}\epsilon + (\delta_n^{\gamma/2} + 2(n b_n)^{-1/2} ) \eta^{-8/(2L\gamma)} \Bigg)^{2L}
\]
for any $0< \gamma < 1, L \in \IN$, $0<\eta<\delta_n$, and  a constant $K_L$ depending on $L$ only.
For appropriately chosen $L$ and $\gamma$, this latter bound is $o(n^{-3})$. Note that the maximum is with respect to a set of cardinality $O(n^3)$, which completes the proof of part (iii).
\hfill\qed

\subsection{Auxiliary Lemmas}
\label{sec:auxlemmas}

In this section we state multivariate versions of the auxiliary lemmas from Section~7.4 in~\cite{kley2014}. Note that Lemma~\ref{lem:Thm224VW} is unaltered and therefore stated without proof. The remaining lemmas are adapted to the multivariate quantities and proofs or directions on how to adapt the proofs in~\cite{kley2014} are collected in the end of this section.

For the statement of Lemma~\ref{lem:Thm224VW}, we define the {Orlicz norm} [see e.g.~\cite{vanderVaartWellner1996}, Chapter 2.2] of a real-valued random variable $Z$ as
\[
\|Z\|_\Psi  = \inf\Big\{C>0: \IE \Psi\Big( {|Z|}/{C}\Big) \leq 1\Big\},
\]
where $\Psi: \IR^+ \to \IR^+$ may be any non-decreasing, convex function with $\Psi(0) = 0$.

For the statement of Lemmas~\ref{lem:SixthMomIncrementHn}, \ref{lem:OrderDftYinA}, and~\ref{lem:BoundDFT} we define, for any Borel set~$A$,
\begin{equation}\label{eqn:djnA}
d^j_n(\omega; A) := \sum_{t=0}^{n-1} I\{X_{t,j} \in A\} \ee^{-\ii t \omega}.
\end{equation}

\begin{lemma}
  \label{lem:Thm224VW}
  Let $\{\IG_t : t \in T\}$ be a separable stochastic process with $\|\IG_s - \IG_t\|_{\Psi} \leq C d(s,t)$ for all $s, t$ with $d(s,t) \geq \bar\eta/2 \geq 0$. Denote by $D(\epsilon,d)$ the packing number of the metric space $(T,d)$. Then, for any $\delta >0 $, $\eta \geq \bar \eta$, there exists a random variable $S_1$ and a constant $K<\infty$ such that
\begin{eqnarray*}
\sup_{d(s,t) \leq \delta} |\IG_s - \IG_t|
&\leq&
S_1 + 2 \sup_{d(s,t) \leq \bar\eta, t \in \tilde T} |\IG_s - \IG_t|\quad\text{and}
\\
\|S_1\|_\Psi &\leq& K \Big[ \int_{\bar\eta/2}^{\eta} \Psi^{-1} \big( D(\epsilon, d)\big) \text{d}\epsilon +  (\delta + 2 \bar\eta) \Psi^{-1} \big( D^2(\eta, d) \big) \Big],
\end{eqnarray*}
where the set $\tilde T$ contains at most $D(\bar\eta,d)$ points. In particular, by Markov's inequality [cf.~\cite{vanderVaartWellner1996}, p.~96],
\begin{multline*}
\IP\Big( |S_1| > x\Big) \leq \Big(\Psi\Big( x \big[ 8 K \big( \int_{\bar\eta/2}^{\eta} \Psi^{-1}\big(D(\epsilon, d)\big) \text{d}\epsilon + (\delta + 2 \bar\eta) \Psi^{-1}\big(D^2(\eta, d) \big) \big) \big]^{-1} \Big) \Big)^{-1} .
\end{multline*}
 for any $x>0$.
\end{lemma}

\begin{lemma} \label{lem:SixthMomIncrementHn}
Let $\bX_{0},...,\bX_{n-1}$, where $\bX_t = (X_{t,1}, \ldots, X_{t,d})$, be the finite realisation of a strictly stationary process with $X_{0,j} \sim U[0,1]$, $j=1,\ldots,d$. Let Assumption~\ref{ass:W} hold. For $x = (x_1, x_2)$ let~$\hat H^{j_1, j_2}_n(x;\omega) := \sqrt{n b_n} (\hat G^{j_1, j_2}_n(x_1,x_2;\omega) - \IE[\hat G^{j_1, j_2}_n(x_1,x_2;\omega)])$. Let $d^j_n(\omega; A)$ be defined as in~\eqref{eqn:djnA}.
Assume that, for $p=1,\ldots,P$, there exist a constant $C$ and a function $g: \IR^+ \to \IR^+$, both independent of $\omega_1,...,\omega_p   \in \IR , n$ and $A_1,...,A_p$, such that
\begin{equation}\label{lem:SixthMomIncrementHn:eqn:Assumption}
\Big| \cum(d_n^{j_1}(\omega_1; A_1), \ldots, d_n^{j_p}(\omega_p; A_p)) \Big| \leq C \Big( \Big| \Delta_n\Big(\sum_{i=1}^p \omega_i \Big) \Big| + 1 \Big) g(\eps)
\end{equation}
for any indices $j_1, \ldots, j_p \in \{1,\ldots,d\}$ and intervals $A_1, \ldots, A_p$ with $\min_{k} \IP(X_{0,j_k} \in A_k) \leq \varepsilon$. Then, there exists a constant $K$ (depending on $C,L,g$ only) such that
\[
\sup_{\omega\in \IR}\sup_{\| a - b \|_1 \leq \varepsilon} \IE | \hat H^{j_1, j_2}_n(a;\omega) - \hat H^{j_1, j_2}_n(b;\omega) |^{2L}
\leq K \sum_{\ell=0}^{L-1} \frac{g^{L-\ell}(\eps)}{(n b_n)^\ell}
\]
for all $\varepsilon$ with $g(\varepsilon)<1$ and all $L = 1, \ldots, P$.
\end{lemma}

\begin{lemma} \label{lem:LipschitzLaplaceSD}
Under the assumptions of Theorem~\ref{thm:AsympDensityRankEstimator}, the derivative
\[\displaystyle{(\tau_1, \tau_2) \mapsto \frac{{\rm d}^k}{{\rm d}\omega^k}\mathfrak{f}^{j_1, j_2}(\omega; \tau_1, \tau_2)}\]
exists and satisfies, for any $k \in \IN_0$ and some constants $C,d$ that are independent of $a=(a_1,a_2),b=(b_1,b_2)$, but may depend on $k$,
\[
\sup_{\omega \in \IR } \Big|\frac{{\rm d}^k}{{\rm d}\omega^k}\mathfrak{f}^{j_1, j_2}(\omega; a_1, a_2)-\frac{{\rm d}^k}{{\rm d}\omega^k} \mathfrak{f}^{j_1, j_2}(\omega; b_1, b_2)\Big| \leq C \|a-b\|_1 (1+|\log\|a-b\|_1|)^D.
\]
\end{lemma}

\begin{lemma} \label{lem:OrderDftYinA}
Let the strictly stationary process $(\bX_t)_{t \in \IZ}$ satisfy condition~\eqref{eq:boundcum}. Let $d^j_n(\omega; A)$ be defined as in~\eqref{eqn:djnA}.
Let $A_1, \ldots, A_p \subset [0,1]$ be intervals, and let
\[\eps := \min_{k=1,...,p} \IP(X_{0,j_k} \in A_k).\] Then, for any $p$-tuple  $\omega_1,...,\omega_p  \in \IR$ and $j_1, \ldots, j_p \in \{1,\ldots,d\}$,
\[
\Big| \cum(d_n^{j_1}(\omega_1; A_1), \ldots, d_n^{j_p}(\omega_p; A_p)) \Big| \leq C \Big( \Big| \Delta_n\Big(\sum_{i=1}^p \omega_i \Big) \Big| + 1 \Big) \varepsilon (|\log \eps|+1)^D,
\]
where $\Delta_n(\lambda) := \sum_{t=0}^{n-1}e^{\ii t\lambda}$ and the constants $C,D$ depend only on $K,p$, and~$\rho$ [with $\rho$ from condition~\eqref{eq:boundcum}].
\end{lemma}

\begin{lemma}\label{lem:BoundCdf}
Let the strictly stationary process $(\bX_t)_{t \in \IZ}$ satisfy condition~\eqref{eq:boundcum} and $X_{0,j} \sim U[0,1]$. Denote the empirical distribution function of $X_{0,j},...,X_{n-1,j}$ by $\hat F_{n,j}$. Then, for any $k \in \IN$, there exists a constant $d_k$ depending only on~$k$, such that
\begin{multline*}
\sup_{x,y\in [0,1], |x-y|\leq \delta_n} \sqrt{n}|\hat F_{n,j}(x) - \hat F_{n,j}(y) - (x-y)|
\\
= O_p\Big((n^2\delta_n+n)^{1/2k}(\delta_n|\log\delta_n|^{d_k}+n^{-1})^{1/2}\Big),
\end{multline*}
as $\delta_n \to 0$.
\end{lemma}

\begin{lemma}\label{lem:quantproc}
Let $\bX_{0},...,\bX_{n-1}$, where $\bX_t = (X_{t,1}, \ldots, X_{t,d})$, be the finite realisation of a strictly stationary process satisfying condition~\eqref{eq:boundcum} and $X_{0,j} \sim U[0,1]$, $j=1,\ldots,d$. Then,
\[
\sup_{j=1,\ldots,d} \sup_{\tau \in [0,1]} |\hat F_{n,j}^{-1}(\tau) - \tau| = O_p(n^{-1/2}).
\]
\end{lemma}

\begin{lemma}\label{lem:BoundDFT}
Let the strictly stationary process $(\bX_t)_{t \in \IZ}$ satisfy condition~\eqref{eq:boundcum} and $X_{0,j} \sim U[0,1]$. Let $d^j_n(\omega; A)$ be defined as in~\eqref{eqn:djnA}. Then, for any $k\in\IN$,
\[
\sup_{j=1,\ldots,d} \sup_{\omega \in \cF_n} \sup_{y \in [0,1]} | d_n^{j}(\omega; [0,y]) | = O_p(n^{1/2+1/k}).
\]
\end{lemma}

\begin{lemma}\label{lem:OrderSmallIncrements}
Under the assumptions of Theorem~\ref{thm:AsympDensityEstimator}, let $\delta_n$ be a sequence of non-negative real numbers. Assume that there exists $\gamma \in (0,1)$, such that~$\delta_n =O( (n b_n)^{-1/\gamma})$. Then,
\[
\sup_{j_1, j_2, \in \{1,\ldots,d\}} \sup_{\omega \in \IR}\sup_{\substack{u,v \in [0,1]^2 \\ \| u - v \|_1 \leq \delta_n }} |\hat H^{j_1, j_2}_n(u;\omega) - \hat H^{j_1, j_2}_n(v;\omega) | = o_p(1).
\]
\end{lemma}

\smallskip\noindent\textbf{Proof of Lemma~\ref{lem:Thm224VW}.} The lemma is stated unaltered as in~\cite{kley2014}. The proof can be found in Section~8.3.1 of the Online Appendix of~\cite{kley2014}.

\smallskip\noindent\textbf{Proof of Lemma~\ref{lem:SixthMomIncrementHn}.} 
Along the same lines of the proof of the univariate version (Section~8.3.2 in~\cite{kley2014}) we can proof
\begin{equation}
\label{lem:SixthMomIncrementHn:eqn:Repr}
  \IE | \hat H^{j_1, j_2}_n(a;\omega) - \hat H^{j_1, j_2}_n(b;\omega) |^{2L}
  = \sum_{\substack{\{\nu_1, \ldots, \nu_R\} \\ |\nu_j| \geq 2, \ j=1,\ldots,R}}
    \prod_{r=1}^R \mathcal{D}_{a,b}(\nu_r)
\end{equation}
with the summation running over all partitions $\{\nu_1, \ldots, \nu_R\}$ of $\{1,\ldots,2L\}$ such that each set $\nu_j$   contains at least two elements, and
\begin{equation*}
\begin{split}
  \mathcal{D}_{a,b}(\xi) & := \sum_{\ell_{\xi_1}, \ldots, \ell_{\xi_{q}} \in \{1,2\}} n^{-3 q/2} b_n^{q/2} \Big(\prod_{m \in \xi} \sigma_{\ell_m} \Big) \\
  & \qquad \times \sum_{s_{\xi_1}, \ldots, s_{\xi_{q}}=1}^{n-1}
    \Big( \prod_{m \in \xi} W_n(\omega - 2\pi s_{m} / n) \Big)
  \cum( D_{\ell_m, (-1)^{m-1} s_m } \, : \, m \in \xi),
\end{split}
\end{equation*}
for any set $\xi := \{\xi_1, \ldots, \xi_q\} \subset \{1,\ldots,2L\}$, $q:=|\xi|$, and
\[D_{\ell,s} := d_n^{j_1}(2\pi s/{n}; M_1(\ell)) d_n^{j_2}(-2\pi s/{n}; M_2(\ell)), \quad
\ell=1,2, \quad s=1,\ldots,n-1,\]
with the sets $M_1(1)$, $M_2(2)$, $M_2(1)$, $M_1(2)$ and the signs $\sigma_{\ell} \in \{-1,1\}$ defined as
\begin{align}
  \sigma_1 & := 2 I\{a_1 > b_1\} - 1,   & \sigma_2 & := 2 I\{a_2 > b_2\} - 1, \nonumber \\
  M_1(1) & := (a_1 \wedge b_1, a_1 \vee b_1], & M_2(2) & := (a_2 \wedge b_2, a_2 \vee b_2], \label{lem:SixthMomIncrementHn:eqn:Sets} \\
  M_2(1) & := \begin{cases}
              [0, a_2] & b_2 \geq a_2 \\
              [0, b_2] & a_2 > b_2,
           \end{cases} &
  M_1(2) & := \begin{cases}
              [0, b_1] & b_2 \geq a_2 \\
              [0, a_1] & a_2 > b_2.
           \end{cases} \nonumber
\end{align}

Employing assumption~\eqref{lem:SixthMomIncrementHn:eqn:Assumption}, we can further prove, by following the arguments of the univariate version, that
\begin{equation*}
  \label{lem:SixthMomIncrementHn:eqn:RateDabxi}
  \sup_{\substack{\xi \subset \{1,\ldots,2L\} \\ |\xi| = q}} \sup_{\|a-b\|_1 \leq \varepsilon} | \mathcal{D}_{a,b}(\xi) | \leq C (n b_n)^{1-q / 2} g(\eps), \quad  2\leq q \leq 2L.
\end{equation*}
The lemma then follows, by observing that
\[
 \Big|  \prod_{r=1}^R \mathcal{D}_{a,b}(\nu_r) \Big| \leq Cg^R(\eps)(nb_n)^{R-L}
\]
for any partition in~\eqref{lem:SixthMomIncrementHn:eqn:Repr} [note that $\sum_{r=1}^R |\nu_r| = 2L$].\hfill\qed

\smallskip\noindent\textbf{Proof of Lemma~\ref{lem:LipschitzLaplaceSD}.}
Note that
\begin{equation*}
\begin{split}
& \cum(I\{X_{0,j_1} \leq q_{j_1}(a_1) \},I\{X_{k,j_2}\leq q_{j_2}(a_2)\}) \\
& \quad - \cum(I\{X_{0,j_1} \leq q_{j_1}(b_1)\},I\{X_{k,j_2}\leq q_{j_2}(b_2)\}) \\
&  = \sigma_1\cum(I\{F_{j_1}(X_{0,j_1}) \in M_1(1)\},I\{F_{j_2}(X_{k,j_2}) \in M_2(1)\}) \\
& \qquad\qquad + \sigma_2\cum(I\{F_{j_1}(X_{0,j_1}) \in M_1(2)\},I\{F_{j_2}(X_{k,j_2}) \in M_2(2)\}),
\end{split}
\end{equation*}
with the sets $M_1(1)$, $M_2(2)$, $M_2(1)$, $M_1(2)$ and the signs $\sigma_{\ell} \in \{-1,1\}$ defined in~\eqref{lem:SixthMomIncrementHn:eqn:Sets}.

\noindent From the fact that $\lambda(M_j(j)) \leq \|a-b\|_1$ for $j=1,2$, we conclude that
\begin{equation*}
\begin{split}
& \Big|\frac{{\rm d}^{\ell}}{{\rm d}\omega^{\ell}}\mathfrak{f}^{j_1, j_2}(\omega; a_1, a_2)-\frac{{\rm d}^{\ell}}{{\rm d}\omega^{\ell}} \mathfrak{f}^{j_1, j_2}(\omega; b_1, b_2)\Big| \\
& \qquad \leq \sum_{k\in\IZ} |k|^{\ell}|\cum(I\{F_{j_1}(X_{0,j_1}) \in M_1(1)\},I\{F_{j_2}(X_{k,j_2}) \in M_2(1)\})| \\
& \qquad\qquad + \sum_{k\in\IZ} |k|^{\ell} |\cum(I\{F_{j_1}(X_{0,j_1}) \in M_1(2)\},I\{F_{j_2}(X_{k,j_2}) \in M_2(2)\})| \\
& \qquad \leq 4 \sum_{k=0}^\infty k^{\ell} \Big((K \rho^{\ell}) \wedge \|a-b\|_1\Big).
\end{split}
\end{equation*}
The assertion then follows by after some algebraic manipulations.\hfill\qed

\smallskip\noindent\textbf{Proof of Lemma~\ref{lem:OrderDftYinA}.} Similar to (8.27) in \cite{kley2014} we have, by the definition of cumulants and strict stationarity, 
\begin{multline}\label{A}
\cum(d_n^{j_1}(\omega_1; A_1), \ldots, d_n^{j_p}(\omega_p; A_p)) \\
= \sum_{u_2,...,u_p = -n}^n \cum(I{\{X_{0,j_1} \in A_1\}}, I{\{X_{u_2,j_2} \in A_2\}} \ldots, I{\{X_{u_p,j_p} \in A_p\}}) \exp\Big(-\ii \sum_{j=2}^p \omega_j u_j \Big) \\
 \times\sum_{t_1 = 0}^{n-1} \exp\Big(-\ii t_1 \sum_{j=1}^p \omega_j\Big) I_{\{0 \leq t_1+u_2 < n\}} \cdots I_{\{0 \leq t_1 + u_p < n\}}.
\end{multline}
By Lemma~8.1 in \cite{kley2014},
\begin{multline}\label{B}
\Big|\Delta_n(\sum_{j=1}^p \omega_j) - \sum_{t_1 = 0}^{n-1} \exp\Big(-\ii t_1 \sum_{j=1}^p \omega_j\Big) I{\{0 \leq t_1+u_2 < n\}} \cdots I{\{0 \leq t_1 + u_p < n\}} \Big| \\
\leq 2\sum_{j=2}^p |u_j|.
\end{multline}
Following the arguments for the proof of (8.29) in~\cite{kley2014}, we further have, for any $p+1$ intervals $A_0, \ldots, A_p \subset \IR$, any indices $j_0, \ldots, j_p \in \{1,\ldots,d\}$, and any $p$-tuple $\kappa := (\kappa_1,...,\kappa_p) \in \IR_+^p$, $p \geq 2$, that
\begin{multline}\label{elle}
	\sum_{k_1, \ldots, k_{p} = -\infty}^{\infty}
	\Big(1 + \sum_{\ell=1}^{p} |k_\ell|^{\kappa_\ell} \Big) \big| \cum\big( I {\{X_{k_1, j_1} \in A_1\}}, \ldots, I {\{ X_{k_{p}, j_p} \in A_{p} \}}, I {\{ X_{0,j_0} \in A_{0}\}} \big) \big| \\
	\leq C \eps (|\log \eps|+1)^d.
\end{multline}
To this end, define $k_0 = 0$, consider the set
\[T_m := \big\{(k_1,...,k_p) \in \IZ^p| \max_{i,j = 0,...,p} |k_i-k_j| = m\big\},\]
 and note that $|T_m| \leq c_p m^{p-1}$  for some  constant~$c_p$. From the definition of cumulants and some simple algebra we get the bound
\[|\cum(I\{X_{t_1, j_1} \in A_1\},...,I\{X_{t_p, j_p} \in A_p\})| \leq C\min_{i=1,...,p} P(X_{0, j_i} \in A_i).\]
With this bound and condition~\eqref{eq:boundcum}, which is implied by Assumption~\ref{ass:exp_alpha_mix}, we obtain, employing the above notation, that
\begin{multline*}
\sum_{k_1, \ldots, k_p = -\infty}^{\infty}
\Big(1 + \sum_{j=1}^{p} |k_{\ell}|^{\kappa_{\ell}} \Big) \big| \cum\big( I {\{X_{k_1, j_1} \in A_1\}}, \ldots, I {\{ X_{k_{p}, j_p} \in A_{p} \}}, I_{\{ X_{0, j_0} \in A_0\}} \big) \big| \\
=  \sum_{m=0}^\infty \sum_{(k_1,...,k_p) \in T_m} \Big(1 + \sum_{\ell=1}^{p} |k_{\ell}|^{\kappa_{\ell}} \Big) \big| \cum\big( I {\{X_{k_1,j_1} \in A_1\}}, \ldots, I {\{ X_{k_{p}, j_p} \in A_{p} \}}, I {\{ X_{0, j_0} \in A_0\}} \big) \big| \\
\leq \sum_{m=0}^\infty \sum_{(k_1,...,k_p) \in T_m} \Big(1 + p m^{\max_j \kappa_j} \Big) \Big( \rho^{m} \wedge \eps \Big) K_p
\leq  C_p \sum_{m=0}^\infty \Big( \rho^{m} \wedge \eps \Big)  |T_m|m^{\max_j \kappa_j}.
\end{multline*}

\noindent For $\eps \geq \rho$, \eqref{elle} then follows trivially. For $\eps<\rho$, set $m_\eps := \log \eps/ \log \rho$ and note that $\rho^{m} \leq \eps$ if and only if $m \geq m_\eps$. Thus,
\begin{equation*}
\sum_{m=0}^\infty \Big( \rho^{m} \wedge \eps \Big) m^u \leq \sum_{m \leq m_\eps} m^u \eps + \sum_{m > m_\eps} m^u \rho^m \leq C\Big(\eps m_\eps^{u+1} + \rho^{m_\eps} \sum_{m=0}^\infty (m + m_\eps)^u \rho^m \Big).
\end{equation*}
The fact that $\rho^{m_\eps}=\eps$ completes the proof   of the desired inequality~(\ref{elle}). The assertion follows from  \eqref{A}, \eqref{B}, \eqref{elle} and the triangle inequality.
\hfill\qed

\smallskip\noindent\textbf{Proofs of Lemmas~\ref{lem:BoundCdf}, \ref{lem:quantproc} and~\ref{lem:BoundDFT}.}
Note that the component processes $(X_{t,j})$ are stationary and fulfill Assumption~(C) in~\cite{kley2014}, for every $j = 1, \ldots, d$. The assertion then follow from the univariate versions (i.\,e., Lemma~8.6, 7.5 and 7.6 in \cite{kley2014}, respectively), as the dimension $d$ does not depend on $n$.\hfill\qed

\smallskip\noindent\textbf{Proof of Lemma~\ref{lem:OrderSmallIncrements}.}
Assume, without loss of generality, that $n^{-1} = o(\delta_n)$ [otherwise, enlarge the supremum by considering $\tilde \delta_n := \max(n^{-1},\delta_n)$]. With the notation $a = (a_1, a_2)$ and $b = (b_1, b_2)$, we have
\begin{equation*}
  \hat H^{j_1, j_2}_n(a;\omega) - \hat H^{j_1, j_2}_n(b;\omega)
  = b_n^{1/2} n^{-1/2} \sum_{s=1}^{n-1} W_n(\omega - 2\pi s / n) (K_{s,n}(u,v) - \IE K_{s,n}(u,v) )
\end{equation*}
where, with $d^j_{n,U}$ defined  in \eqref{eq:defdnu},
\begin{equation*}
\begin{split}
  K_{s,n}(a,b) & := n^{-1} \big( d_{n,U}^{j_1} (2\pi s/n; u_1) d_{n,U}^{j_2}(-2\pi s/n; u_2) - d_{n,U}^{j_1} (2\pi s/n; v_1) d_{n,U}^{j_2}(-2\pi s/n; v_2) \big) \\
    & = d_{n,U}^{j_1} (2\pi s/n; u_1) n^{-1} \big[ d_{n,U}^{j_2}(-2\pi s/n; u_2) - d_{n,U}^{j_2}(-2\pi s/n; v_2) \big] \\
    & \quad + d_{n,U}^{j_2}(-2\pi s/n; v_2) n^{-1} \big[ d_{n,U}^{j_1} (2\pi s/n; u_1) - d_{n,U}^{j_1} (2\pi s/n; v_1) \big].
\end{split}
\end{equation*}
By Lemma~\ref{lem:BoundDFT}, we have, for any $k\in \IN$,
\begin{equation} \label{BoundDFT}
  \sup_{y\in[0,1]} \sup_{\omega \in \mathcal{F}_n} | d_{n,U}^j(\omega; y) | = O_p\Big(n^{1/2+1/k} \Big).
\end{equation}
Employing Lemma~\ref{lem:BoundCdf}, we have, for any $\ell\in\IN$ and $j=1,\ldots,d$,
\begin{equation*}
\begin{split}
  & \sup_{\omega \in \IR} \sup_{y\in[0,1]} \sup_{x : |x-y| \leq \delta_n}
  n^{-1} | d_{n,U}^{j}(\omega; x) - d_{n,U}^{j}(\omega;y) | \\
  &\leq \sup_{y\in[0,1]} \sup_{x : |x-y| \leq \delta_n}
  n^{-1} \sum_{t=0}^{n-1} | I\{F_j(X_{t,j}) \leq x\} - I\{F_j(X_{t,j}) \leq y\} | \\
  & \leq \sup_{y\in[0,1]} \sup_{x : |x-y| \leq \delta_n}
    |\hat F_{n,j}(x \vee y) - \hat F_{n,j}(x \wedge y) - x \vee y + x \wedge y |
    + C \delta_n
  \\
  &= O_p\big(\rho_n(\delta_n,\ell) + \delta_n\big),
\end{split}
\end{equation*}
with $\rho_n(\delta_n, \ell) := n^{-1/2}(n^2 \delta_n+n)^{1/2\ell}(\delta_n|\log \delta_n|^{D_\ell}+n^{-1})^{1/2}$, $\hat F_{n,j}$ denoting the empirical distribution function of $F_j(X_{0,j}), \ldots, F_j(X_{n-1,j})$, and $d_\ell$ being a constant depending only on~$\ell$.
Combining these arguments and observing that
\begin{equation}\label{eqn:propWn}
	\sup_{\omega\in\IR}\sum_{s=1}^{n-1} \Big| W_n(\omega - 2\pi s/n)\Big| = O(n)
\end{equation}
yields
\begin{equation} \label{term1}
\sup_{\omega \in \IR} \!\!\sup_{\substack{u,v \in [0,1]^2 \\ \| u - v \|_1 \leq \delta_n }} \Big|\sum_{s=1}^{n-1} W_n(\omega - 2\pi s / n) K_{s,n}(u,v) \Big|
= O_p\big( n^{3/2+1/k}(\rho(\delta_n,\ell) + \delta_n) \big) .
\end{equation}
With $M_i(j)$, $i,j=1,2$, as defined in~\eqref{lem:SixthMomIncrementHn:eqn:Sets}, we have
\begin{equation}\label{flayosc}
\begin{split}
& \sup_{\|a-b\|_1 \leq \delta_n}  \sup_{s=1,\ldots,n-1}  | \IE K_{s,n}(a,b) | \\
& \qquad \leq n^{-1}\!\! \sup_{\|a-b\|_1 \leq \delta_n} \sup_{s=1,\ldots,n-1} \big|
			\cum(d_{n,U}^{j_1}(2\pi s/{n}; M_1(1)), d_{n,U}^{j_2}(-2\pi s/{n}; M_2(1)))\big| \\
& \quad\quad + n^{-1}\!\! \sup_{\|a-b\|_1 \leq \delta_n} \sup_{s=1,\ldots,n-1} \big|\cum(d_{n,U}^{j_1}(2\pi s/{n}; M_1(2)), d_{n,U}^{j_2}(-2\pi s/{n}; M_2(2))) \big|
\end{split}
\end{equation}
where we have used $\IE d_{n,U}^{j}({2\pi} s/n; M) = 0$.
Lemma~\ref{lem:OrderDftYinA} and $\lambda(M_{j}(j)) \leq \delta_n$, for~$j=1,2$ (with $\lambda$ denoting the Lebesgue measure over $\IR$) yield
\begin{multline*}
\sup_{\|a-b\|_1 \leq \delta_n} \sup_{s=1,\ldots,n-1} | \cum(d_n^{j_1}(2\pi s/{n}; M_1(j)), d_n^{j_2}(-2\pi s/{n}; M_2(j))) | \\ 
\leq C (n+1) \delta_n (1+|\log\delta_n|)^D,
\end{multline*}
It follows that the right-hand side in (\ref{flayosc})  is $O(\delta_n|\log \delta_n|^D)$.
Therefore, by~\eqref{eqn:propWn}, we obtain
\[
  \sup_{\omega\in\IR}\sup_{\|a-b\|_1 \leq \delta_n} \Big| b_n^{1/2} n^{-1/2} \sum_{s=1}^{n-1} W_n(\omega - 2\pi s / n) \IE K_{s,n}(a,b) \Big|
  = O\big((n b_n)^{1/2} \delta_n |\log n|^D\big).
\]
In view of the assumption  that  $n^{-1} = o(\delta_n)$, we have $\delta_n = O(n^{1/2}\rho_n(\delta_n,\ell))$, which, in combination with \eqref{term1}, yields
\begin{equation*}
\begin{split}
  & \sup_{\omega\in\IR}\sup_{\| a - b \|_1 \leq \delta_n} |\hat H^{j_1, j_2}_n(a;\omega) - \hat H^{j_1, j_2}_n(b;\omega) |
\\
  &= O_p\Big( (n b_n)^{1/2} [ n^{1/2+1/k} (\rho_n(\delta_n,\ell) + \delta_n) + \delta_n |\log \delta_n|^D] \Big)\\
  & = O_p\Big ( (n b_n)^{1/2} n^{1/2+1/k} \rho_n(\delta_n,\ell)\Big)
\\	
	&= O_p\Big ( (n b_n)^{1/2} n^{1/k+1/\ell} (n^{-1} \vee \delta_n (\log n)^{D_\ell})^{1/2}\Big)
	= o_p(1).
\end{split}
\end{equation*}
The $o_p(1)$ holds, as we have, for arbitrary $k$ and $\ell$,
\[
O( (n b_n)^{1/2}n^{1/k+1/\ell} \delta_n^{1/2}(\log n)^{D_\ell/2 } ) = O( (n b_n)^{1/2-1/2\gamma}n^{1/k+1/\ell} (\log n)^{D_\ell /2} ).
\]
The assumptions on $b_n$ imply   $(n b_n)^{1/2-1/2\gamma} = o(n^{-\kappa})$ for some~$\kappa>0$, such that this latter quantity is $o(1)$ for $k,\ell$ sufficiently large. The term~$(n b_n)^{1/2}n^{1/k+1/\ell}n^{-1/2}$ is handled in a similar fashion. This concludes the proof. \hfill\qed

\bibliography{quantile_cross}